\documentclass[12pt]{amsart}

\usepackage{amssymb,amsthm,amsmath,amsxtra}
\usepackage{fullpage}
\usepackage{comment}
\usepackage{graphicx}
\usepackage{hyperref}  
\usepackage{color}
\usepackage{array}
\usepackage{lscape}

\usepackage{placeins} 

\usepackage{chngcntr}
\usepackage{booktabs}

\numberwithin{equation}{section}

\theoremstyle{plain}

\newtheorem{prop}[equation]{Proposition}
\newtheorem{lem}[equation]{Lemma}

\newtheorem*{cor*}{Corollary}
\newtheorem*{prob*}{Problem}
\newtheorem*{thm*}{Theorem}
\newtheorem*{thma*}{Theorem A}
\newtheorem*{thmb*}{Theorem B}
\newtheorem*{thmc*}{Theorem C}

\theoremstyle{definition}

\newtheorem{alg}[equation]{Algorithm}
\newtheorem{exm}[equation]{Example}

\theoremstyle{remark}
\newtheorem{rmk}[equation]{Remark}

\newenvironment{enumalgalph}
{\begin{enumerate}}
{\end{enumerate}}

\newenvironment{enumalg}
{\begin{enumerate}}
{\end{enumerate}}

\newcommand{\defi}[1]{\textsf{#1}} 				

\DeclareMathOperator{\Aut}{Aut}

\DeclareMathOperator{\ddiv}{div}
\DeclareMathOperator{\disc}{disc}

\DeclareMathOperator{\Gal}{Gal}
\DeclareMathOperator{\impart}{Im}

\DeclareMathOperator{\PGL}{PGL}
\DeclareMathOperator{\PSL}{PSL}
\DeclareMathOperator{\Real}{Re}
\DeclareMathOperator{\repart}{Re}
\DeclareMathOperator{\sgn}{sgn}

\DeclareMathOperator{\Stab}{Stab}
\DeclareMathOperator{\SL}{SL}
\DeclareMathOperator{\GL}{GL}

\DeclareMathOperator{\PSU}{PSU}

\newcommand{\C}{\mathbb C}
\newcommand{\F}{\mathbb F}
\newcommand{\HH}{\mathbb H}
\newcommand{\PP}{\mathbb P}
\newcommand{\Q}{\mathbb Q}
\newcommand{\R}{\mathbb R}
\newcommand{\Z}{\mathbb Z}
\newcommand{\Qbar}{\overline{\mathbb Q}}

\newcommand{\newt}{\mu}

\newcommand{\Belyi}{Bely\u{\i}}

\newcommand{\calD}{\mathcal{D}}

\newcommand{\calH}{\mathcal{H}}

\newcommand{\la}{\langle}
\newcommand{\ra}{\rangle}

\begin{document}

\title{Numerical calculation of three-point branched covers of the projective line}

\author{Michael Klug, Michael Musty, Sam Schiavone}
\address{Department of Mathematics and Statistics, University of Vermont, 16 Colchester Ave, Burlington, VT 05401, USA}
\email{mklug@uvm.edu, michaelmusty@gmail.com, samuel.schiavone@uvm.edu}
\author{John Voight}
\address{Department of Mathematics and Statistics, University of Vermont, 16 Colchester Ave, Burlington, VT 05401, USA; Department of Mathematics, Dartmouth College, 6188 Kemeny Hall, Hanover, NH 03755, USA}
\email{jvoight@gmail.com}

\date{\today}

\begin{abstract}
We exhibit a numerical method to compute three-point branched covers of the complex projective line.  We develop algorithms for working explicitly with Fuchsian triangle groups and their finite index subgroups, and we use these algorithms to compute power series expansions of modular forms on these groups.
\end{abstract}

\maketitle 
\tableofcontents

In his celebrated work \emph{Esquisse d'un Programme} \cite{Grothendieck}, Grothendieck describes an action of the Galois group $\Gal(\Qbar/\Q)$ on the set of dessins d'enfants, encoding in a combinatorial and topological way the action of $\Gal(\Qbar/\Q)$ on the set of three-point branched covers of the projective line.  In trying to understand this action in an explicit way, we are immediately led to study examples; the subject is considerably enriched by the beauty of specific cases and by the calculation of the smallest of them by hand.  Quickly, though, one runs out of examples that can be determined this way.  Grothendieck asks \cite[p.248]{Grothendieck}: 
\begin{quote}
Are such and such oriented maps `conjugate' or: exactly which are the conjugates of a given oriented map? (Visibly, there is only a
finite number of these.)  I considered some concrete cases (for coverings of low degree) by various methods, J.\ Malgoire considered some others---I doubt that there is a uniform method for solving the problem by computer. 
\end{quote}
Consequently, researchers coming from a wide range of perspectives have put forth significant effort to compute more complicated covers of various types; for a survey, we refer the reader to work of Sijsling--Voight \cite{SijslingVoight}.  

In this article, we provide a general-purpose numerical method for the computation of (connected) three-point branched covers of the projective line.  Although we have not yet established rigorously that this recipe yields an algorithm, nor estimated its running time (though see Remark \ref{rmk:makerigorous}), we can verify that the results of our computations are correct.  Herein, we exhibit the method and present some computed examples of high complexity.  

A \defi{\Belyi\ map} is a map of connected Riemann surfaces $\phi:X \to \PP^1$ that is unramified away from $\{0,1,\infty\}$.  By the theorem of \Belyi\ \cite{Belyi,Belyi2}, a Riemann surface $X$ admits a \Belyi\ map if and only if $X$ can be defined over $\Qbar$ as an algebraic curve.  

The input to our method comes from the combinatorial description of \Belyi\ maps: there is a bijection between \Belyi\ maps of degree $d \geq 1$ up to isomorphism and \defi{transitive permutation triples} of degree $d$,\begin{center} 
$\sigma = (\sigma_0, \sigma_1, \sigma_{\infty}) \in S_d^3$  such that 
$\sigma_{\infty} \sigma_1 \sigma_0 = 1$ and \\
$\langle \sigma_0,\sigma_1, \sigma_{\infty} \rangle \leq S_d$ is a transitive subgroup
\end{center}
up to simultaneous conjugation in $S_d$.  If $\sigma$ corresponds to $\phi$, we say that $\phi$ has \defi{monodromy $\sigma$}.  The \defi{dessin} (or \defi{dessin d'enfant}) corresponding to a \Belyi\ map $\phi:X \to \PP^1$ is the inverse image $\phi^{-1}([0,1])$, a bicolored graph embedded on the surface $X$, with vertices above $0$ and above $1$ labelled white and black, respectively.  


Our first main result is to exhibit a numerical method for drawing dessins conformally (compactly) in the geometry $H$, where $H$ is the sphere, Euclidean plane, or hyperbolic plane: in other words, we are concerned with the specific geometric form of the dessin. We do so by working with the corresponding triangle group $\Delta$, determining explicitly the associated finite index subgroup $\Gamma$, and then drawing the dessin on $H$ together with the gluing relations which define the quotient $X=\Gamma \backslash H$.  

Our second main result is to use this explicit description of the Riemann surface $X$ to compute an equation for the \Belyi\ map $f$ numerically.  The main algorithmic tool for this purpose is the method of power series expansions due to Voight--Willis \cite{VoightWillis} and developed further by the first author \cite{Klug}; the original idea goes back to Hejhal \cite{Hejhal} and Stark.  An important feature of this method is that it works directly on the desired Riemann surface, and hence there are no ``parasitic'' solutions to discard as in other methods \cite{Kreines03,Kreines08}.  

The idea to compute covers using modular forms can also be found in the work of Selander and Str\"ombergsson \cite{SelanderStrombergsson}; our approach is `cocompact' and thus parallel to their `noncocompact' method using $q$-expansions on the modular group.  Indeed, these authors use $\Gamma(2)$, the free group on $2$ generators, also known in this context as the triangle group $\Delta(\infty,\infty,\infty)$, so in all cases, triangle groups figure prominently.  It will be interesting to compare the relative advantages of these two approaches; one advantage of our approach is that we do not have to compute $q$-expansions at each of the cusps.

The paper is organized as follows.  In section \ref{sec:background}, we quickly review the relevant background material on \Belyi\ maps.  In section \ref{sec:embedtriangle}, we realize triangle groups explicitly as Fuchsian groups.  In section \ref{sec:cosets}, we exhibit an algorithm that takes a permutation triple $\sigma$ and realizes the cosets of the point stabilizer subgroup $\Gamma \leq \Delta(a,b,c)$; this algorithm also yields a reduction algorithm and a fundamental domain for the action of $\Gamma$ on the upper half-plane $\calH$ (the generic geometry $H$). We then show how these fit together in our first main result, to conformally draw dessins.  Next, in section \ref{sec:powser}, we discuss power series expansions of differentials on Riemann $2$-orbifolds and recall a method for computing these expansions.  In section \ref{sec:examples}, we compute equations arising from these expansions, refined using Newton's method; along the way, we present several complete examples.

As an illustration of our method, we compute the degree 50 rational \Belyi\ map $f(x)=p(x)/q(x)$ where
\begin{align*}
p(x) &= 2^6 (x^4 + 11x^3 - 29x^2 + 11x + 1)^5 (64x^5 - 100x^4 + 150x^3 - 25x^2 + 5x + 1)^5 \\
&\qquad \cdot (196x^5 - 430x^4 + 485x^3 - 235x^2 + 30x + 4) \\
q(x) &= 5^{10} x^7 (x+1)^7 (2x^2 - 3x + 2)^7 (8x^3 - 32x^2 + 10x + 1)^7 
\end{align*}
realizing the group $\PSU_3(\F_5):2$ regularly over $\Q$ (Example \ref{exm:PSU3(5)}).  Such a polynomial was known to exist by the general theory \cite{MalleMatzat}, but the explicit polynomial is new.  We expect that other covers of large degree and arithmetic complexity can be similarly computed, greatly extending the set of examples known, and consequently that our method can be used to understand more comprehensively the action of $\Gal(\Qbar/\Q)$ on dessins, following Grothendieck.  

The authors would like to thank Joseph Cyr, Jonathan Godbout, Alex Levin, and Michael Novick for their initial contributions to this project.  The authors would also like to thank Hartmut Monien, Richard Foote, and Jeroen Sijsling for useful conversations and the anonymous referee for helpful feedback.  
Our calculations are performed in the computer algebra system \textsf{Magma} \cite{Magma}.  This project was supported by an NSF CAREER Award (DMS-1151047).

\section{Background} \label{sec:background}

In this section, we review three-point branched covers of Riemann surfaces.  For a reference, see the introductory texts of Miranda \cite{Miranda} and Girondo--Gonz\'{a}lez-Diez \cite{GirondoGonzalezDiez}.

\subsection*{Triangle groups and \Belyi\ maps}

By the Riemann existence theorem, Riemann surfaces are the same as (algebraic) curves over $\C$; we will pass between these categories without mention.  

A \defi{\Belyi\ map} over $\C$ is a morphism $\phi:X \to \PP_\C^1$ of curves over $\C$ that is unramified outside $\{0,1,\infty\}$.  Two \Belyi\ maps $\phi,\phi':X,X' \to \PP^1$ are \defi{isomorphic} if there is an isomorphism $\iota:X \to X'$ such that $\phi = \phi' \circ \iota$.  \Belyi\ \cite{Belyi,Belyi2} proved that a curve $X$ over $\C$ can be defined over $\Qbar$, an algebraic closure of $\Q$, if and only if $X$ admits a \Belyi\ map, so we may work equivalently with curves defined over $\Qbar$ or over $\C$.  

A \Belyi\ map $\phi:X \to \PP^1$ admits a short combinatorial description as a \defi{transitive permutation triple} $\sigma = (\sigma_0, \sigma_1, \sigma_{\infty}) \in S_d^3$ such that $\sigma_{\infty} \sigma_1 \sigma_{0} = 1$ and $\langle \sigma_0,\sigma_1, \sigma_{\infty} \rangle \leq S_d$ is a transitive subgroup, according to the following lemma.

\begin{lem} \label{lem:biject}
There is a bijection between transitive permutation triples up to simultaneous conjugacy and isomorphism classes of \Belyi\ maps over $\C$ (or $\Qbar$).
\end{lem}

By the phrase \defi{up to simultaneous conjugation}, we mean that we consider two triples $\sigma,\sigma' \in S_d^3$ to be equivalent if there exists $\rho \in S_d$ such that $\sigma=\rho\sigma\rho^{-1}$, i.e.,
\[ (\sigma_0',\sigma_1',\sigma_\infty') = \rho(\sigma_0,\sigma_1,\sigma_{\infty})\rho^{-1} = (\rho\sigma_0\rho^{-1}, \rho\sigma_1\rho^{-1}, \rho\sigma_{\infty}\rho^{-1}). \]
In this correspondence, the cycles of the permutation $\sigma_0,\sigma_1,\sigma_\infty$ correspond to the points of $X$ above $0,1,\infty$, respectively, and the length of the cycle corresponds to its multiplicity.  Note in particular that there are only finitely many $\Qbar$-isomorphism classes of curves $X$ with a \Belyi\ map of bounded degree.  

The proof of Lemma \ref{lem:biject} has many strands and is a piece of Grothendieck's theory of dessin d'enfants.  The bijection can be realized explicitly using the theory of triangle groups and this constructive proof was the starting point of this paper and is the subject of this section.  For an introduction to triangle groups, including their relationship to \Belyi\ maps and dessins d'enfants, see the surveys of Wolfart \cite{WolfartDessins,WolfartObvious}, and for a more classical treatment see Magnus \cite{Magnus}.

Let $a,b,c \in \Z_{\geq 2} \cup \{\infty\}$ with $a \leq b \leq c$.  We will call the triple $(a,b,c)$ \defi{spherical}, \defi{Euclidean}, or \defi{hyperbolic} according to whether the value
$$
\chi(a,b,c) = 1 - \frac{1}{a} - \frac{1}{b} - \frac{1}{c}  
$$
is respectively negative, zero, or positive, and we associate the geometry
\[ 
H=
\begin{cases} 
\text{sphere $\PP^1_{\C}$},  & \text{ if $\chi(a,b,c)<0$}; \\
\text{plane $\C$}, & \text{ if $\chi(a,b,c)=0$}; \\
\text{upper half-plane $\calH$}, & \text{ if $\chi(a,b,c)>0$}.
\end{cases} \]
accordingly.  The spherical triples are $(2,3,3)$, $(2,3,4)$, $(2,3,5)$, and $(2,2,c)$ for $c \in \Z_{\geq 2}$.  The Euclidean triples are $(2,2,\infty)$, $(2,3,6)$, $(2,4,4),$ and $(3,3,3)$.  All other triples are hyperbolic.

\begin{figure}[h] 
\includegraphics{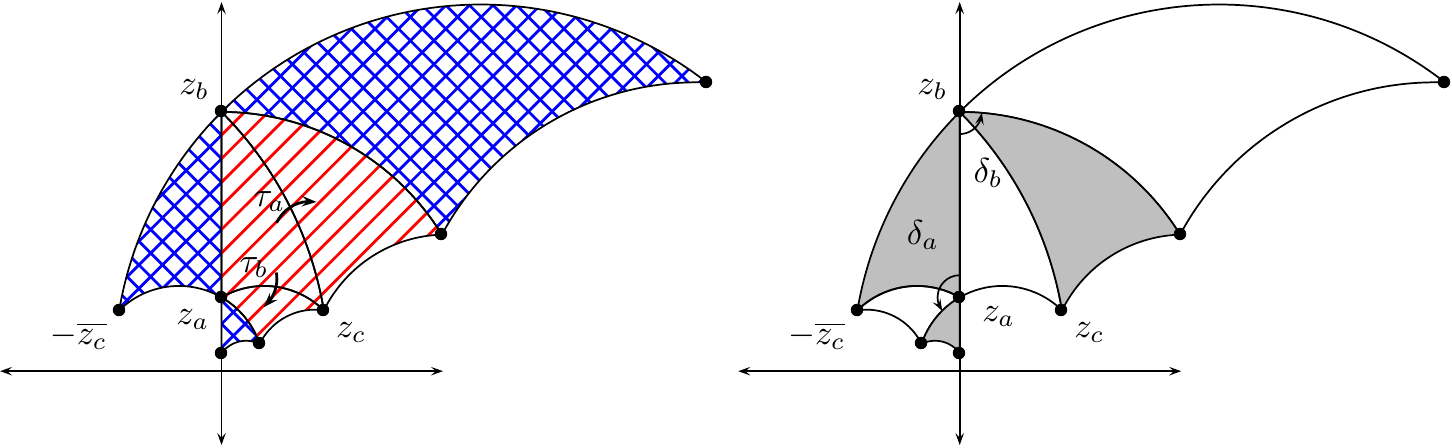}

\caption{A normalized triangle and its images under the reflections $\tau_a$ and $\tau_b$, and under the rotations $\delta_a = \tau_c \tau_b$ and $\delta_b = \tau_a \tau_c$.}  \label{fig:goabc}
\end{figure}

For each triple, we have the associated \defi{triangle group} 
\begin{equation} \label{eqn:Deltaabcpres}
\Delta(a,b,c) = \langle \delta_a, \delta_b, \delta_c \mid \delta_a^a = \delta_b^b = \delta_c^c = \delta_c \delta_b \delta_a =1 \rangle
\end{equation}
analogously classified as spherical, Euclidean, or hyperbolic.  The triangle group $\Delta(a,b,c)$ has a simple geometric description.  Let $T$ be a triangle with angles $\pi/a, \pi/b, \pi/c$ on $H$.  Let $\tau_a, \tau_b,$ and $\tau_c$ be the reflections in the three sides of $T$.  The group generated by $\tau_a, \tau_b,$ and $\tau_c$ is a discrete group with fundamental domain $T$.  The subgroup of orientation-preserving isometries in this group is generated by 
\[ \delta_a = \tau_c \tau_b,\ \delta_b = \tau_a \tau_c,\text{ and }\delta_c = \tau_b \tau_a \] 
and has fundamental domain $D_\Delta$ a copy of $T$ and its mirror.  The above relations can be seen in Figure \ref{fig:goabc}.  For instance, note that $\delta_b(D_\Delta) = \tau_a(D_\Delta)$, the only difference being that the shaded and unshaded triangles are reversed.  Precomposing $\tau_a$ by $\tau_c$ corrects this reversal, hence $\delta_b = \tau_a \tau_c$.  The other relations follow similarly.

From these gluing relations, we see that $\Delta \backslash H \cong \PP^1 \setminus S$ as Riemann surfaces, where $S$ consists of a set of points indexed by indices with $a,b,c=\infty$.  So if $a,b,c<\infty$, the case of interest here, then we simply have $\Delta \backslash H \cong \PP^1$.  From now on we restrict to the case where $a,b,c<\infty$; for the remaining cases, one has a different approach not pursued here.

\begin{rmk}
Mapping elements to their inverses gives an isomorphism from the group $\Delta(a,b,c)$ to the group with the relation with $\delta_c \delta_b \delta_a = 1$ replaced by $\delta_a \delta_b \delta_c = 1$ with the same fundamental domain.  The choice corresponds to an orientation of the Riemann surface $\Delta(a,b,c) \backslash H$; we have chosen the counterclockwise orientation, so that $\delta_s$ acts by counterclockwise rotation by $2\pi/s$ for $s=a,b,c$; and with this choice, we obtain the presentation \eqref{eqn:Deltaabcpres}.
\end{rmk}

\subsection*{Conventions}

First, a convention on permutations (following a convention in computational group theory): we compose permutations from left to right.  This is contrary to the usual functional composition, so to avoid confusion we will write it in the exponentiated form: for example, if $\sigma=(1\ 2)$ and $\tau=(2\ 3)$, then we write $1^{\sigma\tau} = (1^\sigma)^\tau = 2^\tau = 3$, so $\sigma\tau = (1\ 3\ 2)$.  In this way, the action of $S_d$ on $\{1,\dots,d\}$ is on the right.

There is an injective map which, given a transitive permutation triple $\sigma=(\sigma_0,\sigma_1,\sigma_\infty) \in S_d^3$ with orders $a,b,c \in \Z_{\geq 2}$ up to simultaneous conjugation, associates a subgroup $\Gamma \leq \Delta=\Delta(a,b,c)$ of index $[\Delta:\Gamma]=d$ up to conjugation.  Explicitly, to the permutation triple $\sigma$ we associate 
\begin{equation} \label{eqn:gammaperm}
 \Gamma = \pi^{-1}\Stab(1) = \{ \delta \in \Delta : 1^{\pi(\delta)} = 1\},
\end{equation}
where
\begin{equation} \label{eqn:piperm}
\begin{aligned}
\pi : \Delta &\to S_d \\
\delta_a, \delta_b, \delta_c &\mapsto \sigma_0, \sigma_1, \sigma_\infty
\end{aligned}
\end{equation}
is the permutation representation induced by $\sigma$.  

This map has a right inverse: given $\Gamma \leq \Delta$ of index $d$, we define a permutation triple as follows.  We choose cosets $\Gamma\alpha_1,\dots,\Gamma\alpha_d$ of $\Gamma$, and define a permutation representation $\pi:\Delta \to S_d$ by
\begin{equation} \label{eqn:deltaaction}
i^{\pi(\delta)}=j \quad \text{if} \quad \Gamma \alpha_i \delta = \Gamma \alpha_j,
\end{equation}
for $\delta \in \Delta$ and $i,j \in \{1,\dots,d\}$.  We then let $\sigma_0=\pi(\delta_a)$, $\sigma_1=\pi(\delta_b)$, and $\sigma_\infty=\pi(\delta_c)$.  
This defines a group homomorphism, but it changes the left action of $\Delta$ on $\calH$ into a right action of $\Delta$ on $\{1, \dots, d\}$, according to our convention.  Indeed, if $\Gamma \alpha_i \delta = \Gamma \alpha_j$ and $\Gamma \alpha_j \delta' = \Gamma \alpha_k$, so that $i^{\pi(\delta)}=j$ and $j^{\pi(\delta')}=k$, then
\[ \Gamma \alpha_i \delta\delta' = \Gamma \alpha_j \delta' = \Gamma \alpha_k \]
so 
\[ i^{\pi(\delta\delta')} = k = j^{\pi(\delta')} = (i^{\pi(\delta)})^{\pi(\delta')} = i^{\pi(\delta)\pi(\delta')}. \]

If we begin with the permutation triple $\sigma$ and associate the group $\Gamma$ as in (\ref{eqn:gammaperm}) and choose cosets $\Gamma\alpha_i$ with the property that
\begin{equation} \label{eqn:permtocoset}
1^{\pi(\alpha_i)} = i 
\end{equation}
(well-defined, since $1^{\pi(\gamma\alpha_i)}=1^{\pi(\gamma)\pi(\alpha_i)}=1^{\pi(\alpha_i)}=i$ for all $\gamma \in \Gamma$ by definition).  Since $i^{\sigma_0} = j = i^{\pi(\delta_a)}$ if and only if 
\[ \Gamma \alpha_i \delta_a = \Gamma \alpha_j \]
and similarly for $\sigma_1$ (resp.~$\sigma_\infty$) and $\delta_b$ (resp.~$\delta_c$), we recover $\sigma$ as claimed.  

\begin{rmk}
These maps can be made into bijections by artificially restricting the image to those subgroups such that $\sigma_0,\sigma_1,\sigma_{\infty}$ have orders $a,b,c$, i.e., eliminating those subgroups that ``come from a smaller triple''.  
\end{rmk}

\subsection*{Correspondence}

We now prove Lemma \ref{lem:biject}, as sketched in the introduction.  
As mentioned above, $\Delta$ (and also $\Gamma$) act on $H$ on the left and the quotient map  $\phi: X = \Gamma \backslash H \to \Delta \backslash H \cong \PP^1$ is then a \Belyi\ map.

We first show that this association is well-defined on simultaneous conjugacy classes of permutation triples.  Suppose $\sigma,\sigma' \in S_d^3$ are simultaneously conjugate with $\sigma'=\rho\sigma\rho^{-1}$ for some $\rho \in S_d$.  The permutations of these triples have the same respective orders, hence the same triangle group $\Delta=\Delta(a,b,c)$ is associated to both triples.  Let $\pi,\Gamma$ be as in \eqref{eqn:gammaperm}, \eqref{eqn:piperm}.  Then $\Gamma' = \pi^{-1} \Stab(k)$ where $k=1^\rho$: if a word $\upsilon'$ written in $\sigma_0', \sigma_1', \sigma_\infty'$ fixes 1, then $\upsilon' = \rho \upsilon \rho^{-1}$, where $\upsilon$ is the corresponding word written in $\sigma_0, \sigma_1, \sigma_\infty$, so $\upsilon$ fixes $1^\rho = k$.  Since $\langle \sigma_0, \sigma_1, \sigma_\infty \rangle \leq S_d$ is transitive, then there exists $\mu \in \langle \sigma_0, \sigma_1, \sigma_\infty \rangle$ with $1^\mu = k$.  Then
$$
\Gamma' = (\pi')^{-1} \Stab(1) = \pi^{-1} \Stab(k) = \pi^{-1}(\mu \Stab(1) \mu^{-1}) \, .
$$
Choosing $\delta \in \pi^{-1}(\mu)$, then
$$
\Gamma' = \pi^{-1}(\mu \Stab(1) \mu^{-1}) = \delta \pi^{-1}(\Stab(1)) \delta^{-1} = \delta \Gamma \delta^{-1}
$$
so these subgroups are conjugate; and the map
\begin{align*}
\iota : \Gamma \backslash H &\mapsto \Gamma' \backslash H\\
z &\mapsto \delta z
\end{align*}
is an isomorphism of Riemann surfaces with $\phi' \circ \iota = \phi$, where $\phi : \Gamma \backslash H \to \Delta \backslash H$ and $\phi' : \Gamma' \backslash H \to \Delta \backslash H$ are the quotient maps described above.

We now construct the inverse map to the correspondence in Lemma \ref{lem:biject}.  Let $\phi: X \to \PP^1$ be a \Belyi\ map of degree $d$, let $Y = X \setminus \phi^{-1}(\{0, 1, \infty\})$ and let $U = \PP^1 \setminus \{0, 1, \infty\}$.  Then the induced map $\phi|_Y : Y \to U$ is an unramified covering map.  Let $* \in U$ be a basepoint.  We have a presentation of the fundamental group as
\[ \pi_1(U,*)=\langle \eta_0, \eta_1, \eta_\infty \mid \eta_\infty \eta_1 \eta_0 = 1 \rangle \]
where $\eta_0,\eta_1,\eta_\infty$ are homotopy classes represented by loops around 0, 1, and $\infty$, respectively.  For each $x \in \phi^{-1}(*)$, a path $\eta$ in $U$ with initial point $*$ can be lifted to a unique path $\widetilde{\eta}$ in $Y$ with initial point $x$ such that $f \circ \widetilde{\eta} = \eta$.  Thus the terminal point of $\widetilde{\eta}$ will be a unique $x' \in \phi^{-1}(*)$, and this induces a right action of $\pi_1(U,*)$ on $\phi^{-1}(*)$ by $x^\eta = \widetilde{\eta}(1)$.  If we label the $d$ points in $\phi^{-1}(*)$ with $\{1,\dots,d\}$, then the action yields a permutation representation $\sigma : \pi_1(U, *) \to S_d$.  Letting $\sigma_0 = \sigma(\eta_0)$, $\sigma_1 = \sigma(\eta_1)$, and $\sigma_\infty = \sigma(\eta_\infty)$ yields a permutation triple.  Moreover, since $Y$ is path-connected, the group $\langle \sigma_0, \sigma_1, \sigma_\infty \rangle \leq S_d$ is transitive.  

\begin{rmk}
Instead of considering $U=\PP^1 \setminus \{0,1,\infty\}$, one can also consider the orbifold fundamental group of $X(\Delta)=\Delta \backslash H$, which is exactly $\pi_1(X(\Delta),*) \cong \Delta$ (with a choice of basepoint $*$).  Then the map $X(\Gamma) \to X(\Delta)$ is a cover in the orbifold category, and the argument in the previous paragraph becomes a tautology.
\end{rmk}

We now show that this map is well-defined on the set of isomorphism classes of \Belyi\ maps.  Let $\phi : X \to \PP^1$ and $\phi' : X' \to \PP^1$ be isomorphic \Belyi\ maps.  Then there is an isomorphism $\iota:X \to X'$ of Riemann surfaces such that $\phi = \phi' \circ \iota$.  The restriction 
\[ \iota|_{\phi^{-1}(*)} : \phi^{-1}(*)=\{x_1,\dots,x_d\} \to \{x_1',\dots,x_d'\}=(\phi')^{-1}(*) \] 
is a bijection and so if these sets are labelled, we obtain a relabeling $\rho \in S_d$ with $\iota(x_i)=x_j'$ if and only if $i^\rho=j$.  Let $\sigma,\sigma'$ be the permutation triples associated to $\phi,\phi'$ above.  Then 
\[ (\sigma_0, \sigma_1, \sigma_\infty)  = (\rho {\sigma_0}' \rho^{-1}, \rho {\sigma_1}' \rho^{-1}, \rho {\sigma_\infty}' \rho^{-1}) \]
so the triples are simultaneously conjugate, as desired.

To complete the proof of Lemma \ref{lem:biject}, we check directly that these maps are inverses; an explicit version of this is the subject of section \ref{sec:cosets}.

\begin{rmk}
The triple $(\sigma_0, \sigma_1, \sigma_\infty)$ encodes several pieces of information about the associated ramified covering.  The number of disjoint cycles in $\sigma_0$ (resp., $\sigma_1$, $\sigma_\infty$) is the number of distinct points in the fiber above $0$ (resp., $1$, $\infty$).  The ramification indices of the points in the fiber are given by the lengths of these cycles.
\end{rmk}

\subsection*{Genus, passports, and moduli}

The genus of a \Belyi\ map is given by the Riemann-Hurwitz formula.  If we define the \defi{excess} $e(\tau)$ of a cycle $\tau \in S_n$ to be its length minus one, and the excess $e(\sigma)$ of a permutation to be the sum of the excesses of its disjoint cycles (also known as the \defi{index} of the permutation, equal to $d$ minus the number of orbits), then the genus of a \Belyi\ map of degree $d$ with monodromy $\sigma$ is 
\begin{equation} \label{eqn:RH}
g = 1-d + \frac{e(\sigma_0)+e(\sigma_1)+e(\sigma_\infty)}{2}. 
\end{equation}

A \defi{refined passport} consists of the data of a transitive subgroup $G \leq S_d$ and three conjugacy classes $C_0,C_1,C_\infty$ in $G$.  The refined passport of a triple $\sigma$ is given by the group $G= \langle \sigma \rangle$ and the conjugacy classes of the elements of $\sigma$.  The \defi{size} of a refined passport is the cardinality of the set
\[ \{\sigma \in S_d : \sigma_i \in C_i \text{ for $i=0,1,\infty$ and $\langle \sigma \rangle = G$}\}/\!\sim \]
where $\sim$ denotes simultaneous conjugation in $S_d$.  If a refined passport has size $1$, then we call the refined passport (and any associated triple) \defi{rigid}.  For more on rigid triples, with applications to the inverse Galois problem, see e.g.\ Serre \cite{Serre}.  For more on issues related to passports and moduli, see Sijsling--Voight \cite{SijslingVoight}.

Let $X$ be a curve defined over $\C$.  The \defi{field of moduli} $M(X)$ of $X$ is the fixed field of the group $\{\tau \in \Aut(\C) : X^{\tau} \cong X\}$, where $X^{\tau}$ is the base change of $X$ by the automorphism $\tau \in \Aut(\C)$ obtained by applying $\tau$ to any set of defining equations for $X$.  If $F$ is a field of definition for $X$ then clearly $F \supseteq M(X)$.  If $X$ has a minimal field of definition $F$ then necessarily $F=M(X)$.  However, not all curves can be defined over their field of moduli.  For more information, see work of Coombes and Harbater \cite{CoombesHarbater}, D\`ebes and Emsalem \cite{DebesEmsalem}, and K\"ock \cite{Kock}.  Similarly, we define the field of moduli $M(\phi)$ of a \Belyi\ map $\phi$ to be the fixed field of $\{ \tau \in \Aut(C) : \phi^{\tau} \cong \phi\}$.  The \defi{field of rationality} of a refined passport is the field $\Q(\chi(C_i))$ obtained by adjoining the character values on the conjugacy classes $C_0,C_1,C_\infty$ in $G$.  In general, the degree of the field of moduli of a \Belyi\ map over its \defi{field of rationality} is bounded above by the size of its refined passport.  Again, a field of definition of $\phi$ may be bigger than the field of moduli (itself bigger than the field of rationality), and these issues quickly become quite delicate: see Sijsling--Voight \cite[Section 7]{SijslingVoight}.  

\subsection*{Enumerating triples}

To use this bijection to provide tables of \Belyi\ maps, we will need to enumerate triples: it is straightforward, using computational group theory techniques, to enumerate all transitive permutation triples up to simultaneous conjugacy.  Indeed, for $d \leq 30$, the set of transitive groups of degree $d \leq 30$ have been classified \cite{Hulpke}, and these are available in \textsf{Magma} \cite{Magma} and given a unique identifier.  However, the number of such triples grows rapidly with $d$, and so computing a complete database of triples is only practical for small degree.  A simple way to compute these triples using double cosets is as follows.

\begin{lem} \label{lem:c0c1}
Let $G$ be a group, let $C_0,C_1$ be conjugacy classes in $G$ represented by $\tau_0,\tau_1$.  Then the map
\begin{align*}
C_G(\tau_0) \backslash G / C_G(\tau_1) &\to \{(\sigma_0,\sigma_1) : \sigma_0 \in C_0, \sigma_1 \in C_1 \}/\!\sim_{G} \\
C_G(\tau_0) g C_G(\tau_1) &\mapsto (\tau_0, g \tau_1 g^{-1})
\end{align*}
is a bijection, where the latter set is taken up to simultaneous conjugation in $G$.
\end{lem}
\noindent
(Here $C_G(\tau_i)$ denotes the centralizer of $\tau_i$ in $G$.)

One completes the triples arising from Lemma \ref{lem:c0c1} by taking $\sigma_\infty=(\sigma_1\sigma_0)^{-1}$; enumerating over all conjugacy classes of subgroups $G \leq S_d$ and all conjugacy classes $C_0,C_1$ in $G$ then lists all triples up to simultaneous conjugation.  To avoid overcounting, one will want to restrict to those triples $\sigma$ obtained in Lemma \ref{lem:c0c1} which generate $G$ (and not a proper subgroup).  

\begin{proof}
First, surjectivity.  Let $\sigma_0 \in C_0$ and $\sigma_1 \in C_1$.  Conjugating in $G$, we may assume that $\sigma_0=\tau_0$, and then indeed $\sigma_1=g \tau_1 g^{-1}$ for some $g \in G$.  Next, injectivity: if $(\tau_0,g\tau_1 g^{-1})$ and $(\tau_0, g'\tau_1 (g')^{-1})$ are simultaneously conjugate in $G$ then there exists $c \in C_G(\tau_0)$ such that $cg\tau_1 (cg)^{-1} = g' \tau_1 (g')^{-1}$ so $(g')^{-1} cg \in C_G(\tau_1)$ and $C_G(\tau_0) g C_G(\tau_1) = C_G(\tau_0) g' C_G(\tau_1)$ as claimed. 
\end{proof}

There are efficient algorithms for computing double cosets \cite{Linton}, so Lemma \ref{lem:c0c1} has an advantage over more naive enumeration.

\subsection*{Spherical and Euclidean \Belyi\ maps}

To conclude this section, we describe the calculation of the spherical and Euclidean \Belyi\ maps.  The spherical triangle groups $\Delta(a,b,c)$ are finite groups acting on the sphere $\PP^1$, and the corresponding quotients correspond to the classically known finite subgroups of $\PSL_2(\C)$.  The case $(2,2,c)$ corresponds to dihedral groups and they are given explicitly by Chebyshev polynomials.  This leaves the three cases $\Delta(2,3,3)=A_4$, $\Delta(2,3,4)=S_4$, and $\Delta(2,3,5)=A_5$ and these are referred to as the tetrahedral, octahedral, and icosahedral groups respectively because they are the groups of rigid motions of the corresponding Platonic solids.  We refer the reader to Magnus \cite{Magnus} and the references therein for further discussion.  The \Belyi\ maps for these cases are also classical: see Magot--Zvonkin \cite{MagotZvonkin} for references and a discussion of maps arising from Archimedean solids.  (One can also view the spherical dessins using hyperbolic triples, including $\infty$: see recent work of He--Read \cite{HeRead} for a complete description.)

The Euclidean groups $\Delta(a,b,c)$ are infinite groups associated to flat tori.  First, we have that $\Delta(3,3,3) \leq \Delta(2,3,6)$ with index $2$, so it suffices to consider only the triples $(2,3,6)$ and $(2,4,4)$, associated to classical triangulations of the plane.  First consider the triple $(2,3,6)$, to which we associate the elliptic curve $E(2,3,6) : y^2 = x^3 - 1$ with CM by $\Z[(-1+\sqrt{-3})/2]=\Z[\omega]$.  The quotient map $x^3: E(2,3,6) \to \PP^1$ of $E$ by its automorphism group (as an elliptic curve)  $\Aut(E(2,3,6)) = \la -\omega \ra \cong \Z/6\Z$ gives a Galois \Belyi\ map of degree $6$ corresponding to the permutation triple 
\[ \sigma_0=(1\ 4)(2\ 5)(3\ 6), \quad \sigma_1=(1\ 3\ 5)(2\ 4\ 6),\quad \sigma_{\infty}=(1\ 2\ 3\ 4\ 5\ 6) \] 
and subgroup $\Gamma(E) \leq \Delta(2,3,6)$ of index $6$.  For any other finite index subgroup $\Gamma$, the intersection $\Gamma \cap \Gamma(E)$ corresponds to an elliptic curve $X$ equipped with a isogeny $X \to E(2,3,6)$ induced by the inclusion $\Gamma \cap \Gamma(E) \hookrightarrow \Gamma(E)$.  Similar statements hold for $\Delta(2,4,4)$, with the elliptic curve $E(2,4,4): y^2 = x^3 - x$ with CM by $\Z[i]$ and \Belyi\ map $x^2: E(2,4,4) \to \PP^1$ of degree $4$.  In other words, all \Belyi\ maps with these indices can be written down explicitly via CM theory of these two elliptic curves \cite{Silverman,Silverman2}, and so other methods are easier to apply.  See also work of Singerman and Sydall \cite{SingermanSydall}.

\section{Embedding Fuchsian triangle groups} \label{sec:embedtriangle}

Let $(a,b,c)$ be a hyperbolic triple.  In this section, we explicitly construct a hyperbolic triangle $T$ in the hyperbolic plane $H=\calH$ with angles $\pi / a$, $\pi / b$, and $\pi / c$; see also Petersson \cite{Petersson}.

\begin{figure}[h] 
\includegraphics{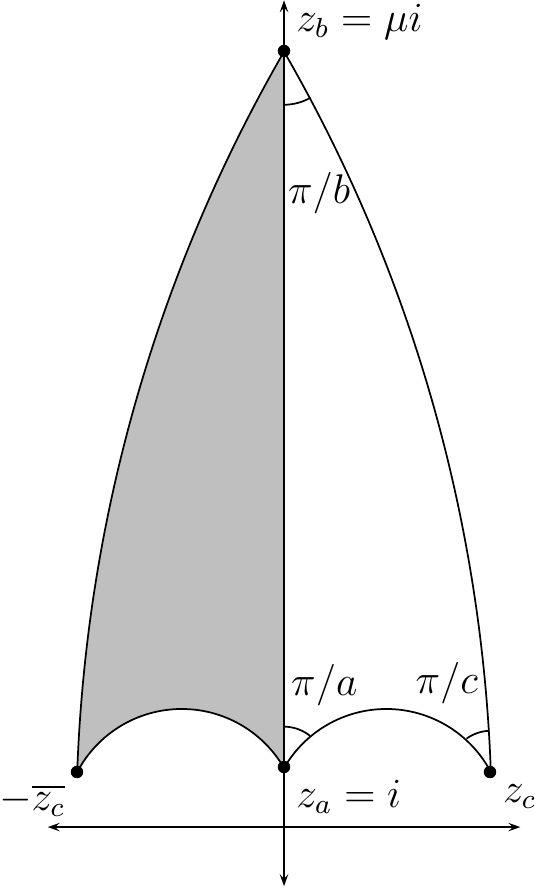}

\caption{Normalized triangle with angles $\pi/a,\pi/b,\pi/c$ and its mirror, comprising $D_\Delta$.} \label{fig:tabc} 
\end{figure}

We  assume that one side of the triangle is on the imaginary axis with one vertex at $i$.  We label the vertices $z_a, z_b,$ and $z_c$ corresponding to the angle at each vertex.  Recall that the hyperbolic metric on $\calH$ is given by
\begin{equation} \label{eqn:hyperbolicmetric} 
\begin{aligned} 
d : \calH \times \calH &\to \R_{\geq 0}  \\
\cosh d(x_1+iy_1, x_2+iy_2) &= 1 + \frac{(x_2-x_1)^2 + (y_2-y_1)^2}{2 y_1 y_2}.
\end{aligned}
\end{equation}
Since $\PSL_2(\R)$ acts transitively on $\calH$ by orientation-preserving isometries, we may take $z_a=i$; further composing by a local rotation around $i$, we may assume that $z_b=\newt i$ for some $\newt \in \R_{\geq 1}$.  Finally, reflecting in the imaginary axis, we may assume that $\Real(z_c) > 0$.  The triangle $T(a,b,c)$ is then unique, and in particular the positions of $z_b$ and $z_c$ are determined by insisting that $T(a,b,c)$ have angles $\pi / a,\pi / b,\pi / c$ at $z_a, z_b,z_c$, as in Figure \ref{fig:tabc}.

We can solve for $\newt$ by recalling the law of cosines from hyperbolic trigonometry (see Ratcliffe \cite[Theorem 3.5.4]{Ratcliffe}).  

\begin{prop} \label{prop:lawcosines}
If $\alpha, \beta, \gamma$ are the angles of a hyperbolic triangle and $C$ is the (hyperbolic) side length of the side opposite $\gamma$, then
\begin{equation*} 
\cosh(C) = \frac{\cos \alpha \cos \beta + \cos \gamma}{\sin \alpha \sin \beta}.
\end{equation*}
\end{prop}

Using (\ref{eqn:hyperbolicmetric}) and Proposition \ref{prop:lawcosines} we have
$$
1 + \frac{(\newt-1)^2}{2\newt} = \frac{\cos \pi/a \cos \pi/b + \cos \pi/c}{\sin \pi/a \sin \pi/b}
$$
so writing
$$
\lambda = \frac{\cos\pi/a\cos\pi/b+\cos\pi/c}{2\sin\pi/a\sin\pi/b}
$$
we have
$$
\newt = \lambda + \sqrt{\lambda^2 -1}.
$$
Note that we are justified in taking the positive square root since we assume that $\newt > 1$.  

Having found $z_b$ we now find equations for the geodesics through $z_a$ and $z_b$ that intersect the imaginary axis with the angles $\pi /a$ and $\pi /b$ respectively. These geodesics are given by the intersection of $\calH$ with the set of points $x+iy \in \C$ satisfying the equations
\begin{equation} \label{eqn:circles}
\begin{aligned}
(x- \cot \pi/a)^2 + y^2 &= \csc^2 \pi/a \\
(x + \newt \cot \pi /b)^2 + y^2 &= \newt^2 \csc^2 \pi/b.
\end{aligned} 
\end{equation}
We then find the unique point of intersection of the two circles \eqref{eqn:circles} in $\calH$ as
\begin{equation} \label{eqn:zc}
z_c = \frac{\newt^2-1}{2(\cot \pi/a + \newt \cot \pi/b)} + i\sqrt{\csc^2 \pi/a - \left( \frac{\newt ^2-1}{2(\cot \pi/a + \newt \cot \pi/b)} - \cot \pi/a \right ) ^2}.
\end{equation}
Note that $z_c \in \Qbar$, in fact, $z_c$ lies in an at most quadratic extension of $\Q(\zeta_{2a},\zeta_{2b},\zeta_{2c})$, where $\zeta_s=\exp(2\pi i/s)$.  

We now find the elements in $\PSL_2(\R)$ that realize $\delta_a$ and $\delta_b$ geometrically, yielding an explicit embedding of $\Delta(a,b,c) \hookrightarrow \PSL_2(\R)$.  Recall from section \ref{sec:background} that the image of $\delta_a$ in $\PSL_2(\R)$ is equal to $\tau_c \tau_b$.

Therefore, the element of $\PSL_2(\R)$ corresponding to $\delta_a$ will fix $z_a$ and $\overline{z_a}$, while sending the vertex $z_c$ to $-\overline{z_c}$.  Similarly, the element of $\PSL_2(\R)$ corresponding to $\delta_b$ will fix $z_b$ and $\overline{z_b}$, while sending the vertex $-\overline{z_c}$ to $z_c$.  Since a linear fractional transformation is uniquely determined by its action on three points, these conditions uniquely determine the elements of $\PSL_2(\R)$ corresponding to $\delta_a$ and $\delta_b$.  

The element $\delta_a$ acts by rotation around $i$ by the angle $2\pi/a$, and conjugating by the matrix $\begin{pmatrix} \sqrt{\mu} & 0 \\ 0 & 1/\sqrt{\mu} \end{pmatrix}$ we obtain a matrix for $\delta_b$, as follows (cf.\ Petersson \cite{Petersson}).

\begin{prop} \label{prop:triangleembed}
We have an embedding
\begin{align*}
\Delta(a,b,c) &\hookrightarrow \PSL_2(\R) \\
\delta_a &\mapsto \begin{pmatrix} \cos(\pi /a) & \sin(\pi /a) \\ -\sin(\pi /a) & \cos(\pi /a)  \end{pmatrix} \\
\delta_b &\mapsto \begin{pmatrix} \cos(\pi / b) & \newt\sin(\pi / b) \\ -\sin(\pi / b )/\newt & \cos(\pi/b) \end{pmatrix}
\end{align*}
where 
\[ \newt = \lambda + \sqrt{\lambda^2 -1} \quad \text{and} \quad \lambda = \frac{\cos\pi/a\cos\pi/b+\cos\pi/c}{2\sin\pi/a\sin\pi/b}.  \]
In this embedding, the fixed points of $\delta_a,\delta_b,\delta_c$ are $z_a=i$, $z_b=\newt i$, and
\[ z_c = \frac{\newt^2-1}{2(\cot \pi/a + \newt \cot \pi/b)} + i\sqrt{\csc^2 \pi/a - \left( \frac{\newt^2-1}{2(\cot \pi/a + \newt \cot \pi/b)} - \cot \pi/a \right ) ^2} \in \overline{\Q}. \]
\end{prop}

By conformally mapping $\calH$ to the unit disc $\calD$ by the map
\begin{align*}
w_p : \calH &\to \calD \\
z &\mapsto w_p(z) = \frac{z-p}{z-\overline{p}} 
\end{align*}
for a point $p \in \calH$, we can extend the action of $\Delta$ to $\calD$ by the action
\begin{align*}
\Delta \times \calD &\to \calD \\
(\delta, w) &\mapsto w_p(\delta w_p^{-1}(w)). 
\end{align*} 

With the identification map $w_{z_a} = w_a : \calH \to \calD$ the element $\delta_a$ then acts on $\calD$ by rotation about $2\pi / a$ about the origin.  To avoid excessive notation, when the identification $w_p : \calH \to \calD$ is clear, we will just identify $D_\Delta$ with its image $w_p(D_\Delta) \subset \calD$.   

\section{Coset enumeration and reduction theory} \label{sec:cosets}

In this section, we exhibit an algorithm that takes as input a permutation triple $\sigma$ and gives as output the cosets of the point stabilizer subgroup $\Gamma \leq \Delta(a,b,c)$; this algorithm also yields a reduction algorithm and a fundamental domain for the action of $\Gamma$ on $\calH$.  Some similar methods were considered in the context of curves with signature $(1;e)$ by Sijsling \cite{Sijsling}.  Throughout, we will use implicitly the conventions of computational group theory for representing finitely presented groups in bits as well as basic algorithms for manipulating them: for a reference, see Holt \cite{Holt} or Johnson \cite{Johnson}.

Throughout, let $\Delta=\Delta(a,b,c)$ be a triangle group with indices $a,b,c \in \Z_{\geq 2}$ and let $\Gamma \leq \Delta$ be a subgroup of finite index $d$.  

\subsection*{Coset graph}

The set $\{\delta_a^{\pm 1}, \delta_b^{\pm 1}\}$ of the generators of $\Delta(a,b,c)$ and their inverses acts on the cosets of $\Gamma$.  We can represent this action as a directed, vertex- and edge-labelled graph with multiple edges allowed.  This multidigraph is similar to the Cayley graph of a group; however, since $\Gamma$ is not assumed to be normal the coset space does not necessarily form a group.  Our algorithm is similar to the classical Todd-Coxeter algorithm (see e.g.\ Rotman \cite{Rotman}) for enumerating the cosets of a subgroup of a group given by a presentation, but it is tailored for our specific application.


The \defi{(symmetric) Cayley graph} for $\Delta$ on the generators $\delta_a,\delta_b$ is a symmetric, directed $4$-regular graph (i.e., every vertex has indegree and outdegree equal to 4) with vertex set $\Delta$ and edge set 
\[ \{\delta \xrightarrow{\epsilon} \delta \epsilon: \delta \in \Delta,\epsilon = \delta_a^{\pm 1}, \delta_b^{\pm 1}\}.  \] 
Consider the quotient of the Cayley graph by the group $\Gamma$, where $\Gamma$ acts on the vertices on the left, identifying $w \sim \gamma w$ for $\gamma \in \Gamma$.  In the quotient, edges map to edges, and so we obtain again a symmetric directed $4$-regular graph with edges labelled by $\delta_a^{\pm 1},\delta_b^{\pm 1}$ and now with vertex set equal to $\Gamma \backslash \Delta$. 

\begin{figure}[h] 
\includegraphics{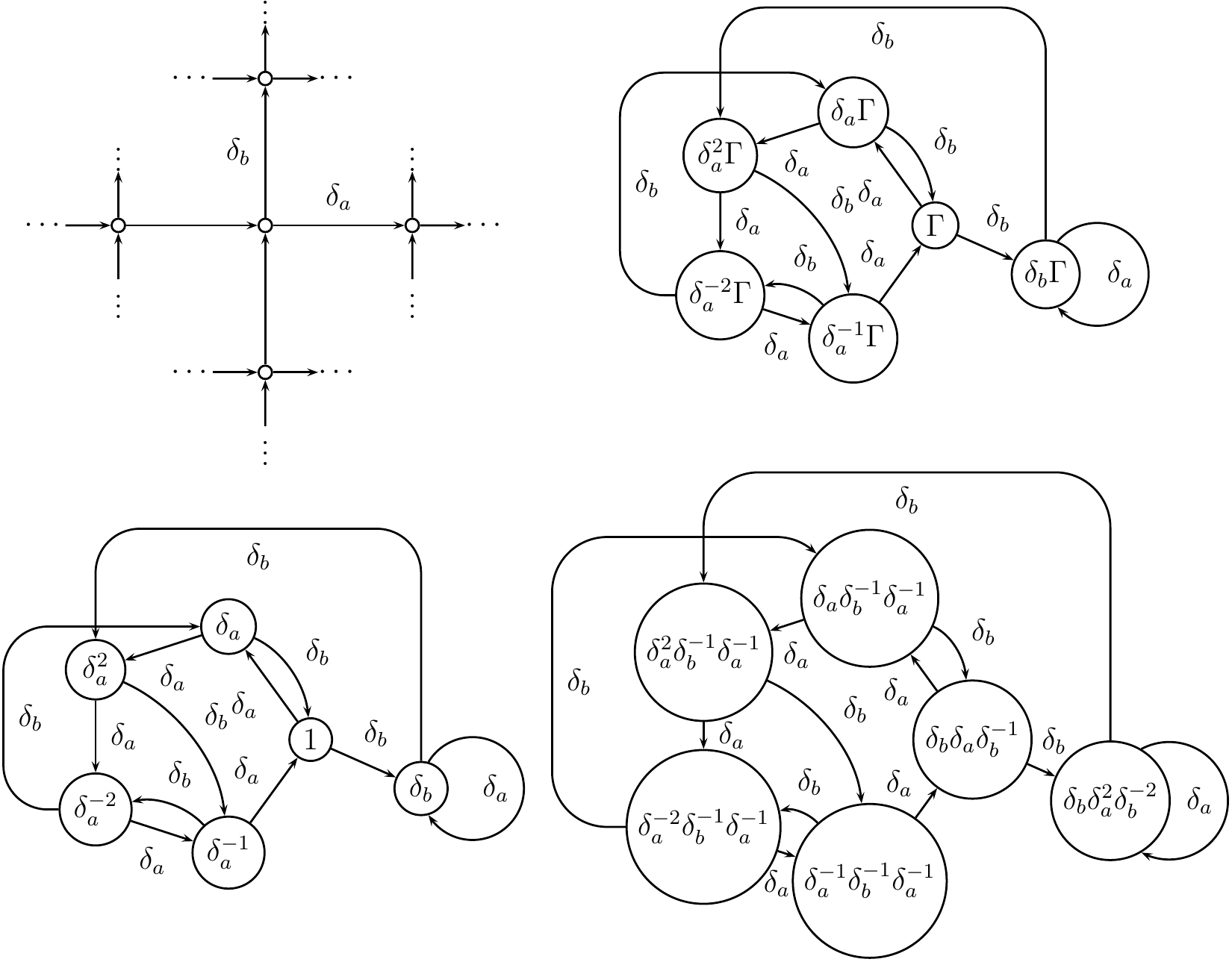}

\caption{Quotient and coset graphs for $\Gamma$ the subgroup of $\Delta(5,6,4)$ corresponding to the permutation triple $(\sigma_0,\sigma_1,\sigma_{\infty})$ with $\sigma_0 = (1\;5\;4\;3\;2)$, $\sigma_1 = (1\;6\;4\;2\;3\;5)$, and $\sigma_\infty = (1\;4\;3\;6)$.} \label{fig:QuotientGraph} 
\end{figure}

\begin{exm}
Let $\Gamma \leq \Delta = \Delta(5,6,4)$ correspond to the permutation triple $\sigma = (\sigma_0,\sigma_1,\sigma_\infty)$ with 
\[ \sigma_0 = (1\;5\;4\;3\;2),\ \sigma_1 = (1\;6\;4\;2\;3\;5),\ \sigma_\infty = (1\;4\;3\;6).  \]
In the top left of Figure \ref{fig:QuotientGraph} is the Cayley graph for $\Delta$ which we quotient by $\Gamma$ to get the quotient graph in the top right. The coset graphs (bottom) exactly resemble the quotient graph. The only difference is that in the quotient graph the vertices are the cosets, whereas in the coset graphs we choose a particular coset representative to label each vertex. In this way there are many different coset graphs for the same quotient graph.
\end{exm}

We will need to choose representatives for the cosets (vertex set) of the aforementioned quotient graph.  So we define a \defi{coset graph} $G=G(\Gamma \backslash \Delta)$ for a permutation representation (\ref{eqn:piperm}) to be the quotient graph of the Cayley graph for $\Delta$ by the group $\Gamma$, as above, but with vertices labelled $\alpha_i$ for $i=1,\dots,d$ representing the cosets $\Gamma \alpha_i$ such that equation (\ref{eqn:permtocoset}) is satisfied.  In both the Cayley graph of $\Delta$ and in its quotient by $\Gamma$ with labels $\delta_a$ or $\delta_b$, we have the corresponding reverse edge labeled with $\delta_a^{-1}$ or $\delta_b^{-1}$ respectively.  

Let $G$ be a coset graph for $\Gamma$ with vertex labels $\alpha_i$.  Then we have a corresponding  fundamental domain $D_\Gamma = D_\Gamma(G)$ for $\Gamma$ given by
\begin{equation} \label{eqn:fundGamma}
 D_\Gamma  = \bigcup_{i=1}^d \alpha_i D_\Delta.
\end{equation}
where $D_\Delta$ is as in Proposition \ref{prop:triangleembed}. The connectedness of the corresponding fundamental domain depends on how we label the vertices in the coset graph. 

If $e:v(\alpha_i) \xrightarrow{\epsilon} v(\alpha_j)$ is an edge of $G$, then we obtain a geodesic segment 
\[ s(e)=\alpha_i D_\Delta \cap \alpha_i\epsilon D_\Delta \subseteq D_\Gamma.  \]
Those edges that correspond to geodesics on the boundary of $D_\Gamma$ we call the \defi{sides} of $G$.  Not all edges of $G$ are sides: those edges that are not sides correspond to geodesics in the interior of $D_\Gamma$. 

\begin{figure}[h]
\includegraphics{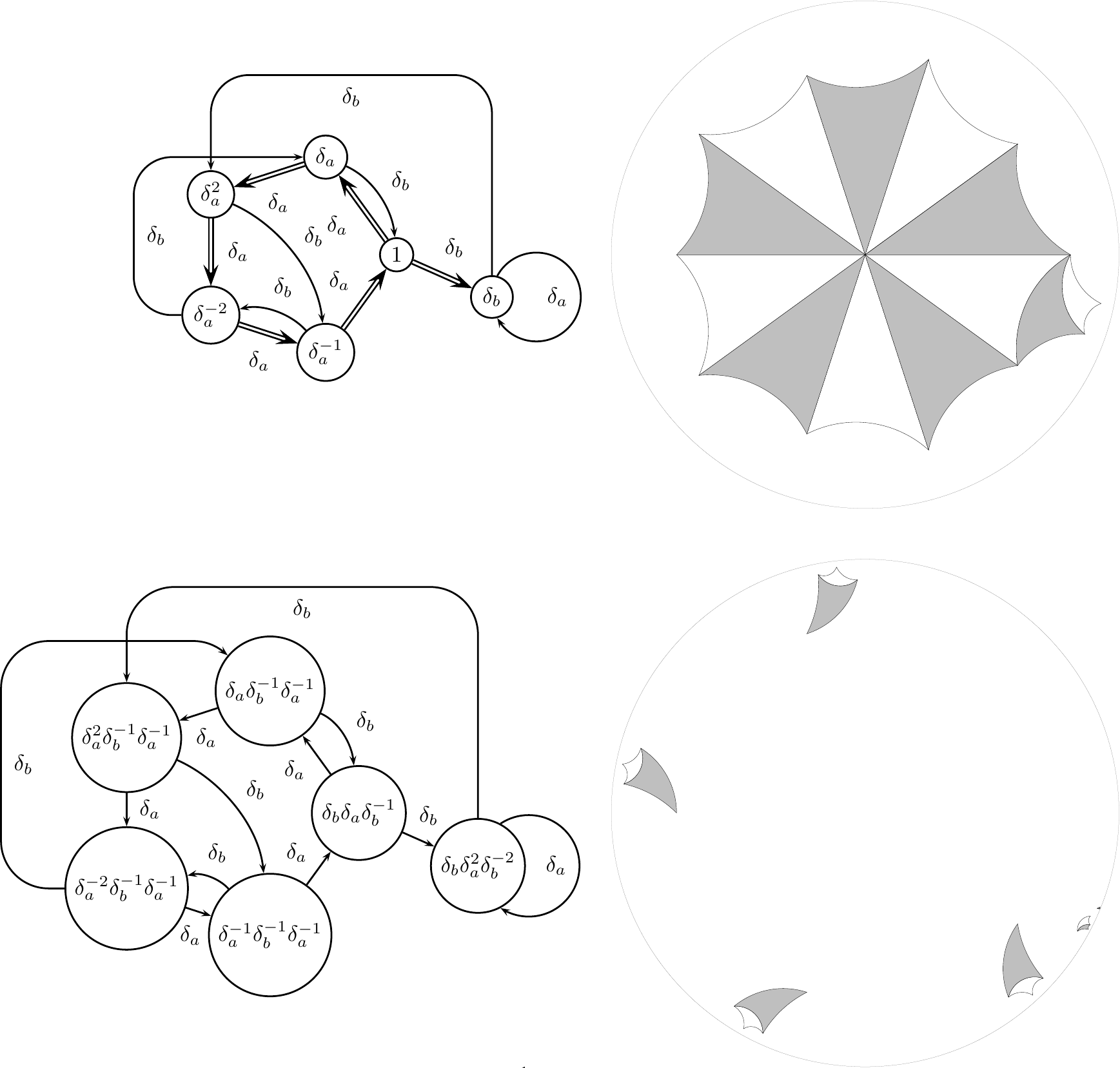}

\caption{The two coset graphs in Figure \ref{fig:QuotientGraph} and their corresponding fundamental domains.}
\label{fig:EdgesSides} 
\end{figure}

\begin{exm}
In Figure \ref{fig:EdgesSides}, the fundamental domain corresponding to the coset graph on the top is connected whereas the fundamental domain corresponding to the coset graph on the bottom is not.  The double arrows are the edges of $G$ that are not sides, and the single arrows are the sides of $G$. In the bottom left of Figure \ref{fig:EdgesSides} is another choice of coset graph for $\Gamma$ as in Figure \ref{fig:QuotientGraph} with the corresponding fundamental domain $D_\Gamma$ on the right.
\end{exm}

Let $S$ be the set of sides of $D_\Gamma$.  A \defi{side pairing} of $G$ is a set of elements 
\begin{equation} \label{eqn:SG}
S(G) = \{ (\gamma, s, s^*) \in \Gamma \times S \times S : \gamma s^* = s \} 
\end{equation}
that defines a labeled equivalence relation on $S$ that induces a partition of $S$ into pairs.  The elements $\gamma$ in (\ref{eqn:SG}) are called the \defi{side pairing elements} of the side pairing $S(G)$.  In bits, we specify a side $s\in S$ as a pair $(i,\epsilon)$ denoting the unique edge $e:v(\alpha_i) \xrightarrow{\epsilon} v(\alpha_j)$ in $G$ corresponding to $s$. 

A side pairing always exists: we give a constructive proof of this fact in the next subsection.

\subsection*{Coset algorithm}

In Algorithm~\ref{alg:cosetalg}, we present an algorithm that to compute a coset graph.  The intuition behind this algorithm is as follows.  We start in the domain $D_\Delta$, corresponding to the identity coset with representative $\alpha_1=1$.  We then consider those translates of the domain $D_\Delta$ which are adjacent to $D_\Delta$ (in the sense that they share a side): with normalizations as in Proposition \ref{prop:triangleembed}, they are
$\delta_a^{\pm 1} D_\Delta$ and $\delta_b^{\pm 1} D_\Delta$.  

\begin{figure}[h] 
\includegraphics[scale=.8]{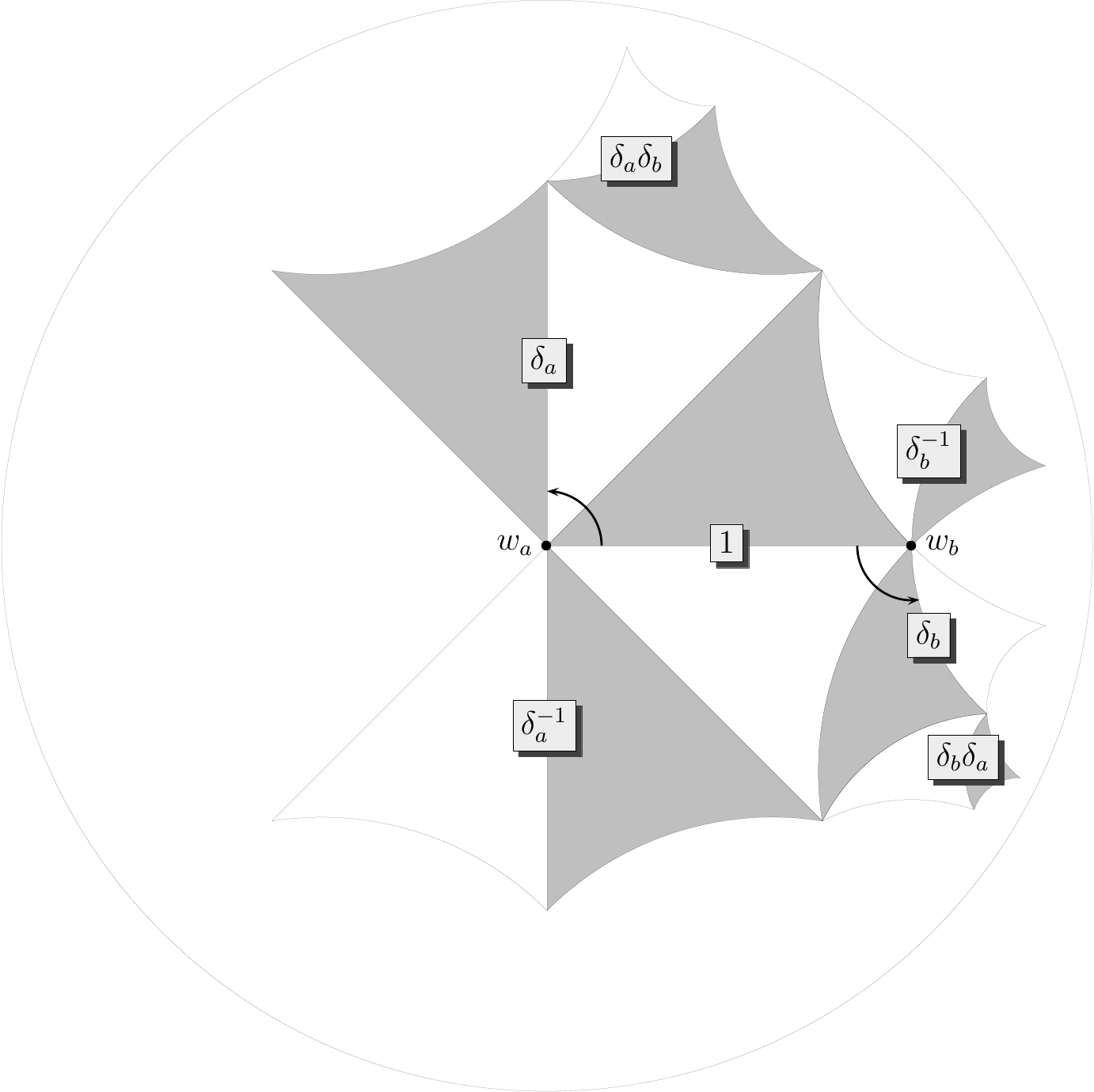}

\caption{The translates $\delta_a^{\pm 1}D_\Delta, \delta_b^{\pm 1}D_\Delta, \delta_a\delta_bD_\Delta, \delta_b\delta_aD_\Delta$.} \label{fig:translates}
\end{figure}

We check each of these cosets if they are old (in the start, equivalent to the identity): those that are not, we add to the \defi{frontier}, labelling them in the graph; those that are, we find an element in $\Gamma$ matching it with an existing coset.  We continue in this manner until the frontier is exhausted.

\begin{alg} \label{alg:cosetalg}
Let $\sigma \in S_d^3$ be a transitive permutation triple.  This algorithm returns a coset graph $G=G(\Gamma \backslash \Delta)$ for the group $\Gamma \leq \Delta$ associated to $\sigma$ and a side pairing for $G$.  

\begin{enumalg}
\item Initialize a graph $G$ with $d$ vertices $v_1,\dots,v_d$, label $v_1$ with $1 \in \Delta$.  Initialize the empty side pairing $S(G)$.
\item Let $v_j$ be the first vertex with no out edges.  If no such vertex exists, return $G,S(G)$.
\item For $\epsilon \in \{\delta_a,\delta_a^{-1},\delta_b,\delta_b^{-1}\}$:
\begin{enumalgalph}
\item Compute $i :=1^{\pi(\alpha_j\epsilon)}$.  Draw a directed edge from $v_j$ to $v_i$ with label $\epsilon$.
\item If vertex $v_i$ is labelled, let $\gamma := \alpha_j\epsilon\alpha_i^{-1}$.  If further $\gamma \neq 1$, then add $(\gamma,(j,\epsilon),(i,\epsilon^{-1}))$ to $S(G)$.  
\item If vertex $v_i$ is not labelled, label $v_i$ with $\alpha_i := \alpha_j \epsilon$.
\end{enumalgalph}
\item Return to Step 2.
\end{enumalg}
\end{alg}

In Step 3b, we only add the side pairing if $\gamma \neq 1$; this removes spurious (non-boundary) sides that are contained inside the fundamental domain $D_\Gamma$.  We can check if $\gamma=1$ exactly by computing the linear transformation corresponding to $\gamma$ as in Proposition \ref{prop:triangleembed} over the number field $K$ such that $\Delta(a,b,c) \hookrightarrow \PSL_2(K)$.  

\begin{rmk}
In fact, the $\Q$-algebra generated by any lift of the matrices $\delta_a,\delta_b,\delta_c$ to $\SL_2(\R)$ is a quaternion algebra over a smaller totally real number field \cite{ClarkVoight}, so the complexity of the check that $\gamma=1$ can be reduced somewhat.  Alternatively, since the group is discrete, we can check if $\gamma z = z$ for some point $z$ in the interior of $D_\Delta$ up to some precision depending on the group.
\end{rmk}

\begin{proof}[Proof of correctness of Algorithm \ref{alg:cosetalg}]
First, we show that the output graph $G$ of the algorithm is indeed a coset graph for $\Gamma$.  Let $k \in \{1,...,d\}$ and consider the vertex $v_k$ of $G$.  Since $\sigma$ is a transitive permutation triple, there exists $\delta \in \Delta$ such that $1^{\pi(\delta)} = k$.  By writing $\delta$ as a word in $\delta_a$ and $\delta_b$, we find that there will be a directed path from the vertex labeled $1$ to the vertex $v_k$.  Specifically, this directed path is given by following the edges labeled by elements $\delta_a$ and $\delta_b$ in the word equaling $\delta$, read from left to right.  Since an in-edge for $v_k$ is constructed by the algorithm, Step 2 of the algorithm will label the vertex $v_k$ with a specific coset representative.  

We have therefore shown that every vertex in $G$ is labeled with a coset representative and all of the cosets of $\Gamma$ in $\Delta$ are represented by exactly one vertex in $G$.   It remains to show that the edges in the graph $G$ that is returned are such that $G$ is a coset graph.  Consider the map from $G$ to the quotient graph of the Cayley graph of $\Delta$ acted on by $\Gamma$ that sends each vertex to the coset that the vertices label representatives and sends each labeled edge in $G$ to the edge with the same label in the quotient graph.  As each vertex is considered in Step 3 of the algorithm, all of the incident edges of that vertex are constructed such that the aforementioned map between graphs is a graph isomorphism.  Therefore the returned graph $G$ is in fact a coset graph.

We now show that the set $S(G)$ that is returned by the algorithm is in fact a side pairing for $G$.  First, assume that $(\gamma, (j, \epsilon), (i, \epsilon^{-1})) \in S(G)$ where $\gamma = \alpha_j \epsilon \alpha_i^{-1}$ is added in Step 3b of the algorithm.  Then $1^{\pi(\alpha_i \epsilon)} = i = 1^{\pi(\alpha_j)}$ so we have 
\[ (1^{\pi(\alpha_j \epsilon)})^{\pi(\alpha_i^{-1})} = 1^{\pi(\alpha_j \epsilon \alpha_i^{-1})} = 1\]
and therefore $\alpha_j \epsilon \alpha_i^{-1} = \gamma \in \Gamma$.  Consider such an element $\gamma \in \Gamma$.  We show that the edges $(j,\epsilon)$ and $(i, \epsilon^{-1})$ are sides of $D_\Gamma$ if and only if $\gamma \neq 1$.  Indeed if $\gamma \neq 1$ then the two edges are distinct subsets of $D_\Gamma$ that are in the same $\Gamma$-orbit.  Since $D_\Gamma$ is a fundamental domain for $\Gamma$, these two edges are on the boundary of $D_\Gamma$ and therefore are sides of $D_\Gamma$.  Conversely, if two distinct sides of $D_\Gamma$ are paired by an element in $\gamma \in \Gamma$, then we must have $\gamma \neq 1$. 

Lastly, we show that $S(G)$ induces a partition of $S$ into pairs.  Since the sides of $G$ are in bijection with a subset of the pairs $(j, \epsilon)$ with $j \in \{1,...,d \}$ and $\epsilon \in \{ \delta_a^{\pm 1}, \delta_b^{\pm 1} \}$ and all such pairs are considered in Step 3 of the algorithm, we know that each side will appear in a side pairing in the output $S(G)$.  Furthermore, each  side can be paired with at most one other side in the partition of $S$ that is induced by $S(G)$, since the value of $1^{\pi(\alpha_j \epsilon)}$ computed in Step 3a of the algorithm is unique.   
\end{proof}

\begin{lem}\label{lem:connected}\notag
Let $D_\Delta$ be a connected fundamental domain for $\Delta$ and let $D_\Gamma$ be the fundamental domain corresponding to a coset graph $G$ computed in Algorithm \ref{alg:cosetalg}.  Then $D_\Gamma$ is connected.
\end{lem}

The conclusion of this lemma is not true in general for an arbitrarily chosen coset graph (cf. Figure \ref{fig:EdgesSides}).

\begin{proof}
Let $\epsilon \in \{\delta_a^{\pm 1}, \delta_b^{\pm 1} \}$.  For any element $\delta \in \Delta$, note that $\delta(D_\Delta)$ and $\delta \epsilon (D_\Delta)$ share an edge and therefore $\delta(D_\Delta) \cup \delta \epsilon(D_\Delta)$ is connected.  For example, if $\epsilon = \delta_a$, then the side between $\delta(w_a)$ and $\delta(\overline{w_c})$ of $\delta(D_\Delta)$ is connected to the side between $\delta \epsilon (w_a)$ and $\delta \epsilon(w_c)$ of $\delta \epsilon(D_\Delta)$.  Note that for each adjacent pair of vertices in $G$, their representatives differ by multiplication on the left by some $\epsilon \in \{ \delta_a^{\pm 1}, \delta_b^{\pm 1} \}$. Therefore, the translates of $D_\Delta$ by the labels of adjacent vertices of $G$ are connected to one another.    

For a chosen coset representative $\alpha_i$ that is a label of a vertex of $G$, we showed that there is a directed path from the vertex labeled $1$ to the vertex labeled $\alpha_i$ in $G$.  By following along this path we are writing $\alpha_i$ as a word in $\delta_a$ and $\delta_b$ from left to right, and all of the translates of $D_\Delta$ by the coset representatives labeling the vertices along this path are a part of $D_\Gamma$.  Therefore, the union of the translates of $D_\Delta$ by all of the labels of the vertices in this path is a connected subset of $D_\Gamma$, that contains $\alpha_i(D_\Delta)$ and $D_\Delta$.  Hence, all $D_\Delta$ translates by the elements that are the labels of our coset graph $G$ are in the connected component of $D_\Delta$ in $D_\Gamma$.  Thus $D_\Gamma$ is connected.
\end{proof}

We propose also three potential variations on Algorithm \ref{alg:cosetalg}; these do not affect its correctness and can be used to obtain fundamental domains that are better (in some way).

First, in Step 2, we can instead choose a vertex $v_j$ with label $\alpha_j$ such that $d(\alpha_j0,0)$ is minimal, where $d$ is hyperbolic metric.  This can potentially yield a fundamental domain $D_\Gamma$ that is contained in a circle of smaller radius.

Second, we can use a \defi{petalling} approach in which we first attempt to extend the fundamental domain by translates consisting of powers of $\delta_a$, before considering those with powers of $\delta_b$. This results in the following modification of Algorithm \ref{alg:cosetalg}.  For $n \in \Z$, we define $\sgn(n)=+1,0,-1$ according as $n>0,n=0,n<0$.

\begin{alg} \label{alg:petal}
Let $\sigma \in S_d^3$ be a transitive permutation triple.  This algorithm returns a coset graph $G=G(\Gamma \backslash \Delta)$ for the group $\Gamma \leq \Delta$ associated to $\sigma$ and a side pairing for $G$.  
\begin{enumalg}
\item Initialize a graph $G$ with $d$ vertices $v_1,\dots,v_d$, label $v_1$ with $1 \in \Delta$.  Initialize the empty side pairing $S(G)$.
\item Let $v_j$ be the first vertex with fewer than $4$ out edges.  If no such vertex exists, return $G,S(G)$.
\item For $k=1,-1,2,-2,\dots$ until both $+1$ and $-1$ cases are done, do:
\begin{enumalgalph}
\item If the $\sgn(k)$ case is done, return to Step 3.  
\item Let $\epsilon := \delta_a^{\sgn(k)}$ and let $j' := j^{\sigma_a^{k-\sgn(k)}}$.
\item Let $i := 1^{\pi(\alpha_{j'} \epsilon)}$.  Draw a directed edge from $v_{j'}$ to $v_i$ with label $\epsilon$.
\item If vertex $v_i$ is labelled, let $\gamma := \alpha_{j'}\epsilon\alpha_i^{-1}$, and declare that the $\sgn(k)$ case is done.  If $\gamma \neq 1$, then add $(\gamma,(j',\epsilon),(i,\epsilon^{-1}))$ to $S(G)$. 
\item If vertex $v_i$ is not labelled, label $v_i$ with $\alpha_i := \alpha_{j'} \epsilon$.
\end{enumalgalph}
\item For $\epsilon\in\{\delta_b,\delta_b^{-1}\}$, repeat Steps 3c--3e with $j' := j$.
\item Return to Step 2.
\end{enumalg}
\end{alg}

In Step 3, we alternate signs in order to make the nicest pictures; it is equivalent to simply take $k=1,2,\dots$.  Since the coset representatives labelling adjacent vertices in the coset graph produced by this approach still differ by multiplication on the right by $\epsilon\in\{\delta_a^{\pm 1},\delta_b^{\pm 1}\}$, the argument in the proof of Algorithm \ref{alg:cosetalg} applies here. By the same reasoning, the translates of $D_\Delta$ by the labels of adjacent vertices of $G$ are connected and Lemma \ref{lem:connected} also applies. 

\begin{exm}
In this example, we compare the results of Algorithms (\ref{alg:cosetalg}) and (\ref{alg:petal}).   Let $\Gamma\leq \Delta(6,6,7)$ corresponding to the permutations 
\[ \sigma_0 = (1\;3\;2\;6\;4\;5),\ \sigma_1 = (1\;4\;5\;3\;7\;6),\ \sigma_\infty = (1\;4\;7\;3\;6\;2\;5) . \]  The advantage to preferring powers of $\delta_a$ is the resulting ``petalling" effect illustrated in Figure \ref{fig:CompareAlg6,6,7}.
\end{exm}

\begin{figure}[h] 
\includegraphics{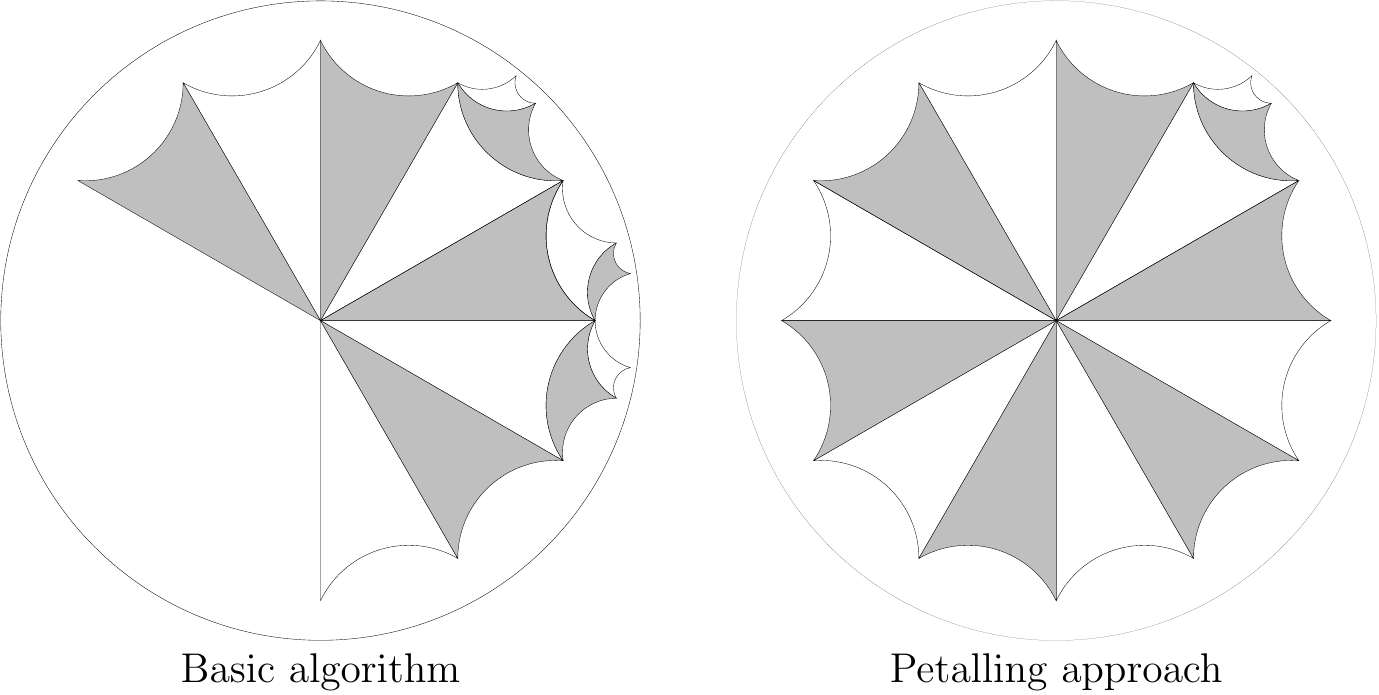}

\caption{Comparing the petalling approach} \label{fig:CompareAlg6,6,7}
\end{figure}

A third possible variation is to find the ``smallest'' cosets in each case.  More precisely, once we have computed a set of generators for $\Gamma$ coming from side pairing elements, we utilize the algorithm of Voight \cite{Voightfd} (the original idea is due to Ford) to compute the \defi{Dirichlet domain} for $\Gamma$, the set
\[ \{w \in \calD : d(0,\gamma w) \geq d(0,w) \text{ for all $\gamma \in \Gamma$}\} \]
that picks out in each $\Gamma$ orbit the points closest to the origin.  From there, we choose a point $p$ in the fundamental triangle $D_\Delta$ and, for each coset of $\Gamma$ in $\Delta$, we choose a representative $\alpha$ such that $\alpha p$ lies in the Dirichlet domain.  In this way, we can compute a fundamental domain which has the property that it lies in a circle of radius $\rho$ with $\rho$ minimal among all possible fundamental domains that are triangulated by $D_\Delta$.  But we choose not to employ this modification in practice for the following reasons.  First, we find the petalling approach to be more natural for our calculations.  Below, we compute power series centered at the elements of the orbit of $0$ and so it is enough to minimize the distance between an arbitrary in the fundamental domain to one of these centers.  Second, we found in practice that the petalling approach produced more aesthetically pleasing dessins.  Finally, it can be a bit expensive to compute the Dirichlet domain.  That being said, it is possible that the Dirichlet domain may be useful in other variants of our method, so we mention it here.

\begin{exm}
Below in Figure \ref{fig:CompareAlgorithms} is a simple example comparing the variations of Algorithm \ref{alg:cosetalg}. The basic algorithm and petalling approach yield a fundamental domain inside the solid blue circle. The Dirichlet domain approach yields a fundamental domain inside the dashed red circle. When we compare these circles we see that the radius of the dashed red circle is smaller than the radius of the solid blue circle. In Figure \ref{fig:CompareAlgorithms}, we have three fundamental domains for $\Gamma \leq \Delta(3,5,3)$ corresponding to the permutations 
\[ \sigma_0 = (1\;3\;4)(2\;6\;5),\ \sigma_1 = (1\;4\;2\;5\;6),\ \sigma_\infty = (3\;6\;4). \]
\end{exm}

\begin{figure}[h] 
\includegraphics[scale=1]{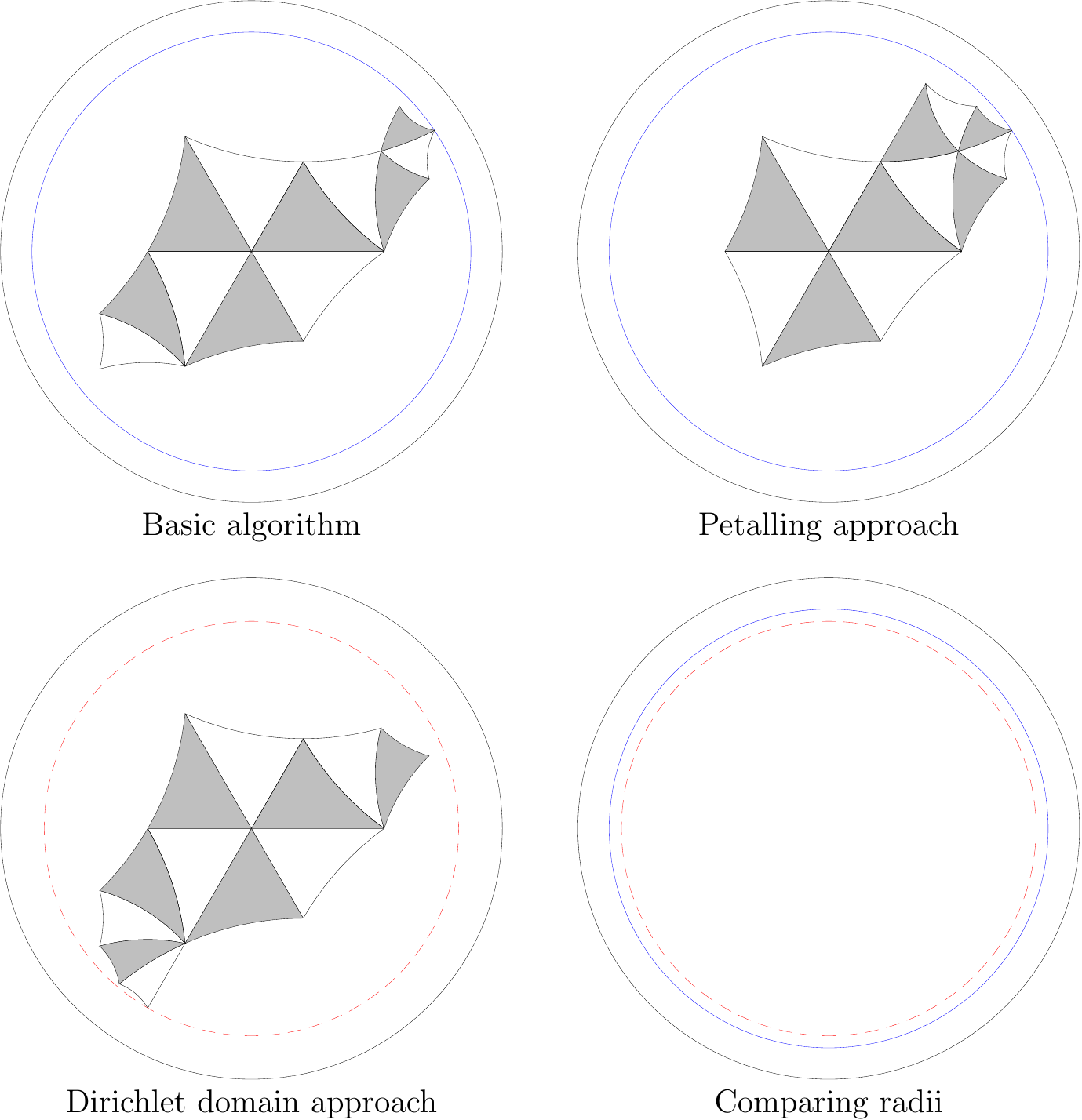}

\caption{Comparison of variations of Algorithm \ref{alg:cosetalg}.} \label{fig:CompareAlgorithms}
\end{figure}

\subsection*{Drawing dessins}\label{subsec:drawdessin}

We now show how the preceding methods fit together in our first main result, to conformally draw dessins.  Compare the work of Schneps \cite[\S III.1]{Schneps} for genus zero dessins (given by an explicit rational function).  Rather than getting bogged down in algorithmic generalities, it is clearer instead to illustrate with an extended example.  

Let 
\[ \sigma_0 =(1\;3\;5\;4)(2\;6),\ \sigma_1 =(1\;5\;2\;3\;6\;4),\ \sigma_\infty = (1\;6\;5\;2\;3\;4) \]
and let $\Gamma \leq \Delta(4,6,6)$ be the corresponding subgroup.  Let $D_\Delta$ be the fundamental domain for $\Delta$ as in Proposition \ref{prop:triangleembed}.  

\begin{figure}[h] 
\includegraphics[scale=.8]{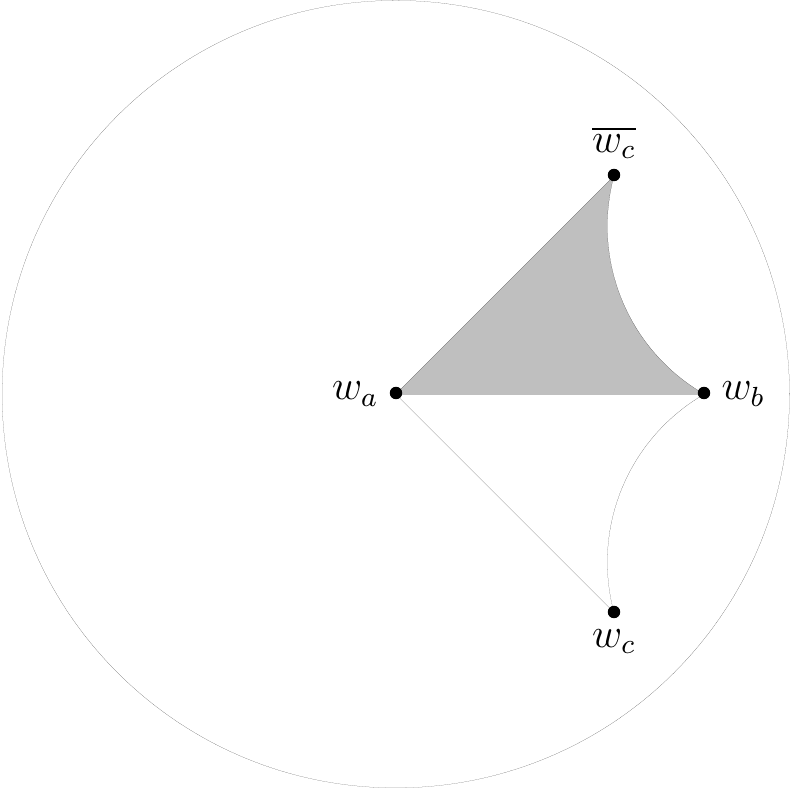}

\caption{$D_\Delta$ for $\Delta = \Delta(4,6,6)$ as in Proposition \ref{prop:triangleembed}.} \label{fig:dessinpic1} 
\end{figure}

\begin{figure}[h] 
\includegraphics[scale=.8]{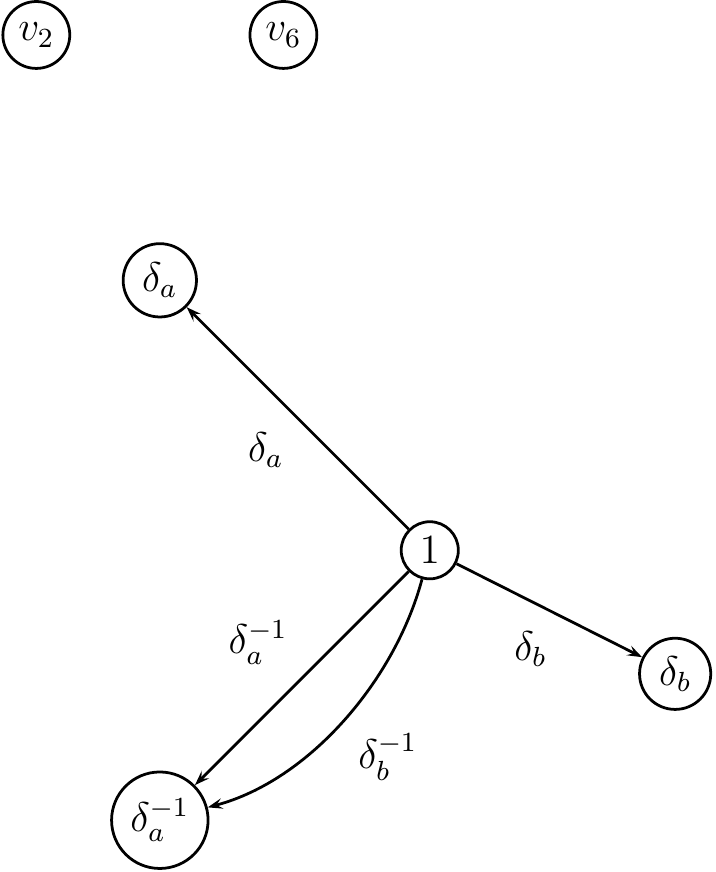}

\caption{A partial coset graph and side pairing at the end of the first iteration of Step $3$ in Algorithm \ref{alg:cosetalg} for our example.} \label{fig:dessinpic2a}
\end{figure}

We first compute a coset graph $G = G(\Gamma \backslash \Delta)$ and a side pairing $S = S(G)$ for $\Gamma$ using Algorithm \ref{alg:cosetalg}. In the first iteration of Step $3$ we have $j=1$ and $\alpha_1=1$. After computing $1^{\pi(\alpha_1\epsilon)}$ for $\epsilon\in\{\delta_a,\delta_a^{-1},\delta_b,\delta_b^{-1}\}$ we have the partial coset graph and side pairing as shown in Figure \ref{fig:dessinpic2a} at the end of the first iteration of Step $3$. Explicitly, since $1^{\pi(\alpha_1\delta_a)}= 3$ we add an edge from $v_1$ to $v_3$ labeled $\delta_a$ and we label $v_3$ with $\alpha_1\delta_a = \delta_a$. Similarly, since $1^{\pi(\alpha_1\delta_a^{-1})}=4$ and $1^{\pi(\alpha_1\delta_b)}=5$, we add edges from $v_1$ to $v_4$ and $v_1$ to $v_5$ labeled $\delta_a^{-1}$ and $\delta_b$ respectively, and we label $v_4$ and $v_5$ with $\delta_a^{-1}$ and $\delta_b$ respectively. Now for $\epsilon = \delta_b^{-1}$ we see that $1^{\pi(\alpha_1\delta_b^{-1})}= 4$. Since $v_4$ is already labeled $\delta_a^{-1}$, we add the edge from $v_1$ to $v_4$ labeled $\delta_b^{-1}$, compute $\gamma = \alpha_1\epsilon\alpha_4^{-1} = \delta_b^{-1}\delta_a$, verify that the linear transformation corresponding to $\delta_b^{-1}\delta_a$ as in Proposition \ref{prop:triangleembed} is not the identity matrix, and append $(\gamma,(j,\epsilon),(i,\epsilon^{-1})) = (\delta_b^{-1}\delta_a,(1,\delta_b^{-1}),(4,\delta_b))$ to the side pairing $S$.

\FloatBarrier

\begin{figure}[h] 
\includegraphics[scale=.8]{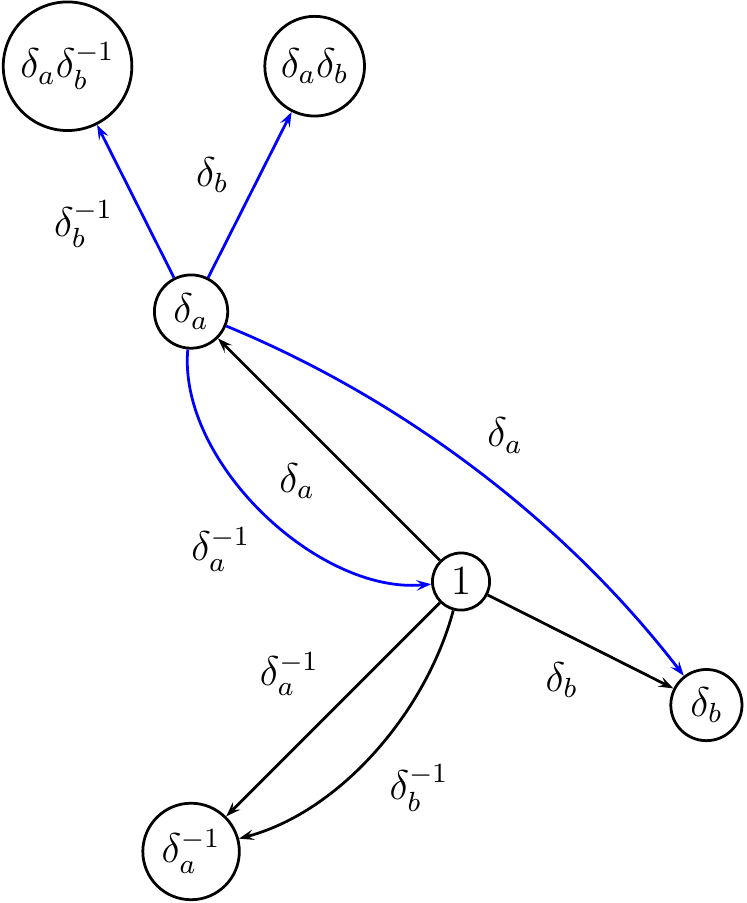}

\begin{align*}
	S = \{(\delta_b^{-1}\delta_a,(1,\delta_b^{-1}),(4,\delta_b)),(\delta_a^2,(3,\delta_a),(5,\delta_a^{-1})) \}
\end{align*}
\caption{A partial coset graph and side pairing at the end of the second iteration of Step $3$ in Algorithm \ref{alg:cosetalg} for our example.} \label{fig:dessinpic2b} 
\end{figure}

By convention, we prefer earlier labels appearing in the coset graph and therefore choose $\delta_a$ in the next iteration of Step $3$, and take $j=3$ and $\alpha_3 = \delta_a$. After computing $1^{\pi(\alpha_3\epsilon)}$ for $\epsilon\in\{\delta_a,\delta_a^{-1},\delta_b,\delta_b^{-1}\}$ we have the following partial coset graph and side pairing in Figure \ref{fig:dessinpic2b} at the end of the second iteration of Step $3$. Explicitly, since $1^{\pi(\delta_a\delta_b)}=6$ and $1^{\pi(\delta_a\delta_b^{-1})}=2$, we label $v_6$ and $v_2$ with $\delta_a\delta_b$ and $\delta_a\delta_b^{-1}$ respectively. Since $1^{\pi(\delta_a^2)}=5$ we append $(\delta_a^2\delta_b^{-1},(3,\delta_a),(5,\delta_a^{-1}))$ to $S$. For $\epsilon=\delta_a^{-1}$ we find that the resulting side pairing element $\gamma$ corresponds to the identity matrix and therefore we do not add this side pairing element to $S$. The new edges added in this iteration are colored blue.

\FloatBarrier

\begin{figure}[h] 
\includegraphics[scale=.8]{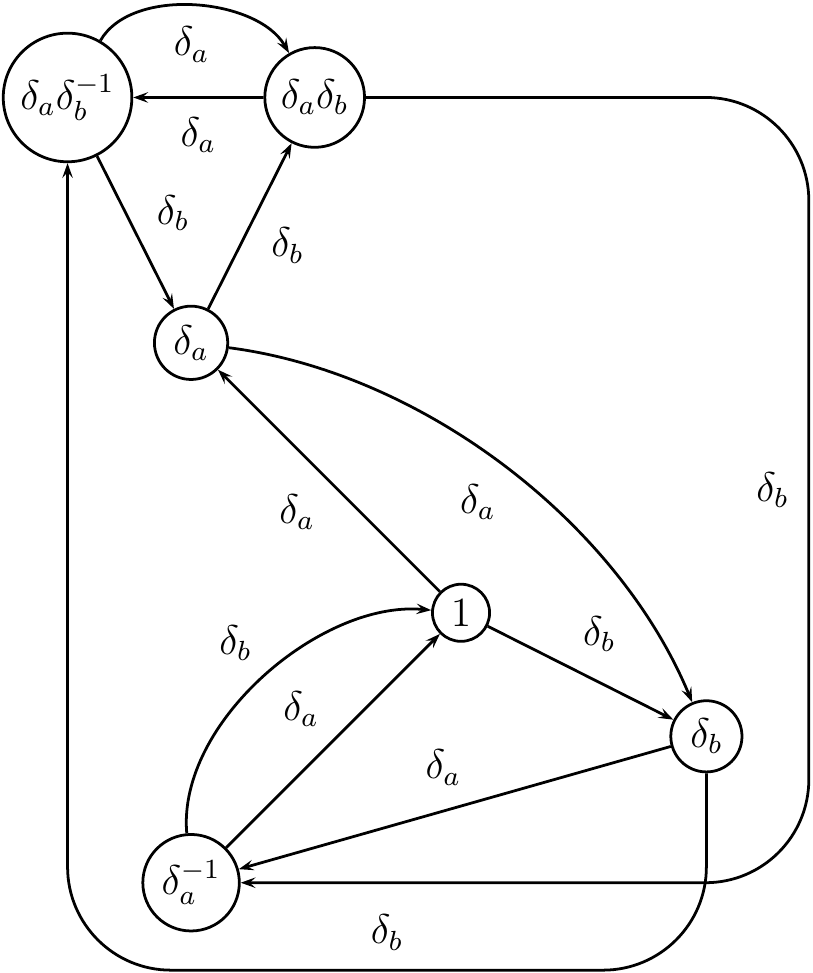}

\[ S = \{(\delta_b^{-1}\delta_a,(1,\delta_b^{-1}),(4,\delta_b)),(\delta_a^2\delta_b^{-1},(3,\delta_a),(5,\delta_a^{-1})),\dots,(\delta_a\delta_b^2\delta_a,(6,\delta_b),(4,\delta_b^{-1}))\} \]

\caption{Output of Algorithm \ref{alg:cosetalg}.} \label{fig:dessinpic2c}
\end{figure}

We continue according to Algorithm \ref{alg:cosetalg} taking the labels $\delta_a^{-1}, \delta_b, \delta_a\delta_b,\delta_a\delta_b^{-1}$ in order for the four remaining iterations of Step $3$. At this point the algorithm terminates and returns $G$ and $S$ shown in Figure \ref{fig:dessinpic2c}. In Figure \ref{fig:dessinpic2c} we have the output of Algorithm \ref{alg:cosetalg}. We have only included the edges of the coset graph labeled $\delta_a$ and $\delta_b$ since for every edge labeled $\delta_a$ or $\delta_b$ there is a corresponding reverse edge labeled $\delta_a^{-1}$ or $\delta_b^{-1}$ respectively.

Next we draw all of the translates of $D_\Delta$ by the coset representatives of $\Gamma \backslash \Delta$ given by the vertices of $G$ as in \eqref{eqn:fundGamma}.  For example, for the vertex of $G$ labeled $\delta = \delta_a\delta_b^{-1}$ we compute the matrix corresponding to $\delta$ as a matrix as in Proposition  \ref{prop:triangleembed}, namely
$$
\delta = 
\begin{pmatrix}
0.6552\dots & -2.3015\dots\\
-0.5694\dots & 3.5262\dots
\end{pmatrix}.
$$
Next we compute $\delta(w_a), \delta(w_b), \delta(w_c), \delta(\overline{w_c})$ and draw in the corresponding translates of the geodesics in $D_\Delta$ between $w_a,w_b,w_c,\overline{w_c}$ as in Figure \ref{fig:dessinpic1}.  By doing this for all of the vertices of $G$ we obtain a fundamental domain for $\Gamma$.  Given a vertex of $G$ labeled $\delta$ such that $1^{\pi(\delta)} = j$ we label the translate of $D_\Delta$ with $j$ and explicitly give the element $\delta$ as a word in $\delta_a^{\pm 1},\delta_b^{\pm 1}$.  For example, again taking $\delta = \delta_a\delta_b^{-1}$ we have $1^{\pi(\delta)} = 2$ so we label $\delta D_\Delta$ with $2$.  

\begin{figure}[h]
\includegraphics{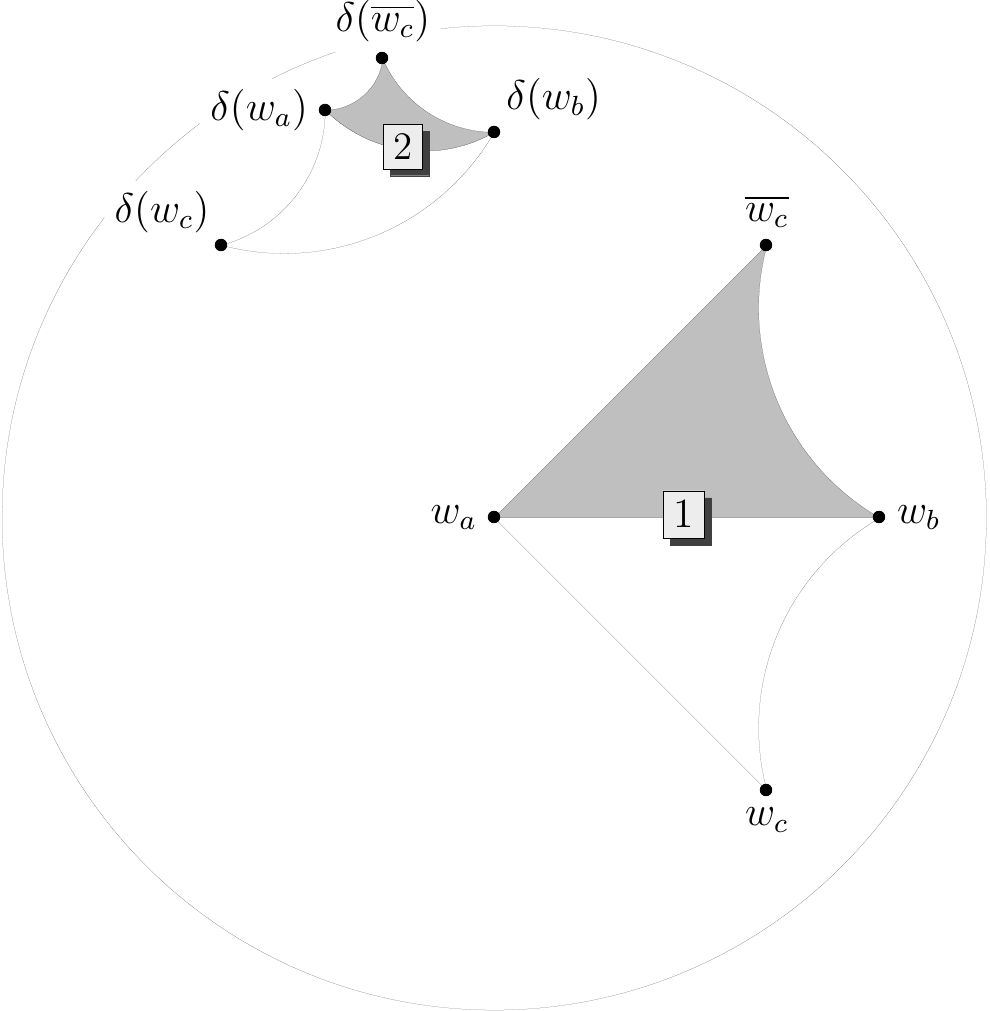}

\caption{$D_\Delta$ and the translate $\delta D_\Delta$ where $\delta = \delta_a\delta_b^{-1}$.}  \label{fig:dessinpic3} 
\end{figure}

\FloatBarrier

Next we must consider the side pairings of $G$ and make them visible on $D_\Gamma$.  For each $(\gamma, (i, \epsilon_1), (j, \epsilon_2)) \in S$ we locate the sides $(i, \epsilon_1)$ and $(j, \epsilon_2)$ of $D_\Gamma$ and give them the same label.  We signify that $\gamma$ maps the side $(j, \epsilon_2)$ to the side $(i, \epsilon_1)$ by coloring the label for side $(j,\epsilon_2)$ red and the label for side $(i,\epsilon_1)$ blue.  For example, in Figure \ref{fig:dessinpic2c} we see that $(\delta_a^2\delta_b^{-1},(3,\delta_a),(5,\delta_a^{-1})) \in S$ so the side $(3,\delta_a)$ which is $\delta_aD_\Delta \cap \delta_a^2D_\Delta$ is identified with the side $(5,\delta_a^{-1})$ which is $\delta_bD_\Delta \cap \delta_b\delta_aD_\Delta$ by the element $\delta_a^2\delta_b^{-1}$.

\begin{figure}[h] 
\includegraphics[scale=.8]{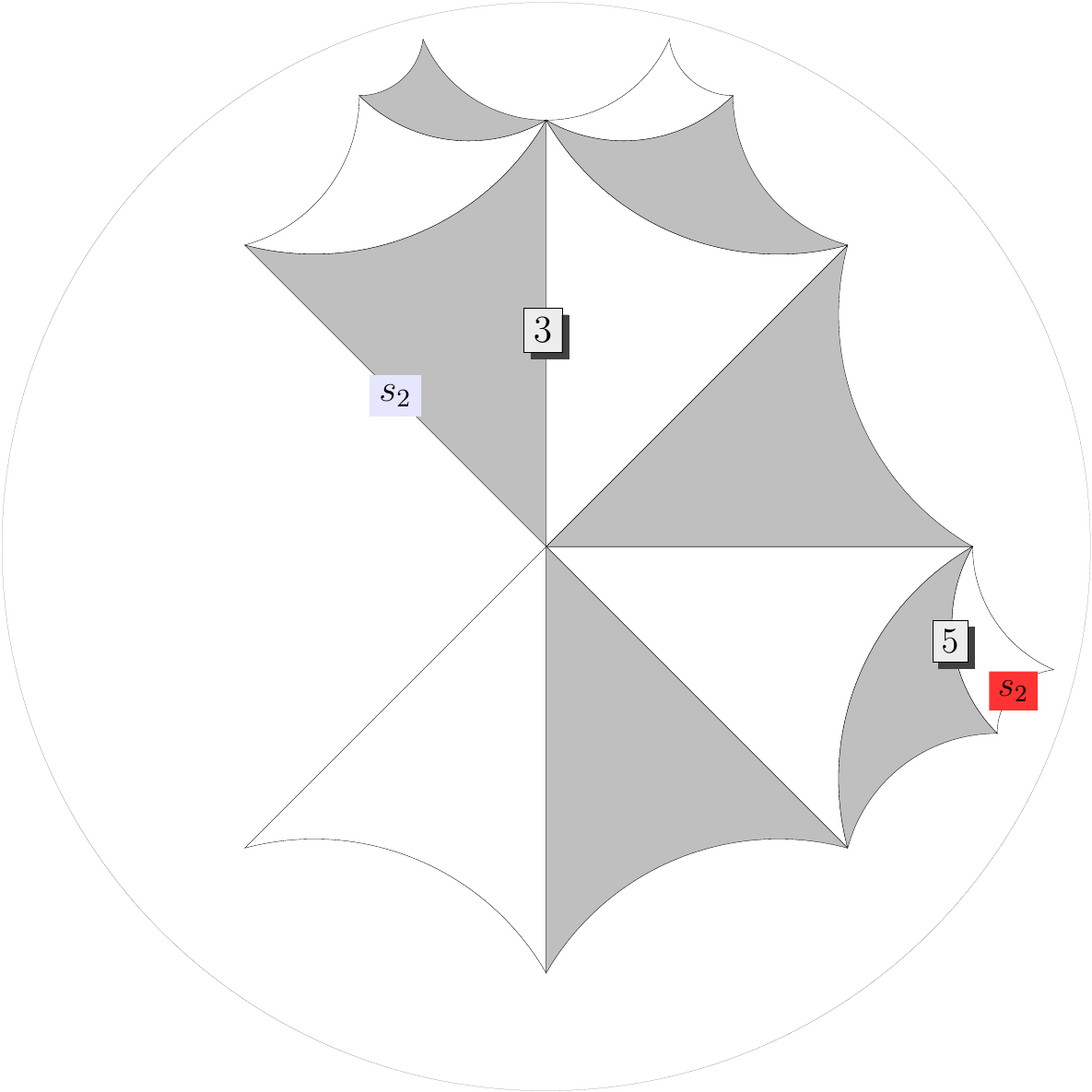}

\caption{An illustration of our convention for labeling the elements of $S$.} \label{fig:dessinpic4}
\end{figure}

In Figure \ref{fig:dessinpic4} we describe our convention for labeling the elements of $S$. For example, the element $(\delta_a^2\delta_b^{-1},(3,\delta_a),(5,\delta_a^{-1}))$ in $S$ tells us that $\gamma=\delta_a^2\delta_b^{-1}$ maps the side $(5,\delta_a^{-1})$ to the side $(3,\delta_a)$. This is indicated in the picture by labeling the two sides $s_2$ (subscript $2$ since we are considering the second element in $S$) and coloring the side $(5,\delta_a^{-1})$ red and the side $(3,\delta_a)$ light blue, i.e.:
\begin{center}
\includegraphics{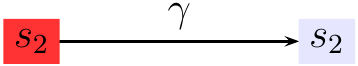}

\end{center}
Also note that $S$ contains the element $(\delta_b\delta_a^{-2},(5,\delta_a^{-1}),(3,\delta_a))$ which is just the same side pairing in the opposite direction.

\FloatBarrier

Finally, for each translate $\delta D_\Delta$ we draw a white dot at $\delta(w_a)$, a black dot at $\delta(w_b)$, an $\times$ at $\delta(w_c)$, and bold the geodesic from $\delta(w_a)$ to $\delta(w_b)$. The white and black dots are the vertices of the dessin, and the bold geodesics are the edges of the dessin. By identifying the sides that are paired together we form a quotient Riemann surface with the dessin conformally embedded.

\begin{figure}[h] 
\includegraphics{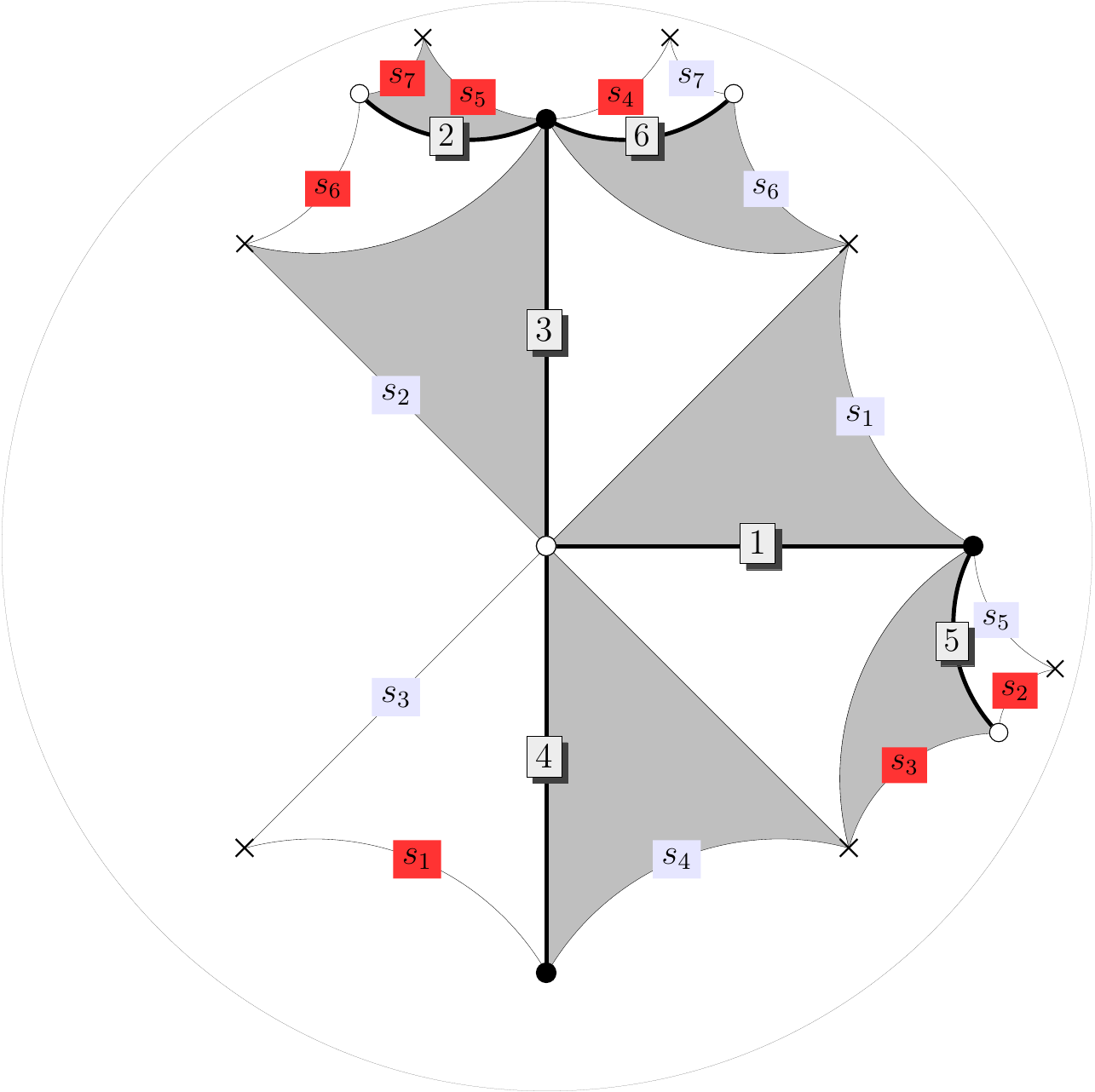}

\caption{A fundamental domain for $\Gamma$ with the embedded dessin.} \label{fig:dessinpic5} 
\end{figure}

In Figure \ref{fig:dessinpic5} is the final fundamental domain for $\Gamma$ which, after identifying the sides, yields a Riemann surface with the dessin conformally embedded. Notice that some of the white and black vertices are identified in the quotient and we are left with $2$ white vertices and $1$ black vertex. From this we can read off each cycle in $\sigma_0$ and $\sigma_1$ by rotating counterclockwise around each of the white dots and black dots respectively.

\FloatBarrier

\subsection*{Reduction algorithm}

We now describe the reduction algorithm for $\Delta$ that takes a point $w \in \calD$ and returns a point $w' = \delta w \in D_\Delta$ and an element $\delta \in \Delta$.  This algorithm \cite{VoightThesis} generalizes the classical reduction theory for $\SL_2(\Z)$ (see e.g.\ Serre \cite[Theorem 7.1]{SerreCourse}); it is a special (and very efficient) case of a full algorithmic theory of Dirichlet fundamental domains for Fuchsian groups whose quotients have finite area \cite{Voightfd}.

\begin{alg} \label{alg:fdreddelta}
Let $z \in \calH$ and let $\Delta=\Delta(a,b,c)$ be a triangle group with fundamental domain $D_\Delta$ given as in Proposition \ref{prop:triangleembed}.  This algorithm returns an element $\delta \in \Delta$ as a word in $\delta_a,\delta_b,\delta_c$ such that $z'=\delta z \in D_\Delta$.  

\begin{enumalg}
\item Initialize $\delta := 1$, $w := w_{z_a}(z)$.
\item Let 
\[ \alpha :=\arg(w) \in [0, 2 \pi). \]  
Let 
\[ i := - \left \lfloor \frac{a \alpha}{2 \pi} + \frac{1}{2} \right \rfloor. \]  
Let $w := \delta_a^i w$ and $\delta := \delta_a^i \delta$.
\item Let
\begin{equation} \label{Amat}
A = \begin{pmatrix} \newt+1 & -\newt+1 \\ -\newt+1 & \newt+1  \end{pmatrix}
\end{equation}
with $\mu$ as in Proposition \ref{prop:triangleembed}.  Let 
\[ \beta :=\arg(-A w_{a}) \in [0, 2 \pi). \] 

Let 
\[ j := -\left \lfloor \frac{b \beta}{2 \pi} + \frac{1}{2} \right \rfloor. \]  
Let $w := \delta_b^jw$ and $\delta := \delta_b^j \delta$.  If $j = 0$, then return $\delta$, otherwise return to Step 2.   
\end{enumalg}
\end{alg}

\begin{proof}
Let $w = w_{z_a}(z)$.  By identifying $\calH$ and $\calD$ with the map $w_{z_a}$, the element $\delta_a$ acts on $\calD$ by a rotation by $2\pi/a$ about the origin.  We can therefore rotate $w$ by a unique element $\delta_a^i$ with $-a < i \leq 0$ so that it is in the region of $\calD$ bounded by the images of the geodesics $z_a z_c$ and $z_a z_c'$.  In fact $i = - \left \lfloor \frac{a \alpha}{2 \pi} + \frac{1}{2} \right \rfloor$ is the unique such integer.  We now show that by applying Step $3$ to the resulting point, when nontrivial we get a point with smaller absolute value.  Since the orbits of the action of $\Delta$ on $\calH$ are discrete, the algorithm will terminate after finitely many steps with a point in $D_\Delta$.

Let $w$ be a point in the region bounded by the images of the geodesics $z_a z_c$ and $z_a z_c'$.  There exists a unique isometry of the unit disc that maps $w_{z_a}(z_b)$ to $0$ and $w_{z_a}(z_a)$ to the negative real axis: this isometry is given by the linear fraction transformation corresponding to $A$ in \eqref{Amat}.  We can therefore rotate $w$ by a unique element $\delta_b^j$ with $0 \leq j < b$ so that it is in the region of $\calD$ bounded by the goedesics $z_b z_c$ and $z_b z_c'$.  In fact $j = \left \lfloor \frac{b \beta}{2 \pi} + \frac{1}{2} \right \rfloor$ is the unique such integer.  Then $\delta_b^j Aw $ will be closer to $A w_{z_a}(z_a)$ than $A w$.  Then, $\delta_b^j w$ has smaller absolute value than $w$ when $j \neq 0$, and the proof is complete.
\end{proof}

\begin{figure}[h] 
\includegraphics{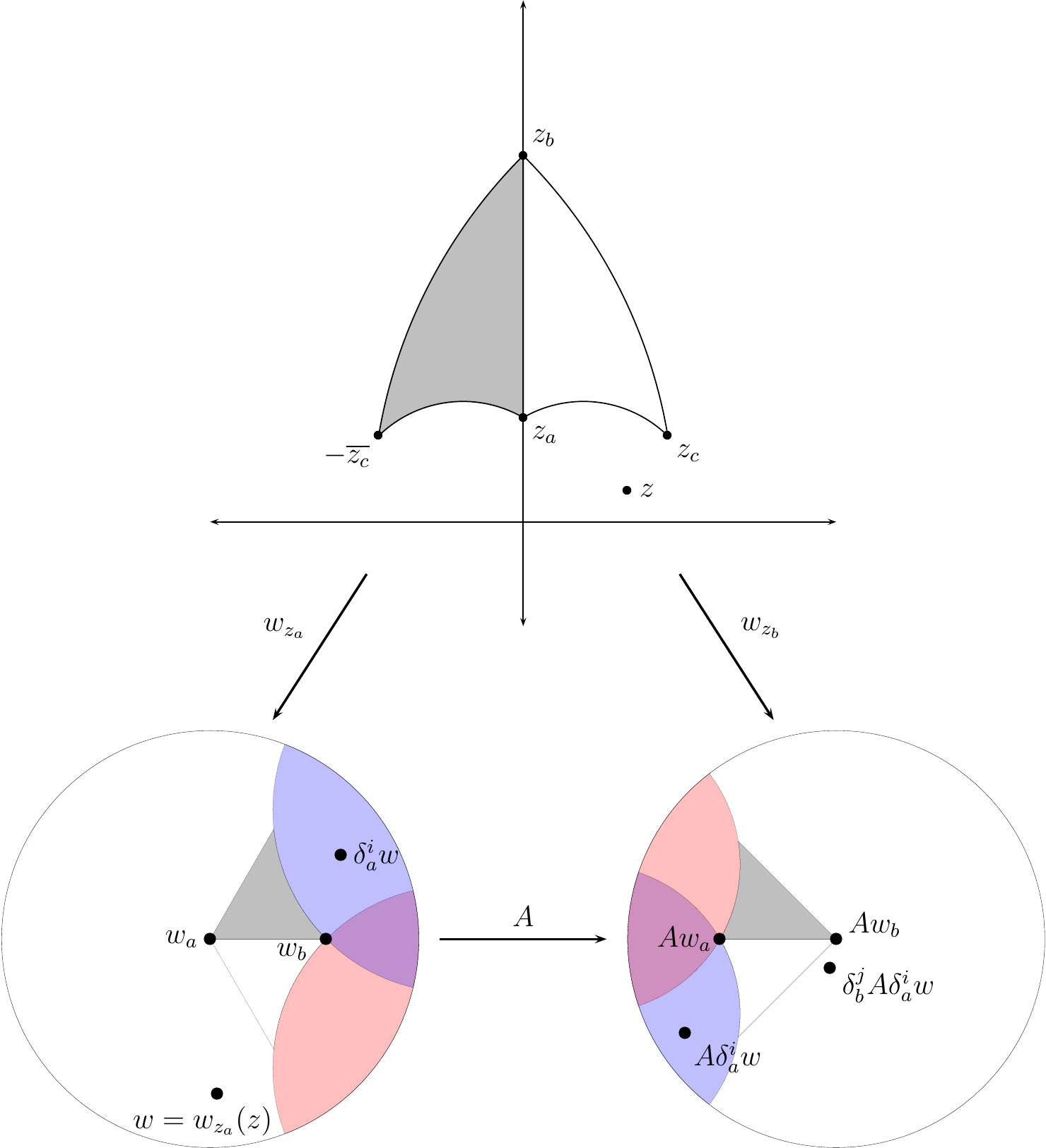}

\caption{An example of the reduction algorithm.  Points in the blue region are moved closer to $D_\Delta$ by $\delta_b$, while points in the red region are moved closer to $D_\Delta$ by $\delta_b^{-1}$.} \label{fig:reduction2} 
\end{figure}

\begin{rmk}
Due to roundoff issues, for points near the boundary of the fundamental triangle, one either needs to pick a preferred side among the possible two paired sides or to slightly thicken the fundamental domain (dependent on the precision) to avoid bouncing between points equivalent under $\Delta$.  For the purposes of our numerical algorithms, any such choice works as well as any other.
\end{rmk}

We now describe how a reduction algorithm for $\Gamma$ is obtained from the reduction algorithm for $\Delta$ together with Algorithm \ref{alg:cosetalg}.

\begin{alg} \label{alg:fdredgamma}
Let $\Gamma \leq \Delta(a,b,c)$ be a subgroup of index $d=[\Delta:\Gamma]$.  Let $\Gamma \alpha_1 , \dots , \Gamma \alpha_d$ be a set of right cosets for $\Gamma$ in $\Delta$ and let $D_\Gamma = \alpha_1 D_\Delta \cup \cdots \cup \alpha_d D_\Delta$.  This algorithm returns an element $\gamma \in \Gamma$ such that $z'=\gamma z \in D_\Gamma$. 

\begin{enumalg}
\item Using Algorithm \ref{alg:fdreddelta}, let $\delta \in \Delta$ be such that $\delta z \in D_\Delta$.  
\item Compute $i := 1^{\pi(\delta^{-1})}$ and return $\gamma=\alpha_i \delta$.
\end{enumalg}
\end{alg}

\begin{proof}[Proof of correctness]
We have $1^{\pi(\delta^{-1})}=i$ if and only if $\Gamma \delta^{-1} = \Gamma \alpha_i$, so we have $\gamma=\alpha_i \delta \in \Gamma$.  But since $\delta z \in D_\Delta$ we have $\gamma z = \alpha_i (\delta z) \in \alpha_i D_\Delta \in D_\Gamma$, as claimed.
\end{proof}

\begin{rmk}
The reason we need to take an inverse in Step 2 of Algorithm \ref{alg:fdredgamma} is because the cosets are labelled by element moving from the center $0$ to the translated piece, but we want to move into the fundamental domain, so we need an inverse.  If the tessellates are labelled by the element that maps to the fundamental triangle, the inverse would be replaced somewhere else (adjacent tessellates would be inverted).
\end{rmk}

\begin{exm}
We illustrate the reduction algorithms (Algorithms \ref{alg:fdreddelta} and \ref{alg:fdredgamma}) with an example. Let $\Delta = \Delta(3,4,5)$ and define $\pi : \Delta \to S_6$ by $\delta_a, \delta_b, \delta_c \mapsto \sigma_0, \sigma_1, \sigma_\infty$, where 
\[ \sigma_0 = (1\;2\;3),\ \sigma_1 = (1\;2)(3\;4\;5\;6),\ \sigma_\infty = (1\;6\;5\;4\;3). \] 
Given the example point $w = 0.34992\ldots + 0.82246\ldots i$ 
, the reduction algorithm for $\Delta$ returns the point $\delta w = 0.08700\ldots + 0.09353\ldots i \in D_\Delta$, where $\delta = \delta_a^{-1} \delta_b^{-1} \delta_a^{-1} \delta_b^{-1} \delta_a^{-1}$.  This $\delta$ corresponds to the linear fractional transformation represented by the matrix
$$
\begin{pmatrix}
1.9532\dots & 2.6172\dots\\
1.0755\dots & 1.9531\dots
\end{pmatrix}.
$$
(See Figure \ref{fig:reduction3}.)

In the reduction algorithm for $\Gamma$, we compute $\pi(\delta^{-1}) = (1\;2\;5)(3\;4\;6)$, so $1^{\pi(\delta^{-1})} = 2$.  This means that $\gamma w$ will be in the coset labeled 2.  From the table in Figure \ref{fig:reduction3}, we see that $\alpha_2 = \delta_a$, i.e.~$\delta_a$ is a representative for the coset labeled 2.  Then we have:
\begin{align*}
\gamma = \alpha_2 \delta = \delta_a  \delta_a^{-1}  \delta_b^{-1}  \delta_a^{-1}  \delta_b^{-1}  \delta_a^{-1} = (\delta_b^{-1} \delta_a^{-1})^2 
\end{align*}
This $\gamma$ corresponds to the linear fractional transformation represented by the matrix
$$
\begin{pmatrix}
1.9080\dots & 3.0001\dots\\
-1.1537\dots & -1.2900\dots
\end{pmatrix}.
$$
Applying $\gamma$ to $w$ yields the point $\gamma w = -0.12450\dots + 0.02857\dots i \in \alpha_2 D_\Delta \in D_\Gamma$. 

\begin{figure}[h] 
\includegraphics{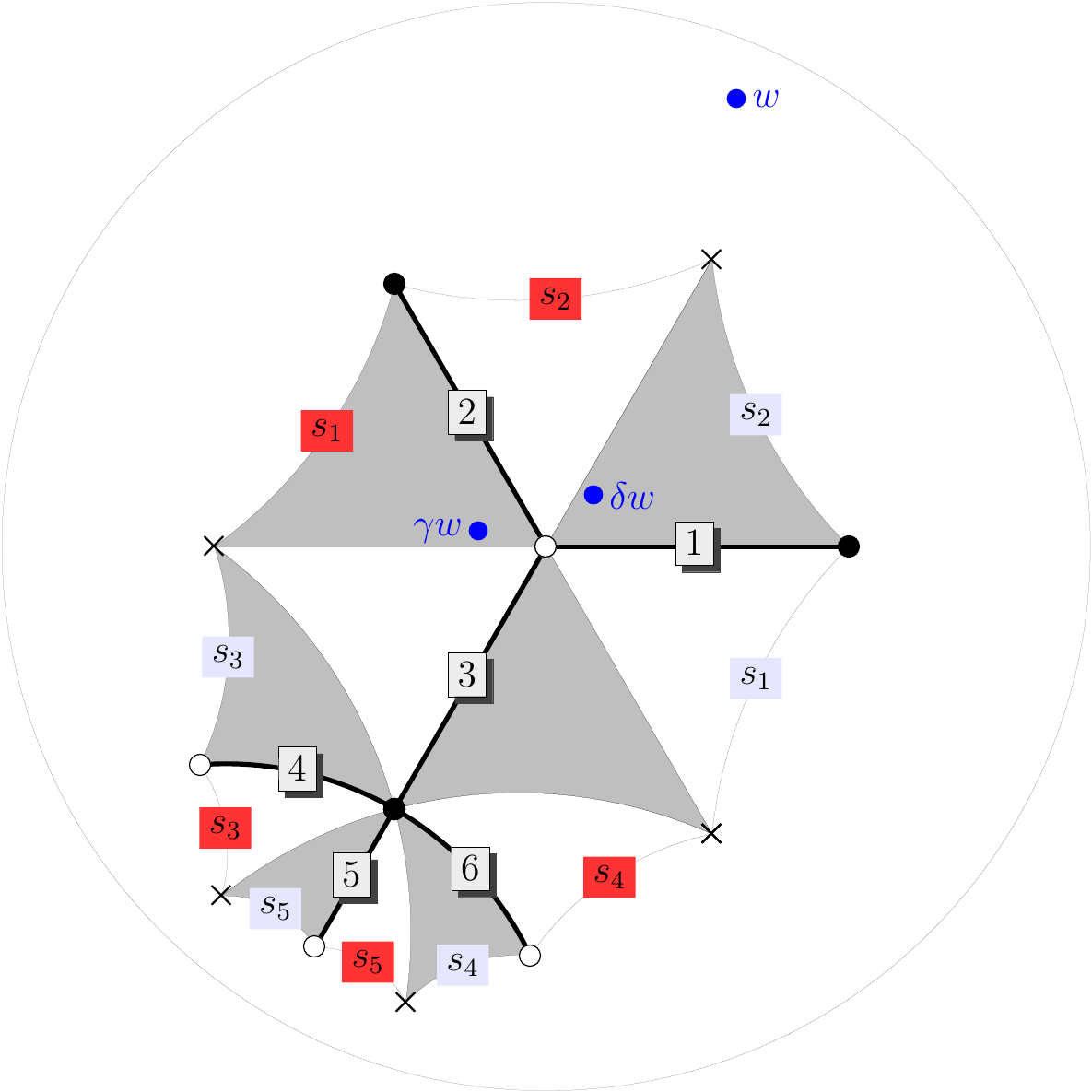}

\begin{center}
\begin{tabular}{ll}
\toprule
Coset Number & Label\\
\midrule
$1$ & 1 \\
$2$ & $\delta_a$ \\
$3$ & $\delta_a^{-1}$ \\
$4$ & $\delta_a^{-1 }  \delta_b$ \\
$5$ & $\delta_a^{-1 }  \delta_b^{-2}$ \\
$6$ & $\delta_a^{-1 }  \delta_b^{-1}$ \\
\bottomrule
\end{tabular}
\hspace{2cm}
\begin{tabular}{ll}
\toprule
Side Pairing & Label\\
\midrule
$s_{1}$ & $\delta_b   \delta_a^{-1}$ \\
$s_{2}$ & $\delta_b^{-1 }  \delta_a^{-1}$ \\
$s_{3}$ & $\delta_a^{-1 }  \delta_b^{-1 }  \delta_a   \delta_b   \delta_a$ \\
$s_{4}$ & $\delta_a^{-1 }  \delta_b   \delta_a   \delta_b^{-1 }  \delta_a$ \\
$s_{5}$ & $\delta_a^{-1 }  \delta_b^{-2 }  \delta_a   \delta_b^{2 }  \delta_a$ \\
\bottomrule
\end{tabular}
\end{center}
\bigskip

\caption{For the example point $w = 0.34992\ldots+0.82246\ldots i$ in the unit disc, the reduction algorithm for $\Gamma$ yields the point $\gamma w = -0.12450\ldots+0.02857\ldots i$ in $D_\Gamma$. The reduction algorithm for $\Delta$ yields the point $\delta w = 0.08700\ldots+0.09353\ldots i$ in $D_\Delta$.} \label{fig:reduction3} 
\end{figure}

\end{exm}

\subsection*{Algorithm to compute a group presentation}

We conclude this section by showing how to determine a presentation for the group $\Gamma$ from the coset graph.  Although this will not figure in our method to compute \Belyi\ maps, it can be computed easily in our setup and may be useful for other applications.  Let $D_\Gamma$ be a fundamental domain for $\Gamma$ as in \ref{eqn:fundGamma} and let $S$ be the associated side pairing. Having computed a side pairing for $\Gamma$, all of the relevant algorithms for obtaining a presentation for $\Gamma$ appear in work of Voight \cite{Voightfd}.

\begin{rmk}
The assumption that the fundamental domains of the Fuchsian groups under examination are hyperbolically convex is verified \cite{Voightfd} in order to apply the Poincar\'{e} polygon theorem. However, this condition is extraneous; the needed general version of the Poincar\'{e} polygon theorem was proven by Maskit \cite{Maskit}. In particular, even when the fundamental domain is disconnected, the Poincar\'{e} polygon theorem can still be applied: see Epstein and Petronio \cite[Theorem 8.1]{EP} (which also includes the colorful history of this theorem).
\end{rmk}

The \defi{vertices} of $D_\Gamma$ are the intersections of the sides of $D_\Gamma$ in $S$ (i.e. the vertices of $D_\Gamma$ as a hyperbolic polygon).  A \defi{pairing cycle} for $D_\Gamma$ is a sequence $v_1,\dots,v_n$ of vertices of $D_\Gamma$ which is the (ordered) intersection of the $\Gamma$-orbit of of $v_1$ with $D_\Gamma$.  To each cycle associate a word $g = g_n g_{n-1} \cdots g_2 g_1$ where $g_i(v_i) = v_{i+1}$.  A cycle is \defi{minimal} if $v_i \neq v_j$ for all $i \neq j$. Any side pairing element $\gamma \in \Gamma$ of $S$ will appear at most once in any word associated to a minimal cycle.  A set of minimal cycles is \defi{complete} if every side pairing element in $S$ occurs in (a necessarily unique) one of the cycles in the set.  Given $S$ there is an algorithm to compute a complete set of minimal cycles for $D_\Gamma$ (see \cite{Voightfd} Algorithm 5.2).  Given such a complete set of minimal cycles of $D_\Gamma$ there is an algorithm to compute a minimal set of generators and relations for $\Gamma$ \cite[Algorithm 5.7]{Voightfd}.

We illustrate these methods with an example.

\begin{exm}
Let $\Gamma \leq \Delta(3,5,3)$ correspond to the permutation triple 
\[ \sigma_0 = (1\;3\;4)(2\;6\;5),\ \sigma_1 = (1\;4\;2\;5\;6),\ \sigma_\infty = (3\;6\;4). \] 
We first obtain a fundamental domain for $\Gamma$ using the petalling approach.

\begin{figure}[h] 
\includegraphics{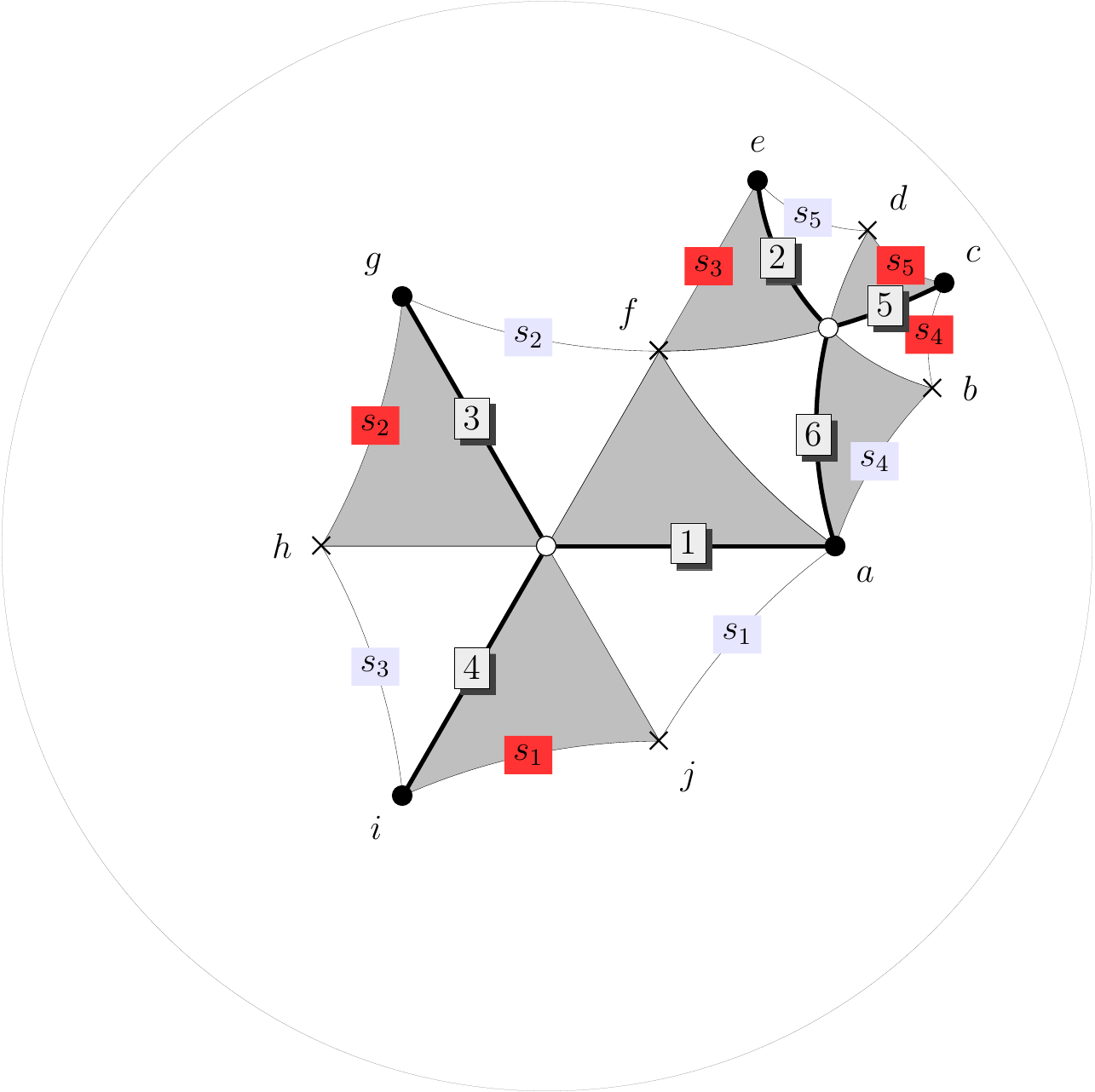}

\caption{A labelled fundamental domain for $\Gamma$.} \label{fig:presentation} 
\end{figure}

Recall that the side pairings follow the convention that we move from red labels to blue labels. We now compute the pairing cycles for the preimages of $1$ (black dots). These pairing cycles are
\begin{align*}
&a\stackrel{s_1^{-1}}{\longrightarrow} i \stackrel{s_3^{-1}}{\longrightarrow}e\stackrel{s_5^{-1}}{\longrightarrow}c\stackrel{s_4}{\longrightarrow}a\\
&g\stackrel{s_2}{\longrightarrow}g
\end{align*}
which correspond to the words $\gamma_4\gamma_5^{-1}\gamma_3^{-1}\gamma_1^{-1}$ and $\gamma_2$ respectively. Some power of each word will be a relation in our presentation where this power is the order of the stabilizer of the associated point. This can be observed by dividing the order of $\sigma_1$ by the number of triangles incident to this point. Since the order of $\sigma_1$ is $5$ and there are $5$ triangles incident to $a$, we obtain the relation $\gamma_4\gamma_5^{-1}\gamma_3^{-1}\gamma_1^{-1}$. Similarly, since only one triangle is incident to $g$, we obtain the relation $\gamma_2^{5}$.

Now we compute the pairing cycles for the preimages of $\infty$ (crosses). These pairing cycles are 
\begin{align*} 
&b\stackrel{s_4}{\longrightarrow} b\\
&d\stackrel{s_5}{\longrightarrow}d\\
&f\stackrel{s_3}{\longrightarrow} h\stackrel{s_2}{\longrightarrow} f\\
&j\stackrel{s_1}{\longrightarrow} j
\end{align*}
which correspond to the words $\gamma_4$, $\gamma_5$, $\gamma_2\gamma_3$, $\gamma_1$ respectively. Again, since the order of $\sigma_{\infty}$ is $3$, by considering the stabilizers we find the relations $\gamma_4^3$, $\gamma_5^3$, $\gamma_2\gamma_3$, $\gamma_1^3$. Thus
$$
\langle \gamma_1,\gamma_2,\gamma_3,\gamma_4,\gamma_5\mid \gamma_2^5=\gamma_4\gamma_5^{-1}\gamma_3^{-1}\gamma_1^{-1}=\gamma_4^3=\gamma_5^3=\gamma_2\gamma_3=\gamma_1^3=1\rangle
$$
is a presentation for $\Gamma$.  This presentation simplifies to
$$
\langle \gamma_1,\gamma_2,\gamma_3\mid \gamma_1^3=\gamma_3^3=\gamma_2^5=(\gamma_3\gamma_1^{-1}\gamma_2)^3=1\rangle
$$
the recognizable $2$-orbifold group of signature $(0;3,3,3,5;0)$ associated to a surface with genus zero and $3$ (resp.~$1$) orbifold points of orders $3$ (resp.~$5$), respectively.  
\end{exm}

\section{Power series expansions of differentials} \label{sec:powser}

In this section, we discuss power series expansions of differentials on Riemann $2$-orbifolds in the language of modular forms on cocompact Fuchsian groups and we describe a method for computing these expansions following Voight--Willis \cite{VoightWillis}; the key ingredients are the tools developed in the previous section involving reduction to a fundamental domain.

\subsection*{Computing power series expansions}

Let $\Gamma < \PSL_2(\R)$ be a cocompact Fuchsian group.  A \defi{modular form} for $\Gamma$ of \defi{weight} $k \in 2 \Z$ is a holomorphic function $f : \calH \to \C$ such that 
\begin{equation} \label{eqn:fgamma}
f(\gamma z) = j(\gamma,z)^k f(z) \text{ for all $z \in \calH$ and all $\gamma \in \Gamma$}
\end{equation} 
where $j(\gamma ,z) = c z+ d$ for $\gamma = \pm \begin{pmatrix} a & b \\ c & d \end{pmatrix} \in \PSL_2(\R)$.  

Equation (\ref{eqn:fgamma}) is equivalent to
\begin{equation} \label{eqn:dgammaz}
f(\gamma z)\, d(\gamma z)^{\otimes k/2} = f(z)\, dz^{\otimes k/2} 
\end{equation}
so that a modular form $f$ of weight $k$ can be thought of as a holomorphic differential $k/2$-form on the quotient $X=\Gamma \backslash \calH$.  Since $X=\Gamma \backslash \calH$ has no cusps, a modular form for $\Gamma$ will also sometimes be called a \defi{cusp form} (extending the classical terminology).  Let $S_k(\Gamma)$ be the $\C$-vector space of cusp forms of weight $k$ for $\Gamma$.  

Let $f \in S_k(\Gamma)$ be a modular form of weight $k$ for $\Gamma$.  Since $\Gamma$ is cocompact, we do not have $q$-expansions available; however, $f$ is holomorphic and therefore has a Taylor series expansion.  Let $p \in \calH$.  By conformally mapping $\calH$ to $\calD$ with the map
\begin{align*}
w_p : \calH & \to \calD \\
z &\mapsto w_p=w = \frac{z-p}{z-\overline{p}} 
\end{align*}
we express $f$ as a Taylor series in $w$ that will be valid for all $z \in \calH$.  So we define a \defi{power series expansion} for $f$ centered at $p$ to be an expression of the form
\begin{equation} \label{eqn:fbnpow}
f(z) = (1-w)^k \sum_{n=0}^\infty b_n w^n
\end{equation}
where $b_n \in \C$ and $w = w_p(z)$.  The factor $(1-w)^k$ is the automorphy factor corresponding to the linear fractional transformation $w$ and is included for arithmetic reasons \cite{VoightWillis}.  

The basic idea of the algorithm is as follows: we use the relation of modularity (\ref{eqn:fgamma}) and reduction to the fundamental domain to solve for the unknown coefficients $b_n$.  Specifically, let $D_\Gamma$ be a fundamental domain for $\Gamma$ acting on $\calD$ as in (\ref{eqn:fundGamma}), and let $\rho > 0$ be the radius of a circle centered at the origin containing $D_\Gamma$, possible since $\Gamma$ is cocompact.  To compute the power series expansion for $f$ centered at $p \in \calH$ to some precision $\epsilon > 0$, we consider a truncation 
\begin{equation} \label{eqn:fNz}
f(z) \approx f_N(z) = (1-w)^k \sum_{n=0}^N b_n w^n 
\end{equation}
valid for all $|w| \leq \rho$ and some $N \in \Z_{\geq 0}$.  For a point $w = w_p(z) \in \calD$ on the circle of radius $\rho$ with $w \notin D_\Gamma$, we use Algorithm \ref{alg:fdreddelta} to find $\gamma \in \Gamma$ such that $w' = \gamma w \in D_\Gamma$.  Letting $z' = w_p^{-1}(w')$, by the modularity of $f$ we have
\[ (1-w')^k \sum_{n=0}^N b_n (w')^n \approx f_N(z') = j(\gamma,z)^kf(z) \approx j(\gamma,z)^k (1-w)^k \sum_{n=0}^N b_n w^n \]
valid to the precision $\epsilon >0$.  For each such point $w$, we obtain such a relation; therefore by taking enough such points we can obtain enough nontrivial linear relations on the unknown coefficients $b_n$ to determine them. 

An alternative method for obtaining linear relations amongst the coefficients $b_n$ which we found to have greater numerical stability uses Cauchy's integral formula.  We identify $\calH$ with $\calD$ via the map $w$, we abuse notation slightly and write simply $f(w) = f(w(z))$.  From the expansion (\ref{eqn:fbnpow}), we apply Cauchy's integral formula to $f(w)/(1-w)^k$: for $n \ge 0$, we have
$$
b_n = \frac{1}{2\pi i} \int_C \frac{f(w)}{w^{n+1}(1-w)^k}\,dw
$$
where $C$ is a simple contour around $0$.  We take $C(\theta) = \rho e^{i \theta}$ for $0 \le \theta \le 2 \pi$ and thus obtain
$$
b_n = \frac{1}{2\pi} \int_0^{2\pi} \frac{f(\rho e^{i\theta})}{(\rho e^{i\theta})^{n}(1-\rho e^{i\theta})^{k}}\,d\theta .
$$
Breaking up $[0, 2 \pi]$ into $Q \in \Z_{\ge0}$ intervals and letting $w_m = \rho \exp(2 \pi m i/Q)$ for $m=1,\dots,Q$, we obtain the approximation (by Riemann summation)
$$
b_n \approx \frac{1}{Q} \sum_{m=1}^{Q} \frac{f(z_m)}{w_m^n(1-w_m)^k}
$$
valid to precision $\epsilon$ for sufficiently large $Q$.  For each $m$, let $\gamma_m \in \Gamma$ be such that $w_m'=\gamma w_m \in D$, and let $z_m=\phi^{-1}(w_m)$ and $z_m'=\phi^{-1}(w_m')$.  Then by the modularity of $f$, we have
$$
b_n \approx \frac{1}{Q} \sum_{m=1}^{Q} \frac{f(z_m)}{w_m^n(1-w_m)^k} \approx \frac{1}{Q} \sum_{m=1}^{Q} \frac{f_N(z_m') j(\gamma_m, z_m)^{-k}}{w_m^n(1-w_m)^k}.
$$
Expanding $f_N(z)$ as in (\ref{eqn:fNz}) and substituting, we obtain an approximate linear equation involving the coefficients $b_n$.  Specifically, we have
\begin{equation} \label{eqn:bigdealrelation}
b_n \approx \sum_{r=0}^N a_{nr} b_r
\end{equation}
where
\begin{equation} \label{eqn:jwmfact}
a_{nr} = \frac{1}{Q} \sum_{m=1}^Q j(\gamma_m,z)^{-k}\frac{(w_m')^r (1-w_m')^k}{w_m^n(1-w_m)^k} .
\end{equation}
Let $A$ be the matrix with entries $a_{nr}$, rows indexed by $n = 0, \dots, N$, and columns indexed by $r = 0, \dots, N$.  If we write $b$ as the column vector with entries $b_n$ for $n = 0, \dots, N$, we have that
\begin{equation} \label{eqn:yayA}
Ab \approx b
\end{equation}
and therefore $b$ is approximately in the eigenspace for the eigenvalue $1$ of the matrix $A$.  

In this way, computing the approximate eigenspace for $1$ of the matrix $K$ yields a basis for the space $S_k(\Gamma)$.  In the next section, we discuss how to take this algorithm to furnish equations for \Belyi\ maps.  

\subsection*{Algorithmic improvements}

We now discuss a few improvements on this algorithm which extend the range of practical computation.  

First, from (\ref{eqn:jwmfact}), conveniently the matrix $K$ factors as 
\begin{equation} \label{matrixmult}
QK = JW'
\end{equation}
where $J$ is the matrix with entries 
$$
J_{nm}=\frac{j(\gamma_m,z_m)^{-k}}{w_m^n(1-w_m)^k} 
$$
with $0 \leq n \leq N$ and $1 \leq m \leq Q$ and $W'$ is the Vandermonde-like matrix with entries
$$
W'_{mr} = (w_m')^r (1-w_m')^k
$$
with $1 \leq m \leq Q$ and $0 \leq r \leq N$.  (The matrix $W'$ further factors as the product of an honest Vandermonde matrix and a diagonal matrix.)  Since $J$ and $W'$ are both fast to compute, the computation of $K$ requires just a single matrix multiplication.  Taking $Q=O(N)$, one can compute the matrix $K$ using $O(N^3)$ ring operations in $\C$ using storage space $O(N^2)$; however, when $N$ is large, one can do better by finding the needed eigenspace for $K$ without ever computing the matrix $K$ directly.

Second, we apply a federalist approach: we consider power series expansions at several points $p_1,\dots,p_s$, and express $f$ as piecewise-defined power series expansions in the associated Voronoi regions (evaluating in the neighborhood of the closest point).  As the number of points grows, the expansion degree $N$ drops as points can be taken in smaller neighborhoods (smaller radius); in general, the geometry of the fundamental domain can lead to several optimal configurations of points.  In our context, where $\Gamma \leq \Delta$ is a subgroup of finite index with cosets $\Gamma \alpha_i$ for $i=1,\dots,d$ and fundamental domain $D_\Gamma = \bigcup_{i=1}^d \alpha_i D_\Delta$ as in (\ref{eqn:fundGamma}), it is simpler to take expansions at the set of points $\{\alpha_i z_a : i=1,\dots,d\}$ as the output, since the reduction algorithm (Algorithm \ref{alg:fdredgamma}) exactly produces a point in a coset translate $\alpha_i D_\Delta$.  In this way, the expansion degree $N$ depends only on that needed for convergence inside $D_\Delta$.  In this approach, we repeat the application of Cauchy's theorem around each center $p_i$, so the resulting matrix is nearly block diagonal, but with some points evaluated in other neighborhoods: examples of what this approach looks like in practice can be found in section \ref{sec:examples}.

Finally, we use the method of Krylov subspaces and Arnoldi iteration, a standard technique in numerical linear algebra \cite{GolubvanLoan}, to compute the eigenspace for eigenvalue $1$.  Given an initial random vector $v$ and $n \in \Z_{\geq 0}$, if the sequence $v, Av, A^2 v, \dots$ converges, it converges to an eigenvector for $A$ with eigenvalue of largest absolute value; this is called the \defi{power method}.  To get finer information, we define the \defi{Krylov subspace} $K_n$ for a matrix $A$ to be the span of the vectors $v, Av, A^2 v, \dots, A^n v$.  The approximate eigenspaces for $A$, when intersected with $K_n$, converge to the full approximate eigenspaces, so we can expect that the intersections are good approximations to the full eigenspaces for the large eigenvalues of $A$ when $n$ is not too large; then another technique (Gram-Schmidt orthogonalization, LU decomposition, SVD decomposition, etc.)\ can be applied to the subspace more efficiently, as this subspace is of much smaller dimension.  However, this method is unstable as described, and \defi{Arnoldi iteration} uses a stable version of the Gram-Schmidt process to find a sequence of orthonormal vectors spanning the Krylov subspace $K_n$.  Many implementations of Arnoldi iteration typically restart after some number of iterations (such as the Implicitly Restarted Arnoldi method, IRAM).  More experiments are needed to fully optimize these techniques in our setting, to select the right variant of these methods and to choose optimal parameters; we leave these experiments for future work.

\subsection*{Hypergeometric series}

To conclude this section, we discuss the connection between hypergeometric series and an explicit power series expansion for a uniformizer $\phi$ for $\Delta \backslash \calH$ where $\Delta$ is a triangle group; this provides an analytic expansion for the \Belyi\ map for a subgroup $\Gamma \leq \Delta$ of finite index.

Let $\Delta=\Delta(a,b,c)$ be a triangle group.  Then $X(\Delta) = \Delta \backslash \calH$ naturally possesses the structure of a Riemann $2$-orbifold with three elliptic points of orders $a,b,c$.  Further, it can be given the structure of a Riemann surface of genus zero, and consequently there is an isomorphism $\phi:X(\Delta) \xrightarrow{\sim} \PP^1(\C)$ of Riemann surfaces that is unique if we insist that the elliptic points map to the points $0,1,\infty$.  The function $\phi$ is not holomorphic (it has a single simple pole at the elliptic point of order $c$), but it is still well-defined as a function on $X(\Delta)$, so that $\phi(\delta z) = \phi(z)$ for all $\delta \in \Delta$, and so we say that $\phi$ is a \defi{meromorphic modular form} for $\Gamma$ of weight $0$.

We consider a power series expansion for $\phi$ in a neighborhood of $p=z_a$ of the form 
\[ \phi(w) = \phi(w(z)) = \sum_{n=0}^{\infty} b_n w^{n}. \]
But since $\delta_a$ fixes $z_a \in \calH$ and so acts by rotation by $\zeta_a=\exp(2\pi i/a)$ in $\calD$ centered at $p$, we have
\[ \phi(\delta_a w) = \phi(\zeta_a w)=\phi(w) \]
so in fact
\[ \phi(w) = \sum_{n=0}^{\infty} b_{an} w^{an}. \]
In particular, $\phi:\calD \to \PP^1(\C)$ is an $a$-to-1 map at $w_a=0 \in \calD$.  We now give an explicit expression for $\phi$ and its coefficients $b_n$.  

We use the classical fact that the functional inverse of $\phi$ can be expressed by the quotient of two holomorphic functions that form a basis for the vector space of solutions to a second order linear differential equation {\cite[Theorem 15, \S 44]{Ford}}.  For the triangle group $\Delta(a,b,c)$, this differential equation has $3$ regular singular points \cite[\S 109]{Ford} corresponding to the elliptic points.  After a bit of simplification \cite[\S 2.3]{AndrewsAskey} we recognize this differential equation as the \defi{hypergeometric ${}_2F_1$ differential equation}
\begin{equation} \label{eqn:hypergeom}
z(1-z)\frac{d^2y}{dz^2} + \left(\gamma-(A+B+1)y\right)\frac{dy}{dz} - AB y =0 
\end{equation}
where 
\begin{equation} \label{eqn:alphabetagamma}
A=\frac{1}{2}\left(1+\frac{1}{a}-\frac{1}{b}-\frac{1}{c}\right),\quad
B=\frac{1}{2}\left(1+\frac{1}{a}-\frac{1}{b}+\frac{1}{c}\right),\quad 
C=1+\frac{1}{a}.
\end{equation}
(For other Fuchsian groups with at least 4 elliptic points or genus $>0$, one must instead look at the relevant Schwarzian differential equation; see Elkies \cite[p.~8]{ElkiesSCC} for a discussion, and Ihara \cite{Ihara} for the general case.)

The solutions to the equation (\ref{eqn:hypergeom}) are well-studied.  For a moment, let $A,B,C \in \R$ be arbitrary real numbers with $C \not\in \Z_{\leq 0}$.  Define the \defi{hypergeometric series} 
\begin{equation} \label{eqn:Fabc}
F(A,B,C;t)= \sum_{n=0}^{\infty} \frac{(A)_n(B)_n}{(C)_n}\frac{t^n}{n!} \in \R[[t]] 
\end{equation}
where $(x)_n = x(x+1) \cdots (x+n-1)$ is the \defi{Pochhammer symbol}.  This series is also known as the \defi{${}_2F_1$ series} or the \defi{Gaussian hypergeometric function} (see Slater \cite[\S 1.1]{Slater} for a historical introduction).  The hypergeometric series is convergent for all $t \in \C$ with $|t|<1$ by the ratio test and is a solution to the hypergeometric differential equation by direct substitution.  A basis of solutions in a neighborhood of $t=0$, convergent for $|t|<1$, is given by \cite[1.3.6]{Slater}
\begin{align*} 
F_1(t) &=F(A,B,C;t) \\
F_2(t) &= t^{1-C}F(1+A-C,1+B-C,2-C;t).
\end{align*}

To obtain locally an inverse to $\phi$ as above, since $1-C=-1/a$, we have 
\[ F_1(t)/F_2(t) \in t^{1/a} \Q[[t]] \]
which uniquely defines this as the correct ratio up to a scalar $\kappa \in \C^\times$.  The value of $\kappa$ can be recovered by plugging in $t=1$: when $C-A-B>0$, we have \cite[1.7.6]{Slater}
\begin{equation} \label{eqn:FABC1}
F(A,B,C;1)=\frac{\Gamma(C)\Gamma(C-A-B)}{\Gamma(C-A)\Gamma(C-B)}
\end{equation}
where $\Gamma$ is the complex $\Gamma$-function.  Since
\[ 2-C-(1+A-C)=C-A-B =\left(1+\frac{1}{a}\right)-\left(1+\frac{1}{a}+\frac{1}{b}\right)=\frac{1}{c}>0 \]
we can use (\ref{eqn:FABC1}) to recover $\kappa$ by solving 
\[\kappa \frac{F_1(1)}{F_2(1)}=w_b=\frac{z_b-z_a}{z_b-\overline{z_a}} \]
where $z_a=i$ and $z_b=\newt i$ are as in Proposition \ref{prop:triangleembed}.  We obtain
\begin{equation} \label{eqn:kdef}
 \kappa=\left(\frac{\newt-1}{\newt+1}\right) \frac{\Gamma(2-C)\Gamma(C-A)\Gamma(C-B)}{\Gamma(1-A)\Gamma(1-B)\Gamma(C)} \in \R.
\end{equation}
Therefore, the function 
\begin{equation} \label{eqn:psit}
\psi(t) = \kappa\frac{F_1(t)}{F_2(t)} 
\end{equation}
maps the neighborhood of $t=0$ with $|t|<1$ to a neighborhood of $w=0 \in \calD$ (intersected with $D_\Delta$) and $\phi$ is an inverse for the map $\psi$.  

\begin{rmk}
If we make a cut in the $t$-plane along the real axis from $t=1$ to $t=\infty$, that is, we insist that $-\pi < \arg(-t) \leq \pi$ for $|t| \geq 1$, then the function $F(A,B,C;t)$ defined by the hypergeometric series can be analytically continued to a holomorphic function in the cut plane.  Similar formulas then hold for explicitly given Puiseux series expansions around the two other elliptic points \cite[\S 5.2]{VoightThesis}.
\end{rmk}

\begin{figure}[h] 
\includegraphics{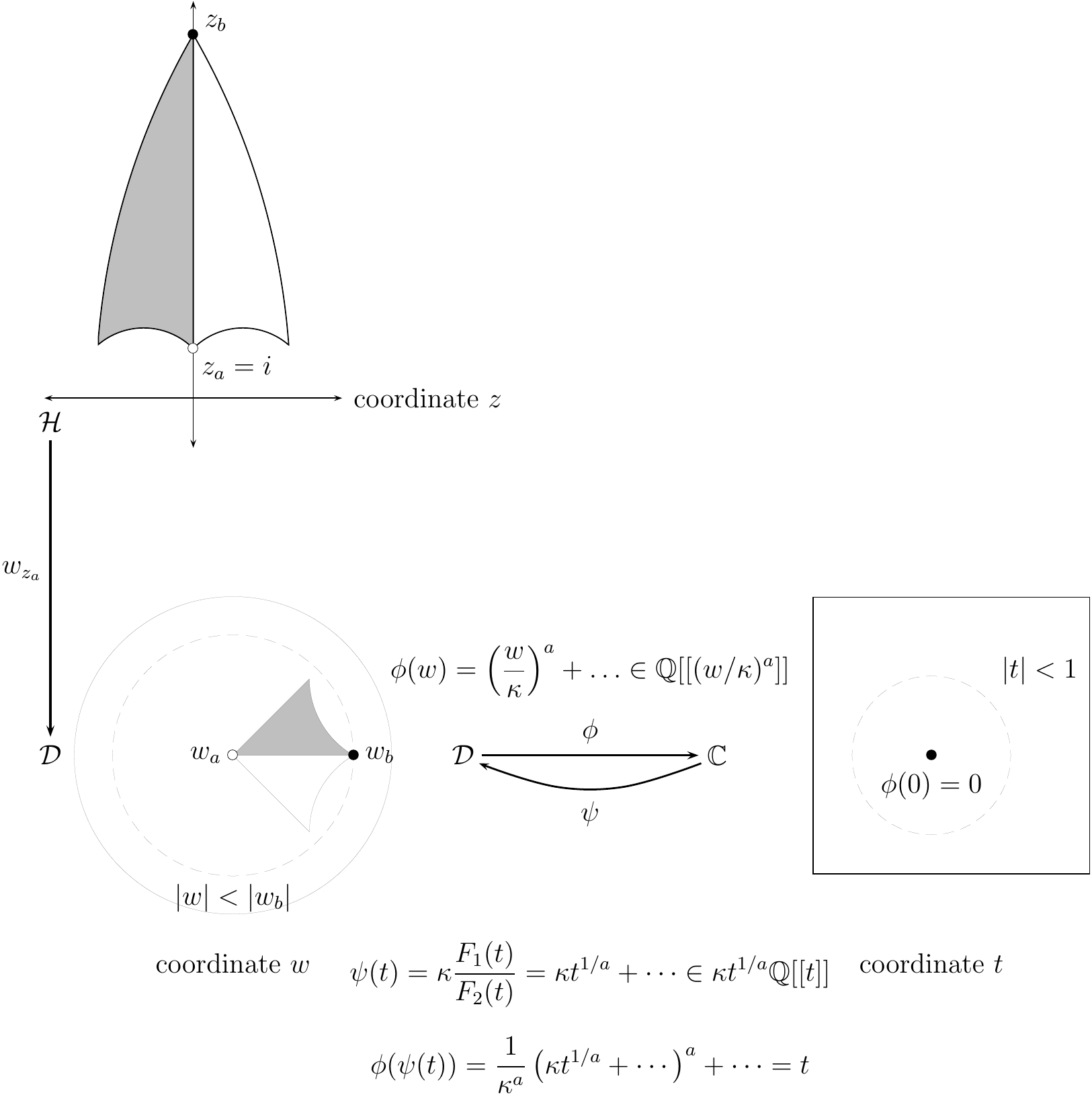}

\caption{Hyperbolic uniformization of the triangle quotient via hypergeometric series.} \label{fig:puiseux}
\end{figure}

Computing the reversion of the Puiseux series $\phi(t)=w$ then gives a power series expansion for $\psi(w)=t$.  We have
\begin{equation} \label{eqn:expandpuiseux}
\phi(w) = \left(\frac{w}{\kappa}\right)^a + \frac{a(C(1+A-C)(1+B-C)-AB(2-C))}{C(2-C)}\left(\frac{w}{\kappa}\right)^{2a} + \dots \in \Q\left[\!\left[\left(\frac{w}{\kappa}\right)^a\right]\!\right].
\end{equation}
Note in particular that the expansion (\ref{eqn:expandpuiseux}) has rational coefficients when considered as a series in $w/\kappa$.

\section{Computing equations, and examples} \label{sec:examples}

In this section, we compute equations for \Belyi\ maps combining the algorithms of the past two sections.  Along the way, we present several explicit examples in full.

Throughout, we retain the notation introduced in section 3: in particular, let $\Gamma \leq \Delta=\Delta(a,b,c)$ be a subgroup of index $d$ corresponding to a transitive permutation triple $\sigma \in S_d^3$.  We compute defining equations for the \Belyi\ map
\[ \phi : X(\Gamma)=\Gamma \backslash \calH \to X(\Delta) = \Delta \backslash \calH \]
defined over a number field.  Our description and method depend on the genus, with low genus cases handled separately.

\subsection*{Genus zero}

First, suppose that $X(\Gamma)$ has genus zero.  Then we have an analytic isomorphism $x:X(\Gamma) \cong \PP_\C^1$ of Riemann surfaces and the \Belyi\ map $\phi:X(\Gamma) \to X(\Delta)$ corresponds to a rational function $\phi(x) \in \C(x)$ of degree $d$.  The isomorphism $x$ is not uniquely specified.  One way to fix this is to specify the images of three points on $X(\Gamma)$, as was done for the function $\phi$; but this may extend the field of definition.  Instead, let $e$ be the order of the stabilizer of $w=0$, so that $x \in \C[[w^e]]$; we require that $x$ has Laurent series expansion
\begin{equation} \label{eqn:xw}
x(w) = (\Theta w)^e + O(w^{3e})
\end{equation}
for some $\Theta \in \C^\times$.  For a given choice of $\Theta$, this normalization specifies $x$ uniquely.  Experimentally, we observe that a choice
\begin{equation} \label{eqn:Theta0}
\Theta = \sqrt[a]{\alpha} \left(\frac{1}{\kappa}\right)
\end{equation}
for some $\alpha \in \overline{\Q}$ yields a series
\[ x(w) = \sum_{\substack{n=e \\ e \mid n}}^{\infty} \frac{c_n}{n!} (\Theta w)^n \]
(weight $k=0$ so no factor $(1-w)^k$ appears) with Taylor coefficients $c_n$ that belong to a number field $K$ (in fact, with small denominators). 

\begin{rmk} \label{rmk:Shimura}
The algebraicity of the coefficients is similar in spirit to results of Shimura \cite{Shimura1,Shimura2,Shimura3}, who proved algebraicity of $c_n$ when the group $\Gamma$ is \defi{arithmetic}, commensurable with the units in the maximal order of a quaternion algebra defined over a totally real field and split at a unique real place.  There are exactly $85$ triples $(a,b,c)$ with $a \leq b \leq c$ that are arithmetic, by a theorem of Takeuchi \cite{Takeuchi,Takeuchi2}; so in the infinitely many remaining cases, one does not know that the coefficients are algebraic, and we merely report the observation that the coefficients seem to be algebraic (and, in fact, nearly integral after suitable normalization).
\end{rmk}

We compute an expansion for $x$ as follows.  First, we find $k \in 2\Z_{\geq 0}$ such that $\dim_\C S_k(\Gamma) \geq 2$.  The dimension of $S_k(\Gamma)$ is given by Riemann-Roch \cite{Shimura}: we have $\dim_\C S_2(\Gamma)=g$ and
\begin{equation} \label{eqn:dimformula}
\dim_\C S_k(\Gamma) = (k-1)(g-1) + \displaystyle{\sum_{i=1}^r \left\lfloor \frac{k(e_i-1)}{2e_i} \right\rfloor}, \quad \text{ if $k \geq 4$}.
\end{equation}
Next, we compute an \defi{echelonized basis} \cite[\S 4.4]{Klug} of functions with respect to the coefficients $b_n$.  Let 
\[ s \equiv \frac{k}{2}(e-1) \pmod e \]
with $0 \leq s < e$.  Then the expansions in $S_k(\Gamma)$ belong to $w^s \C[[w^e]]$.  We find functions
\begin{equation} \label{eqn:gandh}
g(w)=w^{m} + O(w^{m+2e})\quad \text{and} \quad h(w)=w^{m+e} + O(w^{m+2e})
\end{equation}
with $m \equiv s \pmod{e}$ maximal with this property and then take 
\[ x(w)=\frac{\Theta^e h(w)}{g(w) + ch(w)} \]
where $c$ is the coefficient of $w^{m+2e}$ in $h$ to obtain an expansion as in \eqref{eqn:xw}.  The fact that the function $x(w)$ has degree $1$ and thereby gives an isomorphism of Riemann surfaces $x:X=\Gamma \backslash \calH \to \PP_\C^1$ follows from Riemann-Roch: if $K$ is a canonical divisor for orbifold $X$ and $m=\dim H^0(X,dK)$ with $d=k/2$, then $g,h$ span $H^0(X,dK-(m-1)P)$ where $P$ is the point $w=0$.  Finally, using the expansion for $\phi$ in (\ref{eqn:expandpuiseux}), we use linear algebra to find the rational function $\phi(x)$ of degree $d$.

\begin{rmk}
We could instead compute the ring of modular forms for $\Gamma$, making no choices at all: the structure of this ring is described by Voight--Zureick-Brown \cite{JVDZB} and contains complete information about differentials on the orbifold $X(\Gamma)$.  This remark holds in all genera below as well.
\end{rmk}

These calculations are all performed over the complex numbers.  The \Belyi\ map $\phi(x)$ will be defined over a number field of degree at most equal to the size of the corresponding refined passport, and so one can use the LLL algorithm \cite{LLL} to find putative polynomials that these coefficients satisfy.  Alternatively, one can compute the maps $\phi(x)$ for all Galois conjugates (possibly belonging to different refined passports), and then, using continued fractions, recognize the symmetric functions in these conjugates to obtain the minimal polynomial for each coefficient.

\begin{exm}
The smallest degree $d$ for which there exists a hyperbolic, transitive permutation triple $\sigma$ with genus $g=0$ is $d=5$.  There are nine refined passports and each is rigid, having a unique triple up to simultaneous conjugation.  Two triples generate $G=A_5\cong \PSL_2(\F_5)$ and seven generate $G=S_5 \cong \PGL_2(\F_5)$.  We consider the first two triples in detail; two of the resulting maps from the latter collection are merely reported below.

First, consider the triples with $G=A_5$: they both have $(a,b,c)=(5,3,3)$, up to reordering.  The conjugacy class of an element of order $3$ is unique, but there are two conjugacy classes in $A_5$ of elements of order $5$, and the two triples correspond to these two choices.  The two conjugacy classes of order $5$ are interchanged by the outer automorphism of $A_5$ obtained from conjugating by an element in $\PGL_2(\F_5) \setminus \PSL_2(\F_5)$.  We will see (in several ways) below that this automorphism corresponds to complex conjugation on the associated \Belyi\ maps, and that the \Belyi\ map is defined over $\Q$ (and hence over $\R$).

We begin with the first triple
\begin{equation} \label{eqn:sigmaexmgenus0}
\sigma_0=(1\ 5\ 4\ 3\ 2), \quad \sigma_1 = (1\ 2\ 3), \quad \sigma_\infty=(3\ 4\ 5).
\end{equation}
Let $\Gamma \leq \Delta(5,3,3)$ be the group of index $5$ associated to this triple as in (\ref{eqn:gammaperm}).  The ramification diagram then looks as in Figure \ref{fig:ramification}.  The points are labelled with the order of their stabilizer group (no label indicates trivial stabilizer).  

\begin{figure}[h] 
\includegraphics{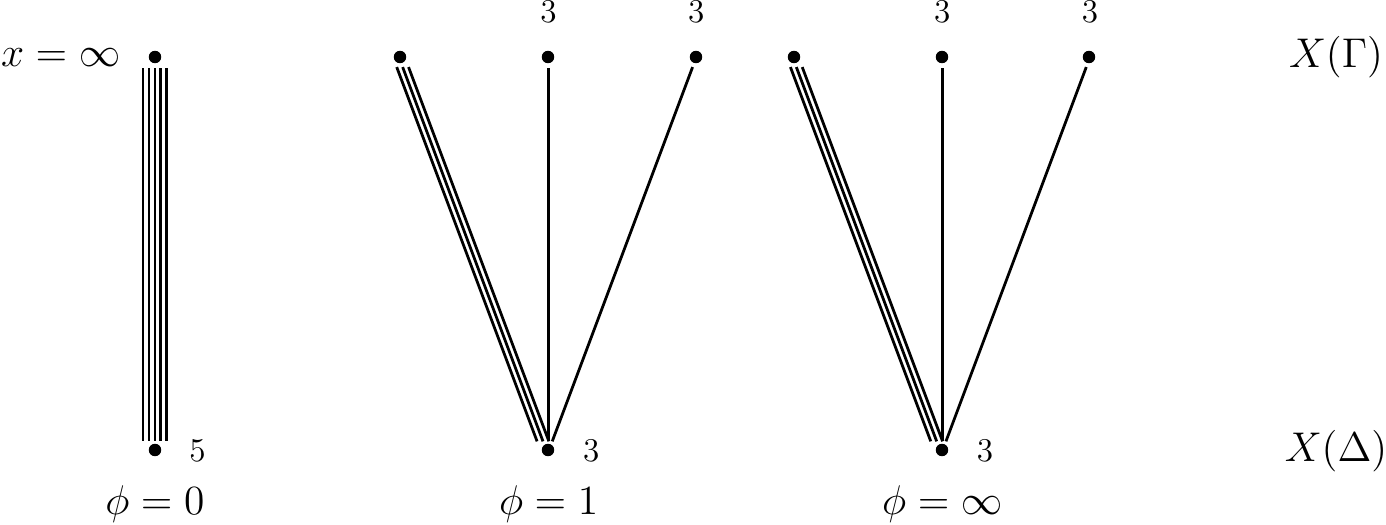}

\caption{The ramification diagram for the triple (\ref{eqn:sigmaexmgenus0}).} \label{fig:ramification}
\end{figure}

From the Riemann-Hurwitz formula (\ref{eqn:RH}) we find that $g(X(\Gamma))=0$.  We have $e=1$ so $s=0$ and  our power series belong to ring $\C[[w]]$ in all weights $k$.

\begin{figure}[h]

\includegraphics{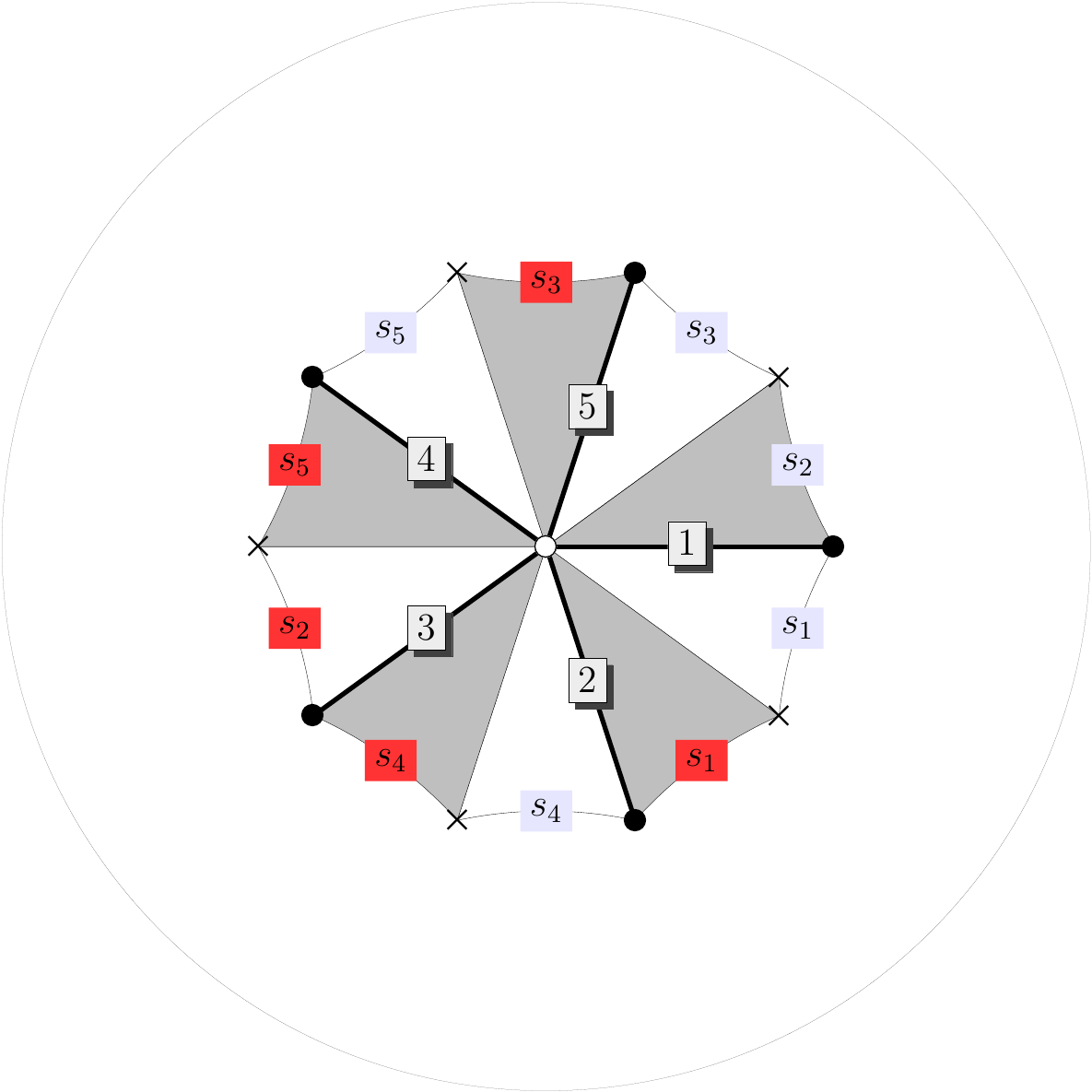}

\begin{center}
\begin{tabular}{ll}
\toprule
Label & Coset Representative\\
\midrule
$1$ & 1 \\
$2$ & $\delta_a^{-1}$ \\
$3$ & $\delta_a^{-2}$ \\
$4$ & $\delta_a^{2}$ \\
$5$ & $\delta_a^{}$ \\
\bottomrule
\end{tabular}
\;\;\;\;\;\;\;\;\;\;
\begin{tabular}{ll}
\toprule
Label & Side Pairing Element\\
\midrule
$s_{1}$ & $\delta_b^{}\delta_a^{}$ \\
$s_{2}$ & $\delta_b^{-1}\delta_a^{2}$ \\
$s_{3}$ & $\delta_a^{}\delta_b^{}\delta_a^{-1}$ \\
$s_{4}$ & $\delta_a^{-1}\delta_b^{}\delta_a^{2}$ \\
$s_{5}$ & $\delta_a^{2}\delta_b^{}\delta_a^{-2}$ \\
\bottomrule
\end{tabular}
\end{center}
\bigskip
\caption{A fundamental domain and conformally drawn dessin for $\Gamma$ corresponding to the permutation triple \eqref{eqn:sigmaexmgenus0} using the methods in  section \ref{sec:cosets}.}  \label{fig:dessin-5,3,3}
\end{figure}

From the dimension formula (\ref{eqn:dimformula}), we compute that
\[ \dim_\C S_k(\Gamma) = 0,1,3 \quad \text{ for $k=2,4,6$} \]
so we look at the space of modular forms of weight $6$ for $\Gamma$.  

We work in precision $\varepsilon=10^{-30}$; the fundamental domain is contained in a circle of radius $\rho=0.528935$, and so (assuming $b_n=O(1)$) we take $N= \lceil \log \varepsilon/\log \rho \rceil = 109$ and $Q=2N$.  We compute the matrix $A$ (\ref{eqn:yayA}) and the singular value decomposition of the matrix $A-1$ and find it has three small singular values, the largest of which is $< 8 \cdot 10^{-29}$, yielding an approximate eigenspace for $A$ with eigenvalue $1$ of dimension $3$.  Computing as in (\ref{eqn:gandh}) and taking
\begin{equation} \label{eqn:exmgen0Theta}
\Theta = 0.3917053\ldots + 1.205545\ldots i= \sqrt[5]{\frac{81}{2}}\exp(2\pi i/5)\left(\frac{1}{\kappa}\right) 
\end{equation}
where $\kappa$ is as in (\ref{eqn:kdef}), we find the (approximate) basis
\begin{align*}
f(w)&=(1-w)^6\left(1 - \frac{81}{4!}(\Theta w)^4 - \frac{6075}{6!} (\Theta w)^6 + \frac{382725}{2\cdot 8!} (\Theta w)^8 \right.  \\
& \qquad\qquad\qquad\qquad \left. - \frac{24111675}{10!} (\Theta w)^{10} + O(w^{14})\right) \\
g(w)&=(1-w)^6\left((\Theta w) + \frac{9}{3!}(\Theta w)^3 - \frac{675}{2 \cdot 5!}(\Theta w)^5 - \frac{1148175}{2\cdot 9!} (\Theta w)^9 \right. \\
&\qquad\qquad\qquad\qquad \left. - \frac{1148175}{11!}(\Theta w)^{11} + \frac{27807650325}{4\cdot 13!}(\Theta w)^{13} + O(w^{15})\right) \\
h(w) &= (1-w)^6\left((\Theta w)^2 + \frac{881798400}{12!} (\Theta w)^{12} + O(w^{22})\right)
\end{align*}
which satisfy
\[ g^2-fh - 3h^2 = 0 \]
to the precision computed.  This gives
\begin{equation} \label{eqn:exmgen0x}
\begin{aligned}
x(w) = \frac{h(w)}{g(w)} &= (\Theta w) - \frac{9}{3!}(\Theta w)^3 + \frac{1215}{2\cdot 5!}(\Theta w)^5 - \frac{-59535}{7!}(\Theta w)^7 \\
&\qquad\quad + \frac{12170655}{9!}(\Theta w)^9 - \frac{-6708786525}{2\cdot 11!}(\Theta w)^{11} + O(w^{13}).
\end{aligned}
\end{equation}

From (\ref{eqn:psit}) we obtain
\[ \frac{w}{\kappa}=\frac{\psi(t)}{\kappa} = t^{1/5} \left(1 + \frac{1}{10} t + \frac{3943}{89100}t^2 + \frac{2161}{81000}t^3 + \frac{23027911}{1246590000}t^4 + \dots \right)  \] 
so
\[ \phi(w)=\psi^{-1}(t) = \left(\frac{w}{\kappa}\right)^5 - \frac{1}{2}\left(\frac{w}{\kappa}\right)^{10} + \frac{637}{3564}\left(\frac{w}{\kappa}\right)^{15} - \frac{383}{7128}\left(\frac{w}{\kappa}\right)^{20}+O(w^{25}) \]
so with $\Theta$ as in (\ref{eqn:exmgen0Theta}) we have
\[ \phi(w)=\psi^{-1}(t) = \frac{81}{2}(\Theta w)^5 - \frac{6561}{8}(\Theta w)^{10} + \frac{4179357}{352} (\Theta w)^{15} +O(w^{20}). \]
From the ramification data, we have $\phi(x)$ is a rational function in $x$ of degree $5$ with numerator $x^5$; comparing power series, we then compute
\begin{equation} \label{eqn:phix0}
\phi(x) = \frac{648x^5}{324x^5+405x^4-120x^2+16} = \frac{648x^5}{(3x+2)^3(12x^2-9x+2)} 
\end{equation}
and
\begin{equation} \label{eqn:phix1}
\phi(x)-1 = \frac{(3x-2)^3(12x^2+9x+2)}{(3x+2)^3(12x^2-9x+2)}. 
\end{equation}
We also verify this numerically: we find that
\[ \{x(\alpha_i \cdot z_c)\}_i = \{0.666\ldots, -0.375000\ldots \pm 0.16137\ldots \sqrt{-1}\} \]
and
\[ \{x(\alpha_i \cdot z_b)\}_i = \{-0.666\ldots, 0.375000\ldots \pm 0.16137\ldots \sqrt{-1}\} \]
which agree with the roots of the numerator and denominator in (\ref{eqn:phix1}).

The expression for $\phi(x)$ in (\ref{eqn:phix0}) is not the simplest one possible, allowing for linear fractional transformations in the $x$ and $\phi$: for example, we have
\[ \frac{1}{\phi(1/3x)} = 6x^5 - 5x^3 + 15/8 x + 1/2 \]
and further substituting $x \leftarrow x-1/2$ we obtain the nicer form
\begin{equation} \label{eqn:phiop}
\phi_{\text{op}}(x) = 6x^5 - 15x^4 + 10x^3 = x^3(6x^2 - 15x + 10) = 1 + (x-1)^3(6x^2 + 3x + 1).
\end{equation}

We can verify that the computed map is correct by computing it using a direct method: beginning with
\[ x^3(a_2x^2 + a_1x+a_0) = 1 + (x-1)^3(b_2 x^2 + b_1x + b_0) \]
we find the system of equations
\[ b_0 -1 = 3b_0 - b_1 = 3b_0-3b_1+b_2 = a_0-b_0+3b_1-3b_2 = a_1-b_1+3b_2=a_2-b_2 = 0\]
which has the unique solution $(a_0,a_1,a_2,b_0,b_1,b_2)=(10,-15,6,1,3,6)$ as in (\ref{eqn:phiop}).  (For more on the direct method and its extensions, see Sijsling--Voight \cite{SijslingVoight}.)

The other triple, companion to \eqref{eqn:sigmaexmgenus0}, is
\begin{equation} \label{eqn:othertriple}
\sigma_0=(1\ 3\ 2\ 4\ 5), \quad \sigma_1 = (1\ 2\ 3), \quad \sigma_\infty=(1\ 5\ 4).
\end{equation}
If we run the above procedure again on this triple, we find the same expression for $x(w)$ as in (\ref{eqn:exmgen0x}) except now with $\Theta=0.3917053\ldots - 1.205545\ldots i$, i.e., the complex conjugate of the previous value (\ref{eqn:exmgen0Theta}).  In other words, the antiholomorphic map given by complex conjugation identifies these two \Belyi\ maps.

The federalist approach in this example does not give anything new, since the origin $0 \in \calD$ does not move $\{\alpha_i \cdot 0 : i=0,\dots,5\}=\{0\}$ under the coset representatives $\alpha_i$.  Using our implementation, the whole calculation takes a few seconds with a standard desktop CPU.  (We will see the benefits of using improved numerical linear algebra techniques in the next example.)

There are seven other refined passports of genus zero with $d=5$, all with monodromy group $G=S_5$.  The tables below list representative permutation triples for each refined passport, as well as polynomials $\phi_0$, $\phi_1$, $\phi_\infty$ such that
\begin{align*}
\phi(x) = \frac{\phi_0(x)}{\phi_\infty(x)} = 1 + \frac{\phi_1(x)}{\phi_\infty(x)} \, .
\end{align*}

\begin{landscape}
$$
\begin{array}{c|c|c|c}
a,b,c & i & \sigma_i & \phi_i\\
\hline
5 & 0 & (1\ 2\ 3\ 4\ 5) & 128x^5\\
2 & 1 & (1\ 2) & (3x-2)^2 (14x^3 + 17 x^2 + 12 x + 4)\\
4 & \infty & (1\ 5\ 4\ 3) & (x+2)^4(2x-1)
\end{array}
$$

$$
\begin{array}{c|c|c|c}
a,b,c, & i & \sigma_i & \phi_i\\
\hline
5 & 0 & (1\ 2\ 3\ 4\ 5) & 108 x^5\\
2 & 1 & (1\ 3) & (x+1)^2 (36 x^3 - 12 x^2 - 2 x + 1)\\
6 & \infty & (1\ 5\ 4)(2\ 3) & (2x - 1)^3 (3x + 1)^2
\end{array}
$$

$$
\begin{array}{c|c|c|c}
a,b,c & i & \sigma_i & \phi_i\\
\hline
4 & 0 & (1\ 2\ 3\ 4) & 84375 x^4 (15x + 16)\\
2 & 1 & (1\ 3)(4\ 5) & (3x^2 + 512x +512)^2 (17x - 16)\\
6 & \infty & (1\ 5\ 4)(2\ 3) & (2x - 1)^3 (3x + 1)^2
\end{array}
$$

$$
\begin{array}{c|c|c|c}
a,b,c & i & \sigma_i & \phi_i\\
\hline
6 & 0 & (1\ 2\ 3)(4\ 5) & -6250x^3 (5x + 1)^2\\
2 & 1 & (1\ 5)(3\ 4) & -(113x^2 + 23x + 2)^2 (13x - 1)\\
6 & \infty & (1\ 4)(2\ 5\ 3) & (19x + 2)^2 (3x - 1)^3
\end{array}
$$

$$
\begin{array}{c|c|c|c}
a,b,c & i & \sigma_i & \phi_i\\
\hline
4 & 0 & (1\ 2\ 3\ 4) & -84375 x^4 (5x - 28)\\
3 & 1 & (1\ 4\ 5) & -(11x - 28)^3 (321 x^2 + 784 x + 784)\\
4 & \infty & (2\ 5\ 4\ 3) & 1792 (x+7)^4 (3x - 4)
\end{array}
$$

$$
\begin{array}{c|c|c|c|c}
a,b,c & i & \multicolumn{2}{c|}{\sigma_i} & \phi_i\\
\hline
4 & 0 & (1\ 2\ 3\ 4) & (1\ 2\ 3\ 4) & 6(26496 + 10719\sqrt{6}) x^4 ((8 - 3\sqrt{6})x -1)) \\
3 & 1 & (2\ 4\ 5) & (3\ 4\ 5) & -6(2244584 + 99851\sqrt{6})((830 - 339\sqrt{6})x^2 - x - 1)((3-\sqrt{6})x - 1)^3\\
6 & \infty & (1\ 2)(3\ 5\ 4) & (1\ 3\ 2)(4\ 5) & (10x - (7 + 3\sqrt{6}))^2 (6x + (12 + 5\sqrt{6}))^3
\end{array}
$$

$$
\begin{array}{c|c|c|c}
a,b,c & i & \sigma_i & \phi_i\\
\hline
3 & 0 & (1\ 2\ 3) & 2x^3 (9x^2 - 5)  \\
6 & 1 & (1\ 2\ 5)(3\ 4) & (x - 1)^3 (3x + 2)^2\\
6 & \infty & (1\ 4\ 3)(2\ 5) & (x + 1)^3 (3x - 2)^2
\end{array}
$$
\end{landscape}
Note in the example where $(a,b,c) = (4,3,6)$ that the size of the refined passport is 2 and correspondingly the field of definition of $\phi$ is the quadratic extension $\Q(\sqrt{6})$.
\end{exm}

\begin{exm} \label{exm:genus0nonewt}
We now consider a more complicated example.  Consider the triple
\begin{equation}\label{eqn:dessin-12,2,5}
\sigma_0=(1\ 2\ 3\ 4)(5\ 6\ 7), \quad \sigma_1=(2\ 3)(4\ 5)(6\ 7), \quad \sigma_\infty=(1\ 5\ 6\ 4\ 2)
\end{equation}
which generates $S_7$.  The size of the refined passport for this triple is $2$; the other triple (up to uniform conjugation) is
\[ \sigma_0=(1\ 2\ 3\ 4)(5\ 6\ 7), \quad \sigma_1=(1\ 2)(4\ 5)(6\ 7), \quad \sigma_\infty=(1\ 5\ 6\ 4\ 3). \]

We find $g(X(\Gamma))=0$ and $\dim_\C S_4(\Gamma) =2$ so we take forms in weight $k=4$.  We work in precision $\varepsilon=10^{-100}$.  The federalist approach has an advantage here: expanding around the vertices $\alpha_i\cdot 0$ for $i=1,5$, we can take  the radius of expansion to be $\rho=0.821032\ldots$ and degree $N=1168$; if we did not, we would need to take radius $0.966411$ and $N=6740$.  Our implementation uses 71 Arnoldi iterations, each requiring about $2$ seconds, to find a basis element, and the total calculation time was about 6 minutes; computing a full numerical kernel instead took over 15 minutes.  

We find with 
\begin{equation} 
\Theta = 1.3259366\ldots + 0.63528305\ldots \sqrt{-1} = \sqrt[12]{\frac{1632000-961536\sqrt{-5}}{4117715}}\left(\frac{1}{\kappa}\right) \end{equation}
that
\[ x = (\Theta w)^3 + \frac{-537128+59872\sqrt{-5}}{194145\cdot 9!}(\Theta w)^9 + \frac{37661184-1520640\sqrt{-5}}{3895843\cdot 12!}(\Theta w)^{12} + O(w^{15}) \]
and
\[ \phi(x) = \frac{p(x)}{q(x)} = 1+\frac{r(x)}{q(x)} \]
where
\begin{align*} 
p(x) &= (379-2102\sqrt{-5})x^4\left((4-\sqrt{-5})x+(-5+\sqrt{-5})\right)^3 \\
q(x) &= 3\left((-2-\sqrt{-5})x+(2+2\sqrt{-5})\right)^5\left((-57+20\sqrt{-5})x^2+(25-19\sqrt{-5})x+60\right) \\
r(x) &= -2(3x-(1+\sqrt{-5}))\left((259+5\sqrt{-5})x^3 + (-356+43\sqrt{-5})x^2\right. \\
&\qquad\qquad\qquad\qquad\qquad\qquad \left. - (144+144\sqrt{-5})x + (240+96\sqrt{-5})\right)^2.
\end{align*}
See below (Example \ref{exm:genus0newt}) for the use of Newton's method on this example.

Finally, we repeat this calculation with the other triple (\ref{eqn:othertriple}) and we find that all coefficients are (to numerical precision) the complex conjugates of the first triples, thus confirming that the Galois action is as expected on this refined passport.
\end{exm}

\begin{figure}[h]

\includegraphics[scale=1]{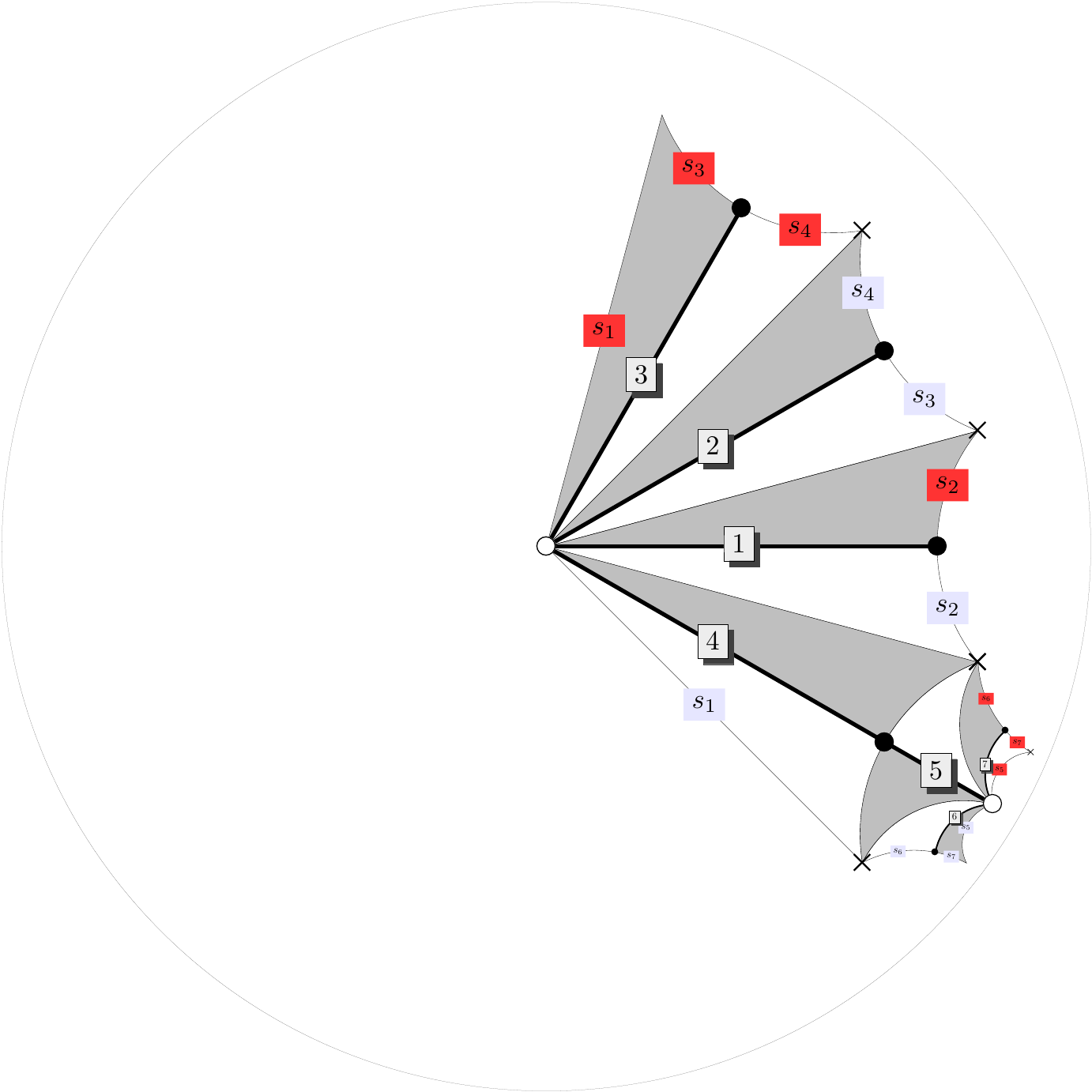}

\begin{center}
\begin{tabular}{ll}
\toprule
Label & Coset Representative\\
\midrule
$1$ & 1 \\
$2$ & $\delta_a^{}$ \\
$3$ & $\delta_a^{2}$ \\
$4$ & $\delta_a^{-1}$ \\
$5$ & $\delta_a^{-1}\delta_b^{}$ \\
$6$ & $\delta_a^{-1}\delta_b^{}\delta_a^{}$ \\
$7$ & $\delta_a^{-1}\delta_b^{}\delta_a^{-1}$ \\
\bottomrule
\end{tabular}
\hfill
\begin{tabular}{ll}
\toprule
Label & Side Pairing Element\\
\midrule
$s_{1}$ & $\delta_a^{-4}$ \\
$s_{2}$ & $\delta_b^{}$ \\
$s_{3}$ & $\delta_a^{}\delta_b^{}\delta_a^{-2}$ \\
$s_{4}$ & $\delta_a^{}\delta_b^{-1}\delta_a^{-2}$ \\
$s_{5}$ & $\delta_a^{-1}\delta_b^{}\delta_a^{3}\delta_b^{-1}\delta_a^{}$ \\
$s_{6}$ & $\delta_a^{-1}\delta_b^{}\delta_a^{}\delta_b^{}\delta_a^{}\delta_b^{-1}\delta_a^{}$ \\
$s_{7}$ & $\delta_a^{-1}\delta_b^{}\delta_a^{}\delta_b^{-1}\delta_a^{}\delta_b^{-1}\delta_a^{}$ \\
\bottomrule
\end{tabular}
\end{center}
\bigskip
\caption{A fundamental domain and conformally drawn dessin for $\Gamma$ corresponding to the permutation triple \ref{eqn:dessin-12,2,5} using the methods in section \ref{sec:cosets}.}  \label{fig:dessin-12,2,5}
\end{figure}
\FloatBarrier

\subsection*{Genus one}

Now suppose that $X(\Gamma)$ has genus $g=1$.  Then the space $S_2(\Gamma)$ of modular forms of weight $2$ has dimension $\dim S_2(\Gamma)=1$ and so is spanned by a form 
\[ f(z)=(1-w)^2\sum_{n=0}^{\infty} b_n w^n = (1-w)^2 \sum_{n=0}^{\infty} \frac{c_n}{n!} (\Theta w)^n \]
where as usual $w=w(z)$ and $\Theta$ is chosen as in (\ref{eqn:Theta0}) to (experimentally) obtain coefficients $c_n$ that are algebraic integers in a number field.  Then $f(z)\,dz$ is the unique (nonzero) holomorphic differential $1$-form on the Riemann surface $X(\Gamma)$ up to scaling by $\C^\times$.  

To compute an equation for $X(\Gamma)$, we compute the Abel-Jacobi map, computing its lattice of periods, and then use the Weierstrass uniformization to obtain a model, as follows.  (For a reference, see Silverman \cite[Chapter VI]{Silverman} or Cremona \cite[Chapters II--III]{Cremona}.)  Let $\{\gamma_i\}_i$ be the side pairing elements arising from the coset graph and fundamental domain $D_\Gamma$ for $\Gamma$ as computed in Algorithm \ref{alg:cosetalg}.  Suppose that $\gamma_i$ pairs the vertex $v_i$ with $v_i'=\gamma_i v_i$ in $D_\Gamma$.  In particular, these vertices are within the domain of convergence of $D_\Gamma$.  Then we compute the set of  normalized periods by
\begin{equation} \label{eqn:varpi}
\varpi_i = \int_{z(v_i)}^{z(v_i')} \Theta f(z)\,dz = \int_{v_i}^{v_i'} f(w) \frac{d(\Theta w)}{(1-w)^2} \approx \sum_{n=0}^N \frac{c_n}{(n+1)!} (\Theta w)^{n+1} \bigg|_{v_i}^{v_i'}.
\end{equation}
The span 
\[ \Lambda = \sum_i \Z \varpi_i = \Z \omega_1 + \Z\omega_2 \subset \C \]
is a (full) lattice in $\C$, and the Abel-Jacobi map
\begin{equation} \label{eqn:AJ}
\begin{aligned}
X(\Gamma)=\Gamma \backslash \calH &\xrightarrow{\sim} \C/\Lambda \\
z &\mapsto u(z) = \int_{0}^{w(z)} f(w) \frac{d(\Theta w)}{(1-w)^2} 
\end{aligned}
\end{equation}
is an isomorphism of Riemann surfaces.  We may assume without loss of generality, interchanging $\omega_1,\omega_2$ if necessary, that $\tau=\omega_1/\omega_2 \in \calH$.  (Further, by application of the Euclidean algorithm, i.e., fundamental domain reduction for $\SL_2(\Z)$, we may assume $|\repart \tau\,| \leq 1/2$ and $|\tau| \geq 1$, so that $\impart \tau \geq \sqrt{3}/2$.) Let $q=\exp(2\pi i\tau)$.  Then $X(\Gamma)$ has $j$-invariant $j(\Lambda)=j(\tau)$
\[ j(q) = \frac{1}{q} + 744 + 196884q + 21493760q^2 + \dots, \]
where $j(q)$ is the modular elliptic $j$-invariant (not to be confused with the automorphy factor by the same name!), and we have a uniformization
\begin{equation} \label{eqn:AJtoE}
\begin{aligned}
\C/\Lambda &\xrightarrow{\sim} E(\Gamma) : y^2 = x^3 - 27c_4(\Lambda)x - 54c_6(\Lambda) \\
u &\mapsto \left(\wp(\Lambda;u),\frac{\wp'(\Lambda;u)}{2}\right) 
\end{aligned}
\end{equation}
where
\[ \wp(\Lambda;u) = \frac{1}{u^2} + \sum_{k=2}^{\infty} (2k-1)\frac{2\zeta(2k)}{\omega_2^{2k}}E_{2k}(\tau) u^{2k-2} \]
is the Weierstrass $\wp$-function, 
\[ E_{2k}(\tau) = 1+(-1)^k\frac{4k}{B_{2k}}\sum_{n=0}^{\infty} \sigma_{2k-1}(n) q^n \]
is the normalized Eisenstein series with $B_{2k}$ the Bernoulli numbers, so that
\begin{equation} \label{eqn:c4c6}
\begin{aligned}
E_4(\tau) &= 1 + 240\sum_{n=1}^{\infty} \sigma_3(n)q^n= 1 + 240q + 2160q^2 + 6720q^3 + \dots \\
E_6(\tau) &= 1 - 504\sum_{n=1}^{\infty} \sigma_5(n)q^n = 1 - 504q - 16632q^2 - 122976q^3 + \dots,
\end{aligned}
\end{equation}
and finally
\[ c_{2k}(\Lambda) = \left(\frac{2\pi}{6\omega_2}\right)^{2k} E_{2k}(\tau). \]
Implicitly, we are taking the point $w=0$ to be the origin of the group law, so we always get an elliptic curve (instead of a genus one curve).

The \Belyi\ map is computed from the explicit power series $\phi(w) \in \Q[[(w/\kappa)^a]]$ as in (\ref{eqn:expandpuiseux}).  A general function on $E(\Gamma)$ is of the form
\[ \phi(x,y) = \frac{p(x) + q(x) y}{r(x)} \]
Since $x(w)=\wp(\Lambda;u(w))$ and $y(w)=\wp'(\Lambda;u(w))/2$ are also computed power series in $w$, we can solve for the coefficients of $p(x),q(x),r(x)$ using linear algebra.

\begin{exm}
The smallest degree $d$ for which there exists a hyperbolic refined passport of genus one is $d=4$, and there is a unique such refined passport with $(a,b,c)=(4,3,4)$ and representative triple
\[ \sigma_0 = (1\ 2\ 3\ 4), \quad \sigma_1=(1\ 2\ 3), \quad \sigma_\infty = (1\ 4\ 2\ 3). \]

\begin{figure}[h] 
\includegraphics{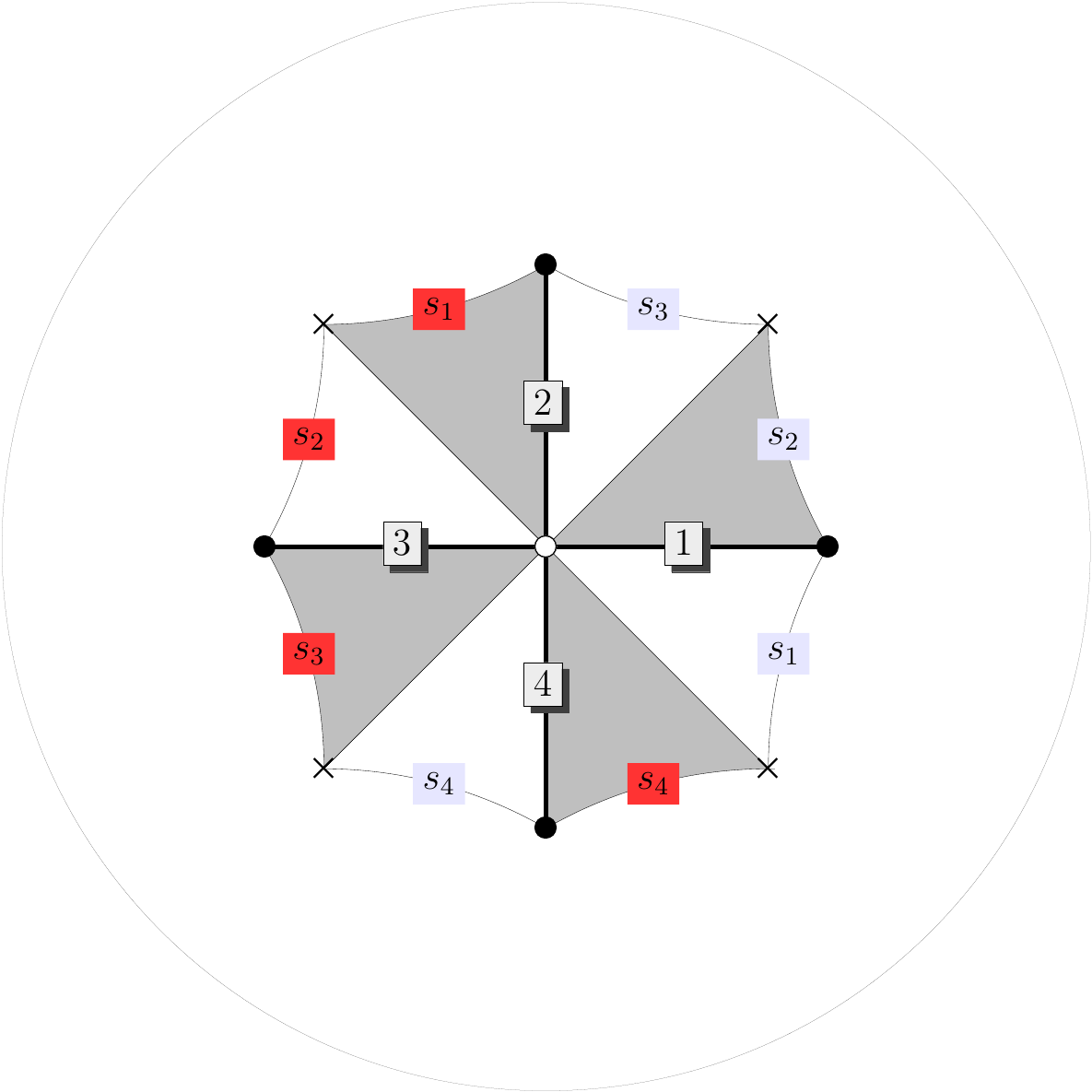}

\begin{center}
\begin{tabular}{ll}
\toprule
Label & Coset Representative\\
\midrule
$1$ & 1 \\
$2$ & $\delta_a^{}$ \\
$3$ & $\delta_a^{2}$ \\
$4$ & $\delta_a^{-1}$ \\
\bottomrule
\end{tabular}
\;\;\;\;\;\;\;\;\;\;
\begin{tabular}{ll}
\toprule
Label & Side Pairing Element\\
\midrule
$s_{1}$ & $\delta_b^{}\delta_a^{-1}$ \\
$s_{2}$ & $\delta_b^{-1}\delta_a^{-2}$ \\
$s_{3}$ & $\delta_a^{}\delta_b^{}\delta_a^{-2}$ \\
$s_{4}$ & $\delta_a^{-1}\delta_b^{}\delta_a^{}$ \\
\bottomrule
\end{tabular}
\end{center}
\bigskip
\caption{A fundamental domain and conformally drawn dessin for $\Gamma$ corresponding to the permutation triple
$$
\sigma_0 = (1\;2\;3\;4),\quad\sigma_1=(1\;2\;3),\quad\sigma_\infty = (1\;4\;2\;3)
$$
 using the methods in \ref{sec:cosets}.} \label{fig:dessin-4,3,4}
\end{figure}

In under a second, we compute in precision $\varepsilon=10^{-30}$ the normalized differential
\[ f(z) = 1 + (\Theta w) - \frac{1}{2!}(\Theta w)^2 + \frac{6}{4!}(\Theta w)^4 - \frac{6}{5!}(\Theta w)^5 + \frac{126}{6!}(\Theta w)^6 + O(w^8) \]
with
\[ \Theta = -1.786853\ldots i = \sqrt[4]{\frac{27}{4}} i \left(\frac{1}{\kappa}\right). \]

We find two classes which span the homology of $X(\Gamma)$, as in Figure \ref{fig:homology}.

\begin{figure}[h]
\includegraphics{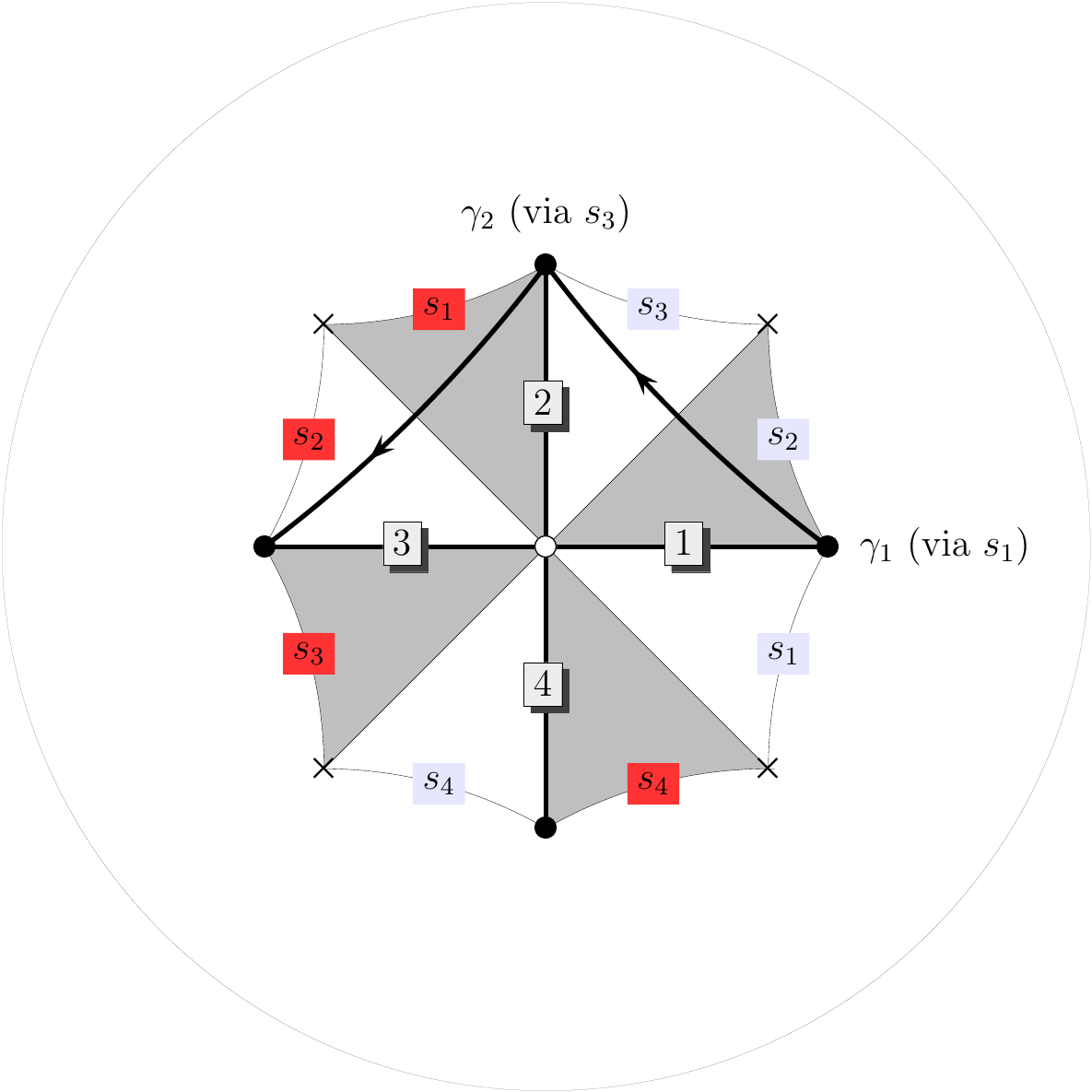}

\caption{Representatives $\gamma_1, \gamma_2$ of the homology classes of $X(\Gamma)$.}  \label{fig:homology}
\end{figure}

We integrate as in (\ref{eqn:varpi}): we find the two periods
\[ \omega_1,\omega_2=\pm 1.68575035481259604287120365776 - 1.07825782374982161771933749941 i \]
and numerically we verify that all other periods lie in the lattice spanned by these two periods.  We find
\[ \tau=\omega_1/\omega_2 = -0.4193\ldots + 0.9078\ldots\sqrt{-1} \]
and $\tau$ belongs to the fundamental domain for $\SL_2(\Z)$.  We then compute
\[ j(q)=31.648529\ldots=\frac{207646}{6561} \]
and from (\ref{eqn:c4c6}) that
\[ c_4=-0.03626\ldots = -\frac{47}{1296}, \qquad c_6=-0.05056\ldots = -\frac{2359}{46656} \]
so $X(\Gamma)$ has equation
\begin{equation} \label{eqn:y2x3}
y^2 = x^3 + \frac{47}{48} x + \frac{2359}{864}
\end{equation}
with minimal model
\begin{equation} \label{genus1X}
y^2 = x^3 + x^2 + 16x + 180, 
\end{equation}
(via the map $(x,y) \mapsto (4x-1/3, 8y)$), the elliptic curve with conductor $48$ with Cremona label \textsf{48a6}.  

We compute the \Belyi\ map by consideration of power series.  From (\ref{eqn:AJ}) we have
\[ u(w) = (\Theta w) + \frac{1}{2}(\Theta w)^2 - \frac{1}{6}(\Theta w)^3 + O(w^5) \]
and from (\ref{eqn:AJtoE}), with
\[ \wp(\Lambda;u) = \frac{1}{u^2} - \frac{47}{240}u^2 - \frac{337}{864} u^4 + \frac{2209}{172800} u^6 + O(u^8) \]
we compute that 
\begin{align*}
 x(w) &= \frac{1}{(\Theta w)^2} - \frac{1}{\Theta w} + \frac{13}{12} - (\Theta w) + \frac{3}{5} (\Theta w)^2 + O(w^3) \\
 y(w) &= -\frac{1}{(\Theta w)^3} + \frac{3}{2} \frac{1}{(\Theta w)^2} - 2\frac{1}{(\Theta w)} + \frac{9}{4} - \frac{12}{5}(\Theta w) + O(w^2) 
\end{align*}
and confirm again that $x(w),y(w)$ satisfy \eqref{eqn:y2x3} to the precision computed.  Then from (\ref{eqn:expandpuiseux}) we find
\[ \phi(w) = \frac{27}{4}(\Theta w)^4 - \frac{108}{5}(\Theta w)^8 + \frac{8181}{175}(\Theta w)^{12} + O(w^{16}). \]
Then we compute a linear relation among the functions 
\[ \phi,\phi x, \phi y, \phi x^2, 1\]
and we find the \Belyi\ map
\[ \phi(x,y) = \frac{139968x^2 - 23328x + 490860 + 279936y}{(12x-13)^4} \]
of degree $4$ with divisor $4(13/12,9/4)-4\infty$ and $\phi-1$ with divisor 
\[ \ddiv(\phi-1) = 3(-5/12,-3/2)+(67/12,27/2)-4\infty \]
as desired.
\end{exm}

\subsection*{Hyperelliptic}

Now suppose that $X(\Gamma)$ has genus $g=g(X(\Gamma)) \geq 2$.  Then we define the \emph{canonical map}
\begin{equation} \label{eqn:canhyp}
\begin{aligned}
X(\Gamma) &\xrightarrow{\sim} \PP^{g-1} \\
z &\mapsto (f_1(w) : \dots : f_g(w))
\end{aligned}
\end{equation}
where $f_1,\dots,f_g$ is a basis for $S_2(\Gamma)$.  We have two possibilities.

We suppose first, in this subsection, that $X(\Gamma)$ is \defi{hyperelliptic}, meaning that there exists a (nonconstant) degree $2$ map $X(\Gamma) \to \PP^1$.  We test if $X$ is hyperelliptic as follows: we compute the function $x(w)=g(w)/h(w)$ where $g,h \in S_2(\Gamma)$ are as in (\ref{eqn:gandh}), compute $y=(dx/dw)/h(w)$, and test if there is an equation of the form 
\[ y^2 + u(x)y = v(x) \]
with $\deg u(x) \leq g+1$ and $\deg v(x) \leq 2g+2$; if so, we have found our embedding as a hyperelliptic curve.  One can use Riemann--Roch to verify that this test is correct (see e.g.\ \cite[\S 4]{Galbraith}).  Analogous to the case in genus $g=1$, a general function on $X$ is of the form
\[ \phi(x,y) = \frac{p(x)+q(x) y}{r(x)} \]
and linear algebra yields an expression for $\phi$.  

\begin{exm}
Consider the triple
\[ \sigma_0=\sigma_1=(1\ 2\ 3 \ 4\ 5), \quad \sigma_{\infty}=(1\ 4\ 2\ 5\ 3) \]
generating $\Z/5\Z$.  Then $g=2$ so $X$ is hyperelliptic.  We compute a basis for $S_2(\Gamma)$ to precision $10^{-30}$ in $15$ seconds.  With
\[ \Theta = \sqrt[5]{\frac{1}{24}} \left(\frac{1}{\kappa}\right) \]
we find the following forms in $S_2(\Gamma)$:
\[ \begin{aligned} 
\frac{f(w)}{(1-w)^2} &= 1 + \frac{120}{5!}(\Theta w)^5 + \frac{21}{10!}(\Theta w)^{10} - \frac{45864}{15!}(\Theta w)^{15} + \frac{237414996}{20!}(\Theta w)^{20} + O(w^{25}) \\
\frac{g(w)}{(1-w)^2} &= (\Theta w) - \frac{3}{6!}(\Theta w)^6 - \frac{126}{11!}(\Theta w)^{11} - \frac{366912}{16!}(\Theta w)^{16} 
+O(w^{21})
\end{aligned}
\]
We then compute
\[ x(w) = \frac{f(w)}{g(w)} = (\Theta w)^{-1} + \frac{3}{10\cdot 4!}(\Theta w)^4 + \frac{1218}{55 \cdot 9!}(\Theta w)^9 + O(w^{14}) \]
and $y(w)=x'(w)/g(w)$ where $x'$ denotes the derivative with respect to $\Theta w$.  Then $x$ has a pole of order $1$ at $w=0$ and $y$ has a pole of order $3$, therefore we obtain the equation 
\[ y^2 = x^6 - \frac{1}{6} x. \]
We then compute using (\ref{eqn:expandpuiseux}) the expansion
\[ \phi(w) = \frac{1}{24}(\Theta w)^5 - \frac{1}{1152}(\Theta w)^{10} + \frac{23}{1824768}(\Theta w)^{15} + O(w^{20}) \]
and hence
\[ \phi(x,y) = \frac{y+1}{2x^3}. \]
This dessin corresponds to a Galois \Belyi\ map and so it can be constructed using other means; however, it makes for a nice example and our method did not use this property.
\end{exm}

\subsection*{General case}

Now suppose still that $g \geq 2$ but now that $X(\Gamma)$ is not hyperelliptic; then the canonical map (\ref{eqn:canhyp}) is a closed embedding.  If $g=3$ then $f_1,f_2,f_3$ satisfy an equation of degree $4$; otherwise, the ideal of equations vanishing on the image is generated in degree $2$ and $3$: in fact, in degree $2$ as long as $X$ is neither trigonal nor isomorphic to a plane quintic.  So again using linear algebra and the power series expansions for $f_1,\dots,f_g \in S_2(\Gamma)$, we can compute a generating set for the ideal of quadrics and cubics vanishing on the image to obtain equations for $X(\Gamma)$.  

Once equations for $X(\Gamma)$ are obtained, we again write down using Riemann-Roch and the explicit representation of the function field of $X(\Gamma)$ a general function of degree $d$ and use linear algebra using explicit power series to compute an equation for the \Belyi\ map $\phi$; alternatively, we compute the ramification values and use Riemann-Roch to find a function with specified divisor.

\begin{exm}
Consider the rigid permutation triple
\[
\sigma_0 = (1\ 2\ 3\ 4\ 5\ 6\ 7), \quad \sigma_1=(1\ 6\ 2\ 5\ 7\ 3\ 4), \quad \sigma_{\infty}=(1\ 5\ 3\ 6\ 2\ 4\ 7), \]
a $(7,7,7)$-triple generating $\PSL_2(\F_7)\cong \GL_3(\F_2)$ and with refined passport of genus $g=3$.  With
\[ \Theta = 0.9280540\ldots + 0.4469272\ldots\sqrt{-1} = \sqrt[7]{-16}\left(\frac{1}{\kappa}\right) \]
we find a basis of $3$ forms $f,g,h \in S_2(\Gamma)$ with coefficients in the field $\Q(\nu)$ where $\nu=(-1+\sqrt{-7})/2$.  The canonical embedding of a genus three curve is a quartic in $\PP^2$, which we compute in the echelonized basis to be
\begin{align*}
&11f^3h - 11f^2g^2 + (5\nu + 2)f^2h^2 + (-10\nu + 29)fg^2h + (-25\nu - 10)fh^3 \\
&\qquad + (5\nu + 2)g^4 + (-8\nu - 12)g^2h^2 + (-21\nu + 18)h^4 = 0
\end{align*}
valid to the precision computed.

\begin{figure}[h] 
\includegraphics{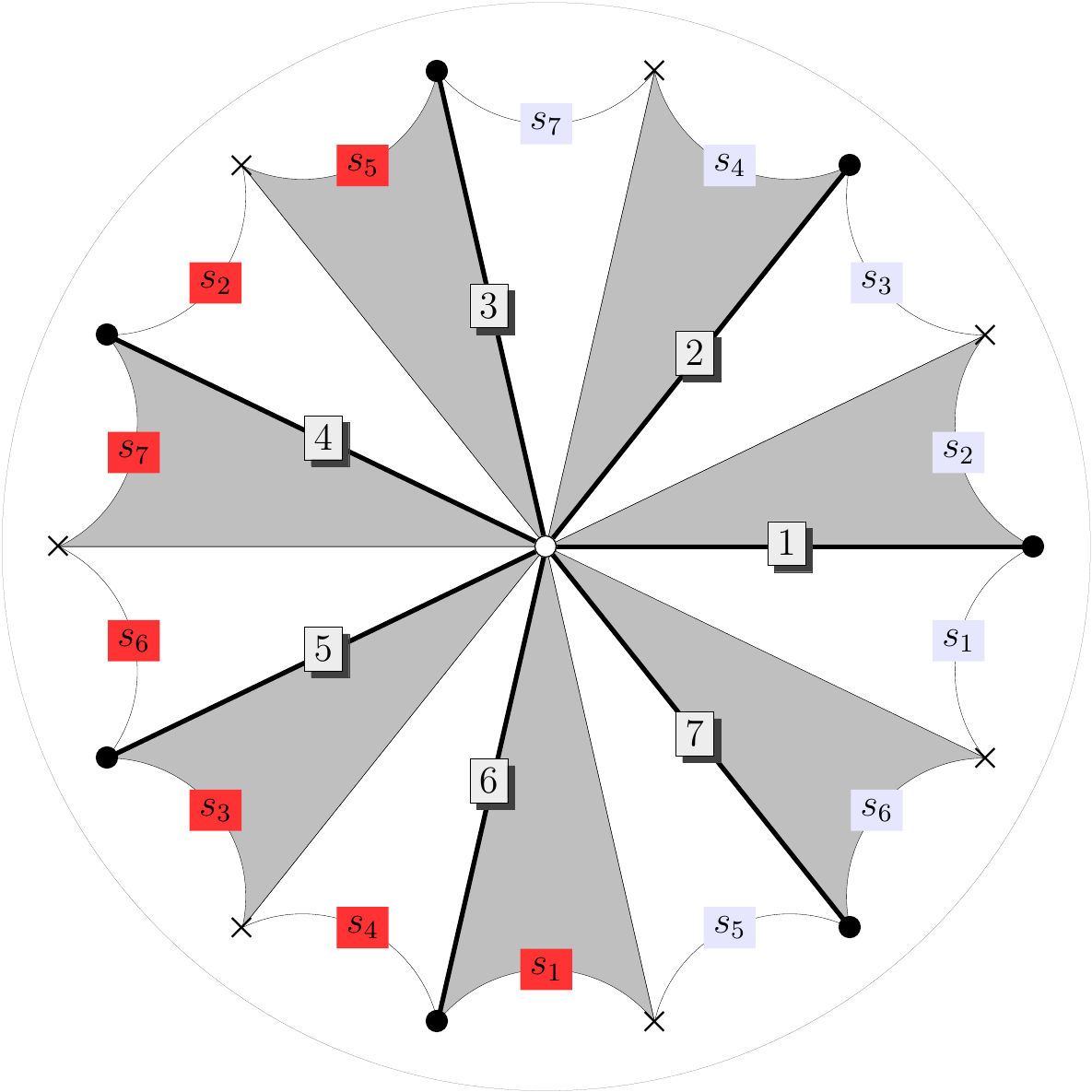}

\begin{center}
\begin{tabular}{ll}
\toprule
Label & Coset Representative\\
\midrule
$1$ & 1 \\
$2$ & $\delta_a^{}$ \\
$3$ & $\delta_a^{2}$ \\
$4$ & $\delta_a^{3}$ \\
$5$ & $\delta_a^{-3}$ \\
$6$ & $\delta_a^{-2}$ \\
$7$ & $\delta_a^{-1}$ \\
\bottomrule
\end{tabular}
\;\;\;\;\;\;\;\;\;\;
\begin{tabular}{ll}
\toprule
Label & Side Pairing Element\\
\midrule
$s_{1}$ & $\delta_b^{}\delta_a^{2}$ \\
$s_{2}$ & $\delta_b^{-1}\delta_a^{-3}$ \\
$s_{3}$ & $\delta_a^{}\delta_b^{}\delta_a^{3}$ \\
$s_{4}$ & $\delta_a^{}\delta_b^{-1}\delta_a^{2}$ \\
$s_{5}$ & $\delta_a^{-1}\delta_b^{}\delta_a^{-2}$ \\
$s_{6}$ & $\delta_a^{-1}\delta_b^{-1}\delta_a^{3}$ \\
$s_{7}$ & $\delta_a^{2}\delta_b^{}\delta_a^{-3}$ \\
\bottomrule
\end{tabular}
\end{center}
\bigskip
\caption{A fundamental domain and conformally drawn dessin for $\Gamma$ corresponding to the permutation triple
$$
\sigma_0 = (1\;2\;3\;4),\quad\sigma_1=(1\;2\;3),\quad\sigma_\infty = (1\;4\;2\;3)
$$
 using the methods in Section \ref{sec:cosets}.} \label{fig:dessin-7,7,7}
\end{figure}

We find by evaluating the normalized power series the values $\phi^{-1}(0)=\{(1:0:0)\}$, $\phi^{-1}(1)=\{(0:-\nu:1)\}$, and $\phi^{-1}(\infty)=\{(0:\nu:1)\}$.  Computing with Riemann-Roch (linear algebra), we find a generator for the $1$-dimensional space of functions with divisor $7(1:0:0)-7(0:\nu:1)$, containing $\phi$; computing similarly with $7(0:-\nu:1)-7(0:\nu:1)$, containing $\phi-1$, we find that
\[ \phi(x,y) = \frac{p_2(y)x^2 + p_1(y)x + p_0(y)}{44(y-\nu)^7} \]
where
\begin{align*}
p_2(y) &= (-55\nu - 154)y^4 + (-77\nu - 110)y^3 + (187\nu - 638)y^2 + (385\nu + 550)y \\
p_1(y) &= (55\nu + 154)y^6 + (-44\nu + 440)y^4 + (202\nu + 380)y^3 \\
&\qquad + (-715\nu + 990)y^2 + (-570\nu - 668)y
\end{align*}
and
\begin{align*}
p_0(y) &= (77\nu + 110)y^7 + (-308\nu - 264)y^6 + (-1049\nu - 1590)y^5 \\
&\qquad + (2112\nu - 3520)y^4 + (3299\nu + 3634)y^3 \\
&\qquad + (1980\nu + 4840)y^2 + (-2535\nu + 1846)y - 616\nu - 880
\end{align*}
where $x=f/h$ and $y=g/h$.  
\end{exm}

\subsection*{Refinement using Newton's method}

We have shown how to obtain equations for the \Belyi\ map $\phi:X(\Gamma) \to X(\Delta)$ defined over a number field.  As the complexity of the examples increases, it is easier to use Newton's method to refine the solution: we set up a system of equations with rational coefficients that define the \Belyi\ map and then refine the approximate solution.  

\begin{exm} \label{exm:genus0newt}
We return to Example \ref{exm:genus0nonewt}.  Computing initially in precision $\varepsilon=10^{-30}$, we find an approximation to $\phi(x)$
\[ \frac{(8.265 - 2.345\sqrt{-1})x^4(x-1.190 - 0.1064\sqrt{-1})^3}{(x-1.555- 0.4969\sqrt{-1})^5(x^2 - (0.6334-0.2483\sqrt{-1})x-0.6515-0.5111\sqrt{-1})} \]
in about $10$ seconds.  We then set up the equations that describe rational functions with the same ramification pattern \cite{SijslingVoight}: we have
\[ \phi(x) = \frac{a_9x^4(x+a_1)^3}{(x+a_6)^5(x^2+a_8x+a_7)} = 1+(a_9-1)\frac{(x+a_5)(x^3+a_4x^2+a_3x+a_2)^2}{(x+a_6)^5(x^2+a_8x+a_7)}. \]
Expanding and setting the coefficients equal to zero, we obtain a set of $7$ equations in $9$ unknowns.  By scaling $x$ we may assume $a_2=a_3$.  The condition that $x$ is normalized as in (\ref{eqn:xw}) implies the further condition that
\[ a_1a_6a_8 + 5a_1a_7 - 3a_6a_7 = 0, \]
and this now gives in total $9$ equations in $9$ unknowns.  We apply Newton's method to the approximate solution above: the solution is correct to error $10^{-20}$, so Newton iteration converges after $20$ iterations in about $4$ seconds to a solution which is correct to $10^{-500}$.  Now the coefficients of $\phi$ are very easy to recognize!
\end{exm}

\begin{exm} \label{exm:PSU3(5)}
We now consider a much larger example, for which the combined efficiency gains introduced in this paper are essential to make the calculation practical.  We consider the exceptional permutation representation of $G=\PSU_3(\F_5) \hookrightarrow S_{50}$ arising from the action on the cosets of the subgroup $A_7 \leq G$ of index $50$.  This exceptional representation can be seen on the Hoffman-Singleton graph  \cite{Hafner}, a 7-regular undirected graph with 50 vertices and 175 edges with automorphism group equal to the semidirect product $\PSU_3(\F_5):2$ (extension by Frobenius).

We look for rigid refined passports and refined passports of genus zero.  There are $6$ rigid refined passports, all with orders $(5,2,7)$: there are $2$ of genus zero and $4$ of genus two.  There are $5$ refined passports of genus zero: the $2$ rigid $(5,2,7)$-refined passports, along with $2$ $(5,2,8)$-triples with refined passport of size 2, and one $(4,2,10)$-triple with refined passport of size 6.  

We first consider the rigid $(5,2,7)$-refined passports.  They arise from the two conjugacy classes of order $7$ in $G$, interchanged by an automorphism of $G$ (also an automorphism of the associated dessin), so as above it suffices to consider just one of these triples:

\begin{equation}\label{eqn:sigmasPSU3(5)}
\begin{aligned}
\sigma_0 &= (2\ 44\ 11\ 7\ 28)(3\ 42\ 4\ 46\ 17)(5\ 32\ 6\ 34\ 16)(8\ 12\ 43\ 33\ 25) \\
&\qquad (10\ 38\ 37\ 31\ 50) (13\ 49\ 40\ 22\ 23)(15\ 48\ 24\ 29\ 26) \\
&\qquad (18\ 21\ 27\ 30\ 45)(19\ 20\ 39\ 41\ 36) \\
\sigma_1 &= (1\ 5)(2\ 25)(4\ 20)(6\ 33)(7\ 47)(8\ 17)(9\ 42)(10\ 41)(11\ 35) \\
&\qquad (13\ 30)(14\ 18)(16\ 37)(19\ 40)(22\ 27)(24\ 29)(26\ 46) \\
&\qquad (28\ 39)(31\ 49)(34\ 48)(38\ 44) \\
\sigma_{\infty} &=
            (1\ 5\ 37\ 44\ 25\ 6\ 32)(2\ 39\ 4\ 9\ 42\ 3\ 8)(7\ 35\ 11\ 38\ 41\ 28\ 47) \\
            &\qquad (10\ 50\ 49\ 30\ 22\ 19\ 36)(12\ 17\ 26\ 24\ 34\ 33\ 43) \\
&\qquad (13\ 23\ 27\ 21\ 14\ 18\ 45)(15\ 46\ 20\ 40\ 31\ 16\ 48) \\
\end{aligned}
\end{equation}

In Figure \ref{fig:dessin-PSU3(5)}, we display the fundamental domain, using the methods in Section \ref{sec:cosets}.  The picture is in the unit disc (centered at $i$) and with the left most point of the fundamental domain at the center. The labels are omitted because of the small triangles close to the boundary.

\begin{figure}[h] 
\includegraphics{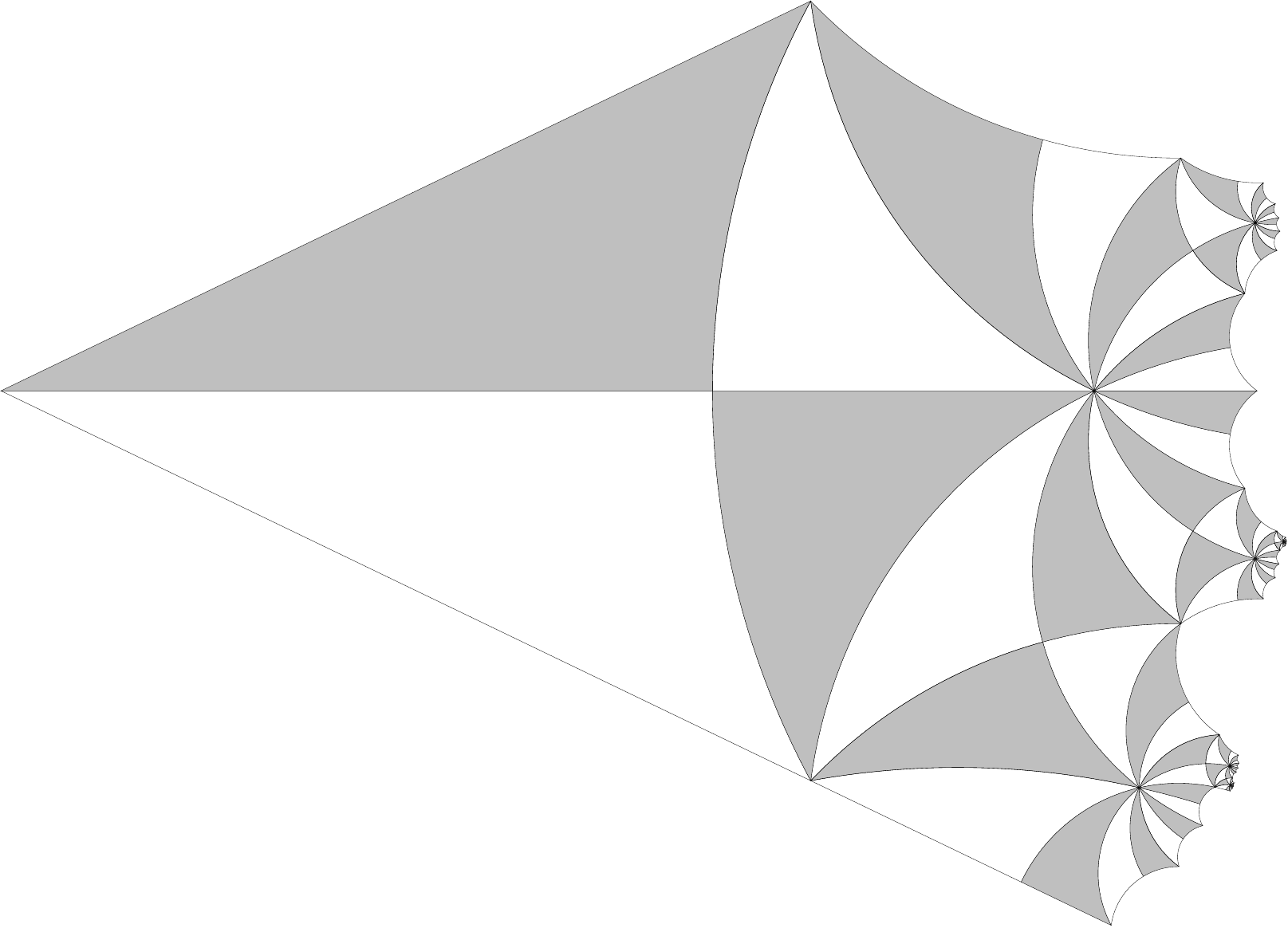}

\caption{Fundamental domain for $\Gamma$ corresponding to the permutation triple \ref{eqn:sigmasPSU3(5)}.} \label{fig:dessin-PSU3(5)}
\end{figure}

We work in precision $\varepsilon=10^{-60}$, and it takes about $20$ minutes (with twice $172$ Arnoldi iterations) for each triple followed by another $15$ minutes of Newton iteration to obtain $1000$ digits of precision.  Initially, we find a \Belyi\ map defined over $\Q(\zeta_7)^+$, the cubic totally real subfield of the cyclotomic field $\Q(\zeta_7)$; this is also the field of definition of the normalized power series we computed.  After some creative simplification, involving a change of coordinates as in (\ref{eqn:phiop}), we descend to the \Belyi\ map 
\[ \phi(x) = \frac{p(x)}{q(x)} = 1+ \frac{r(x)}{q(x)} \]
where
\begin{align*}
p(x) &= 2^6 (x^4 + 11x^3 - 29x^2 + 11x + 1)^5 (64x^5 - 100x^4 + 150x^3 - 25x^2 + 5x + 1)^5 \\
&\qquad \cdot (196x^5 - 430x^4 + 485x^3 - 235x^2 + 30x + 4) \\
q(x) &= 5^{10} x^7 (x+1)^7 (2x^2 - 3x + 2)^7 (8x^3 - 32x^2 + 10x + 1)^7 \\
\end{align*}
and
\begin{align*}
r(x) &= (28672x^{20} - 2114560x^{19} + 13722240x^{18} - 65614080x^{17} + 245351840x^{16} \\
&\qquad\qquad
        - 660267008x^{15} + 1248458280x^{14} - 1700835920x^{13} + 1704958640x^{12} \\
&\qquad\qquad -
        1267574420x^{11} + 690436992x^{10} - 257110380x^9 + 52736995x^8 - 
        948040x^7 \\
&\qquad\qquad - 1171555x^6 - 246148x^5 + 86660x^4 + 11060x^3 - 1520x^2  - 240x - 8)^2 \\
&\qquad\cdot (16384x^{10} - 34960x^9 - 160960x^8 + 620820x^7 - 792960x^6 + 416087x^5 \\
&\qquad\qquad        + 57435x^4 + 935x^3 + 705x^2 + 110x + 4).
\end{align*}

Remarkably, the Galois group of the degree $20$ polynomial $r(x)$ is a group of order $240$, an extension of the normal subgroup $A_5$ by the Klein 4-group $V_4$, and the associated number field has discriminant $2^{14} 5^{26} 7^{15}$.  Similarly, the other factors have non-generic Galois group with small ramification.  

Computing monodromy (as in the next subsection), we verify that the polynomial $f(x) = p(x) - tq(x) \in \C[x;t]$ indeed has Galois group $\PSU_3(\F_5)$ and this descends to a $\PSU_3(\F_5)$-extension of the field of rationality $K=\Q(\sqrt{-7})$.  As a Galois extension of $\Q(t)$, therefore, we find the group $\PSU_3(\F_5):2$; the additional descent to $\Q$ follows from the fact that the dessin has an automorphism of order $2$ identifying it with its complex conjugate, which allows us to descend the cover defined by $\phi$ to $\Q$ as a \defi{mere cover} (see D\`ebes and Emsalem \cite{DebesEmsalem} for the definition), but the full Galois group is only defined over $K$.  The existence of the rational function $f(x)$ was assured by the general theory of rigidity \cite{MalleMatzat}, but the explicit polynomial itself is new.  

Independently of this, we verify that the Galois group $G$ over $K=\Q(\sqrt{-7})$ is indeed $\PSU_3(\F_5)$ following a suggestion of Kl\"uners.  First, we prove that $G$ is primitive.  The Galois group of a specialization is a subgroup of $G$.  Specializing at $t=2$ and reducing modulo a prime of norm $29$, we find an element $\sigma$ of order $7$ with cycle type $7^7 1$.  Then $\sigma \in G$ has a unique fixed point and so fixes any block that contains it.  Then the cardinality $s$ of this block has $s \mid 50$ and $s \equiv 1 \pmod{7}$, so $s=1,50$, hence $G$ is primitive.  (See more generally Kl\"uners \cite[\S 3.3]{KluenersSubfields}.)  Alternatively, we check that the field extension $K(f)/K(t)$ has no proper subfields.  Now there are $9$ primitive permutation groups of degree $50$.  We eliminate the smaller groups $H$ by finding a conjugacy class that belongs to $\PSU_3(\F_5)$ that does not belong to $G$; alternately, using the algorithm of Fieker--Kl\"uners \cite{FiekerKlueners}, we compute that the Galois group of the specialization at $t=2$ is $\PSU_3(\F_5)$.  Using the classification of primitive subgroups of $S_{50}$, this leaves $4$ possibilities: $\PSU_3(\F_5)$, $\PSU_3(\F_5):2$, $A_{50}$, or $S_{50}$.  The discriminant of $f(t)$ is 
\[ \disc(f(t)) = \frac{5^{560} 7^{1092}}{2^{1918}} t^{36}(t-1)^{20}, \]
which is a square, so immediately $G \leq A_{50}$.  To rule out $A_{50}$ (and again $S_{50}$ at the same time), we show that $G \leq S_{50}$ is not $2$-transitive: the polynomial of degree $50\cdot 49$ with roots $x_i-x_j$ with $i \neq j$, where $x_i$ are the roots of $f(x)$---computed symbolically using resolvents---factors over $\Q(t)$.  Thus $G \leq \PSU_3(\F_5):2$.  (Since the prime $p=53 \nmid \#S_{50}$ is a prime of good reduction for the cover by a theorem of Beckmann \cite{Beckmann}, it is enough to compute this factorization over $\F_p$.)  Computing another relative resolvent to distinguish between $\PSU_3(\F_5)$ and $\PSU_3(\F_5):2$, we verify that $G = \PSU_3(\F_5)$; alternatively, we construct the fixed field of $\PSU_3(\F_5)$ in the Galois extension over $\Q$ and find that it is equal to $K$.

To conclude with the remaining examples at the other extreme, for the refined passport of size $6$ with orders $(4,2,10)$, we repeat these steps with all $6$ triples, and in about $3$ CPU hours we have computed computed the result to $1000$ digits of precision.  Computing the minimal polynomial of each coefficient, we find a \Belyi\ map defined over the number field $L$ with defining polynomial
\[ x^{12} - x^{11} - 6x^{10} - 10x^9 - 15x^8 - 16x^7 - 19x^6 - 24x^5 - 35x^4 - 30x^3 - 6x^2 - 9x - 9. \]
The field $L$ is a quadratic extension of the field $K=\Q(\sqrt{5},\sqrt[3]{10})$ and $L$ (like $K$) is ramified only at $2,3,5$.  We were not able to find a descent of this map to $K$, and we guess that there might be an obstruction of some kind, for example, the curve $X(\Gamma)$ may be a conic defined over $K$ that is not $K$-isomorphic to $\PP^1$.  For more on the issue of descent, see Sijsling--Voight \cite[\S 7]{SijslingVoight} and
the references therein.
\end{exm}

\subsection*{Verification}

Once a \Belyi\ map has been computed, it remains to verify that it is in fact correct.  There are a number of methods to achieve this: see the survey by Sijsling--Voight \cite{SijslingVoight}.  

In our situation, it is enough to verify that the cover computed is a three-point cover; this can be checked by a discriminant calculation.  Once this has been checked, the way that the cover was constructed guarantees that it has the correct monodromy, computed as they were to sufficient precision to separate triangles.  This verifies that the cover is correct over $\C(t)$ (and thus over $\Qbar(t)$ and identifies the embedding of the number field yielding the specific monodromy triple).

Example \ref{exm:PSU3(5)} was also independently verified over $\C(t)$ using code due to Bartholdi, Buff, Graf von Bothmer, and Kr\"oker \cite{Bartholdi}.

\begin{rmk} \label{rmk:makerigorous}
At the present time, our method is not rigorous, and we have not analyzed its running time.  However, our computations are quite successful in practice, and we believe that it may be possible to furnish this rigorous analysis.  A bound on the height of a \Belyi\ map, in the context of Arakelov geometry, has been established by Javanpeykar \cite{Javanpeykar}, and so it would suffice to establish that the linear system yielding power series expansions of modular forms is well-conditioned.  In this vein, it would also be helpful to establish the integrality of coefficients of the power series expansions under a suitable choice of $\Theta$, as indicated in Remark \ref{rmk:Shimura}.  We leave these investigations for future work.
\end{rmk}

\end{document}